%% file: main.tex
\documentclass{article}
\pdfoutput=1
\usepackage{xparse}
\usepackage[utf8]{inputenc}
\usepackage{lmodern}
\usepackage{textgreek}
\usepackage{csquotes}
\usepackage{subfiles}
\usepackage{csquotes}
\usepackage{xr-hyper}
\usepackage[backend=biber,style=alphabetic, doi = false, isbn=false]{biblatex}
\renewbibmacro{in:}{}
\AtEveryBibitem{%
  \clearfield{primaryClass}%
  \clearfield{EDITOR}
   \clearfield{eid}
}
\usepackage{amsmath}
\usepackage{amsthm}
\usepackage{amssymb}
\usepackage{mathtools}
\usepackage{amsfonts}
\usepackage{wasysym}
\usepackage{bbm}
\usepackage{thmtools} %new math ambient - also for cross references
\usepackage{faktor}
\let\originallhook\lhook
\let\originalrhook\rhook
\let\amslrcorner\lrcorner
\let\amsllcorner\llcorner
\let\amsulcorner\ulcorner
\let\amsurcorner\urcorner

\usepackage{MnSymbol}
\let\lrcorner\amslrcorner
\let\llcorner\amsllcorner
\let\ulcorner\amsulcorner
\let\urcorner\amsurcorner
\newcommand{\lhook}{\mathrel{\raise.018ex\hbox{$\originallhook$}}}
\newcommand{\rhook}{\mathrel{\raise.018ex\hbox{$\originalrhook$}}}
\usepackage{quiver}
\usepackage{float}
\usepackage{xargs}
\usepackage{enumitem}
\usepackage[hidelinks]{hyperref} %with hidden color
\usepackage[capitalize, nameinlink]{cleveref}

\usepackage{graphicx}
\usepackage{xcolor}

\usepackage{pdfpages}

\usepackage[margin=3.4cm]{geometry}
 
\definecolor{Bnavy}{RGB}{0, 66, 128}
\definecolor{Bdust}{RGB}{140,179,217}
\definecolor{Bsugarpaper}{RGB}{198, 217, 236}
\definecolor{Bgreen}{RGB}{142, 183, 114}
\definecolor{Blimegreen}{RGB}{202, 222, 189}
\definecolor{Bgreentheme}{RGB}{36, 87, 1}

\usepackage{makecell}
\usepackage{stackengine} 

\usepackage{microtype}

\usepackage{tikz} 
\usepackage{tikz-cd}
\usetikzlibrary{ %snake arrows
	calc,%
	arrows,%
	shapes,
	shapes.geometric,
    positioning,
    fit
} 
\usepackage{tikz-3dplot}
\usetikzlibrary{backgrounds, calc, math, patterns, arrows.meta}
\usetikzlibrary{decorations.markings}

\tikzcdset{arrow style=tikz, diagrams={>=stealth}}
\usepackage[labelsep=quad,indention=10pt]{subfig}
\theoremstyle{plain}
\newtheorem{theorem}{Theorem}[subsection]

\newtheorem{lemma}[theorem]{Lemma}
\newtheorem{corollary}[theorem]{Corollary}
\newtheorem{proposition}[theorem]{Proposition}
\newtheorem{construction}[theorem]{Construction}
\newtheorem{innercustomthm}{Theorem}
\newenvironment{ctheorem}[1]
  {\renewcommand\theinnercustomthm{#1}\innercustomthm}
  {\endinnercustomthm}

\theoremstyle{definition}
\newtheorem{definition}[theorem]{Definition}
\newtheorem{example}[theorem]{Example}
\newtheorem{recollection}[theorem]{Recollection}
\crefname{recollection}{Recollection}{Recollections}
\newtheorem{remark}[theorem]{Remark}
\newtheorem{notation}[theorem]{Notation}

\newtheorem{conjecture}[theorem]{Conjecture}
\newtheorem{question}[theorem]{Question}
\newcounter{diagram}  % Create counter for diagrams
\crefname{diagram}{Diagram}{Diagrams}
\creflabelformat{diagram}{#2\textup{(#1)}#3}
\newenvironment{diagram}[1][]{%
    \crefalias{equation}{diagram}
    \begin{equation}%
    \begin{tikzcd}[#1]%
}{%
    \end{tikzcd}%
    \end{equation}%
}

\usepackage{notation}
    \title{Presenting the topological stratified homotopy hypothesis}
\author{Lukas Waas}
\date{January 2024}
    \newlist{RhoProps}{enumerate}{4}
    \setlist[RhoProps]{label*=(\roman*)${}_\rho$, leftmargin = 1.3cm}
    \Crefname{RhoPropsi}{Property}{Properties}
    \crefname{RhoPropsi}{Property}{Properties}
     \newlist{phiProps}{enumerate}{4}
    \setlist[phiProps]{label*=(\roman*)${}_\phi$, leftmargin = 1.3cm}
    \Crefname{phiPropsi}{Property}{Properties}
     \newlist{PhiProps}{enumerate}{4}
    \setlist[PhiProps]{label*=(\roman*)${}_\Phi$, leftmargin = 1.3cm}
    \Crefname{PhiPropsi}{Property}{Properties}
    \newlist{fProps}{enumerate}{4}
    \setlist[fProps]{label*=(\roman*)${}_{f''}$, leftmargin = 1.3cm}
    \Crefname{fPropsi}{Property}{Properties}
    \newlist{fibrancyProp}{enumerate}{4}
    \setlist[fibrancyProp]{label*=(\roman*), leftmargin = 1.3cm}
    \Crefname{fibrancyPropi}{Property}{Properties}
    \newlist{QuinProp}{enumerate}{4}
    \setlist[QuinProp]{label*=(\roman*), leftmargin = 1.3cm}
    \Crefname{QuinProp}{Property}{Properties}
    \newlist{TransfLem}{enumerate}{4}
    \setlist[TransfLem]{label*=(\roman*), leftmargin = 1.3cm}
    \Crefname{TransfLemi}{Assumption}{Assumptions}
    \newlist{IntroQuestions}{enumerate}{4}
    \setlist[IntroQuestions]{label*=Q(\arabic*), leftmargin = 1.3cm}
    \Crefname{IntroQuestionsi}{Question}{Questions}
\begin{document}
\maketitle
\begin{abstract}
 This article is concerned with three different homotopy theories of stratified spaces: The one defined by Douteau and Henriques, the one defined by Haine, and the one defined by Nand-Lal.
 One of the central questions concerning these theories has been how precisely they connect with geometric and topological examples of stratified spaces, such as piecewise linear pseudomanifolds, Whitney stratified spaces, or more recently Ayala, Francis and Tanaka's conically smooth stratified spaces. More precisely, so far, it has been an open question whether there exist (semi-)model structures on stratified topological spaces that present these theories, in which such relevant examples of stratified spaces are bifibrant. Here, we prove an affirmative answer to this question. As a consequence, we obtain a model categorical interpretation of a stratified homotopy hypothesis. Specifically, we show that Lurie's stratified singular simplicial set functor induces a Quillen equivalence between the semimodel category of stratified topological spaces presenting Nand-Lal's homotopy theory of stratified spaces and the left Bousfield localization of the Joyal model structure that corresponds to such $\infty$-categories in which every endomorphism is an isomorphism. We then perform a detailed investigation of bifibrant objects in these model structures of stratified spaces, proving a series of detection criteria and illuminating the relationship to Quinn's homotopically stratified spaces.
\end{abstract}

\subfile{1Introduction}

\subfile{2Preliminaries}

\subfile{3StratifiedHomotopyTheories}

\subfile{4Modelstructures}

\subfile{5BifibrantObjects}
\section*{Acknowledgments}
The author is supported by a PhD-stipend of the Landesgraduiertenförderung Baden-Württemberg. Furthermore, he wants to thank Sylvain Douteau and Peter Haine for the insightful discussions concerning their work.
\printbibliography
\appendix
\subfile{6Appendix}
\end{document}

%% file: 1Introduction.tex
\section{Introduction}
Conceptually speaking, Grothendieck's homotopy hypothesis refers to the following statement:
\begin{displayquote}
Assigning to a topological space its $\infty$-groupoid of paths induces an equivalence between the \textit{homotopy theory of spaces} (more precisely, CW-complexes) and the \textit{homotopy theory of $\infty$-groupoids.}
\end{displayquote}
Whether this statement is regarded as a theorem or a conjecture, of course, strongly relies on the precise model of $\infty$-categories - and consequently of $\infty$-groupoids - one has in mind. Nowadays, the following result due to Kan and Quillen is often taken as a formal interpretation of the homotopy hypothesis:
\begin{ctheorem}{HH}\cite{Quillen}\label{PSHH:Thm_Kan_Quillen_Equ} The geometric realization and singular simplicial set adjunction 
\[
\real{-} \colon  \sSetN \rightleftharpoons \TopN : \Sing
\]
induces a Quillen equivalence between the Quillen model structure on topological spaces and the Kan-Quillen model structure on simplicial sets.
\end{ctheorem}
That this is a formal interpretation of the homotopy hypothesis can be argued as follows.
\begin{enumerate}
    \item The Quillen equivalence induces an equivalence 
    \[
    \TopN[W^{-1}] \simeq \textnormal{\textbf{Kan}}[H_k^{-1}]
    \]
    between topological spaces localized at weak homotopy equivalences and Kan complexes localized at homotopy equivalences of Kan complexes.
    \item Quasi-categories (as introduced by Joyal and popularized by Lurie in \cite{HigherTopos}) have proven to be a powerful and versatile model for the theory of $(\infty,1)$-categories ($\infty$-categories henceforth).
    Kan complexes are the $\infty$-groupoid within the framework of quasi-categories, and it follows from the existence of the Kan-Quillen and the Joyal model structure that the right-hand side thus defines the homotopy theory of $\infty$-groupoids $\iGrpd$. The singular simplicial set then provides a model for the $\infty$-groupoid of paths in this interpretation of $\infty$-categories.
    This justifies the usage of Kan complexes as a model for $\infty$-groupoids.
    \item It follows from the existence of the Quillen model structure on $\TopN$, or more classically an argument involving Whitehead's theorem, that the inclusion of CW complexes into topological spaces, $\textnormal{\textbf{CW}} \hookrightarrow \TopN$ induces an equivalence $\textnormal{\textbf{CW}}[H^{-1}] \simeq \TopN [W^{-1}]$, where $H$ is the class of homotopy equivalences.
    Since most spaces of geometric interest at least have the homotopy type of a CW complex (see, for example, \cite{Milnor1959OnSH}),
    it follows that if one is interested in studying the homotopy types of such classical spaces, one may perform such an investigation in $\TopN[W^{-1}]$. Thus, the latter can rightfully be called the homotopy theory of spaces.\\
\end{enumerate}
Combining these insights, one obtains an equivalence 
\[
\iSpaces = \textnormal{\textbf{CW}}[H^{-1}] \simeq \TopN[W^{-1}] \simeq \textnormal{\textbf{Kan}}[H_k^{-1}]  = \iGrpd,
\]
as asserted in the homotopy hypothesis.
\\
One of the central assertions of \cite{AFRStratifiedHomotopyHypothesis} is a smooth stratified version of the homotopy hypothesis.
Recall that, roughly speaking and in the broadest sense, a stratified space is a topological space together with a decomposition into disjoint pieces, the so-called strata. A stratified map is, again roughly speaking, a continuous map between such objects that has the property that the image of each stratum in the source is completely contained in a stratum in the target. Stratifications of topological spaces often arise naturally when investigating spaces with singularities, by decomposing a singular space into manifold pieces (see, for example, \cite{Whitney,mather1970notes,mather1973strat,thom1969ensembles}) and in these scenarios the set of strata tends to naturally inherit the structure of a poset from the topological closure relation.
The homotopy theory of such stratified spaces, using homotopies that also preserve the strata, was first investigated in detail by Quinn in \cite{quinn1988homotopically}. Quinn focused on a class of stratified spaces with excellent homotopical properties, the so-called homotopically stratified spaces (called homotopically stratified sets in \cite{quinn1988homotopically}), proving, among other results, a stratified version of the s-cobordism theorem for certain homotopically stratified spaces with manifold strata.
Following a more differential topological framework, which they established in \cite{LocalStructOnStrat}, in \cite{AFRStratifiedHomotopyHypothesis} Ayala, Francis, and Rozenblyum developed a homotopy theory of stratified spaces within the context of differential topology. 
One of their central assertions was the following statement: 
\begin{displayquote}
    Let $\StratN_{C^\infty}$ denote the category of conically smooth stratified spaces (with conically smooth stratified maps), and $H_s$ the class of stratified homotopy equivalences (in the conically smooth sense). Then there is a fully faithful embedding
\[
\StratN_{C_\infty}[H_s^{-1}] \hookrightarrow \iCat
\]
into the homotopy theory of $\infty$-categories\footnote{At the current point in time, it appears that \cite{AFRStratifiedHomotopyHypothesis} is missing a definition of what exactly the category of conically smooth stratified spaces they consider is. Note that one cannot simply use the category of conically smooth stratified spaces of \cite{LocalStructOnStrat}, since the latter does not contain the stratified simplices. Without a definition of the stratified smooth category accommodating such spaces with corners (and the development of the necessary theory), it is currently not possible to verify the truth of the smooth stratified homotopy hypothesis.}.
\end{displayquote}
 Furthermore, they conjectured a topological analogue of this result. This relies on the exit path construction of McPherson, Treuman, Woolf, and Lurie (\cite{Woolf,treumann_2009,HigherAlgebra}). Roughly speaking, this construction associates to a stratified space an $\infty$-category in which morphisms are given by paths that either remain within one stratum, or start in one stratum and immediately exit into a higher one.
\begin{conjecture}\cite{AFRStratifiedHomotopyHypothesis}\label{con:stratified_homotopy_hypothesis}
     Topological exit paths define a fully faithful functor
          \[  \textnormal{\textbf{Exit}} \colon \mathcal{S}\textnormal{\textbf{trat}} \hookrightarrow \iCat \]
        from a homotopy theory of topological stratified spaces $\mathcal{S}\textnormal{\textbf{trat}}$ into $\infty$-categories.
\end{conjecture}
Again, just as in the case of the classical homotopy hypothesis, any answer to this conjecture must first provide a formal interpretation of the homotopy theories in question. In this case, the difficulty lies with the left-hand side, specifically the question after a homotopy theory of topological stratified spaces.
In recent years, three of such theories have been independently proposed (not all with the intent to tackle the topological stratified homotopy hypothesis).
In the following $\StratN$ will always denote the category of all poset-stratified spaces, with stratified maps between them (see, for example, \cite{douteauwaas2021} for an overview). Similarly, we denote by $\TopPN$, for a partially ordered set $\pos$, the category of poset-stratified spaces over a fixed poset, with morphisms given by stratified maps that descend to the identity on $\pos$, so-called \define{stratum-preserving} maps.
One way of defining a homotopy theory of stratified spaces is, of course, to specify a category of topologically stratified spaces, in these cases, the category of all poset-stratified spaces with stratified maps, and then to localize the latter at an appropriate class of \textit{stratified weak equivalences}. 
\begin{enumerate}
    \item In \cite{douteauEnTop,Henriques} Douteau and  Henriques independently built on a result of Miller's in \cite{miller2013}, concerning Quinn's homotopically stratified spaces. In \cite{miller2013}, it was shown that stratum-preserving homotopy equivalences between homotopically stratified spaces are precisely such maps that induce homotopy equivalence on the strata and on the so-called \textit{$[p<q]$-homotopy links} (where $p<q$ are elements of the stratifying poset). Recall, that these are given by the spaces of exit-paths, starting in the $p$-stratum and immediately exiting into the $q$-stratum. Since for general poset-stratified spaces pairwise homotopy links are insufficient to even guarantee that the underlying map is a weak equivalence of topological spaces, \cite{douteauEnTop,Henriques} additionally considered so-called generalized homotopy links, obtained by replacing the stratified interval by a stratified simplex. Then a weak equivalence of stratified spaces is defined to be a stratified map that induces an isomorphism on stratifying posets and weak equivalences on all generalized homotopy links (which include the strata).
    We call such stratified maps \define{poset-preserving diagrammatic equivalences} and denote the resulting $\infty$-category obtained by localizing poset-preserving diagrammatic equivalences by $\AltStratD$. Analogously, given a fixed poset $\pos$, we denote the localization of $\TopPN$ at stratum-preserving (poset-preserving) diagrammatic equivalences over $\pos$ by $\AltTopPD$.
    Then taking generalized homotopy links induces an equivalence of $\infty$-categories
        \[
        \HolIP[] \colon \AltTopPD \simeq \FunC (\sd(\pos)^\op, \iSpaces)
        \]
    where the right hand side denotes the category of $\iSpaces$ valued presheaves on the subdivision on the poset $\pos$.
    Explicitly, $\sd(\pos)$ is the category whose objects are finite increasing sequences $[p_0 < \cdots < p_n]$ in $\pos$ (so-called regular flags) and whose morphisms are inclusions of subsequences.  
    \item In \cite{haine2018homotopy}, Haine followed the insight that the categories of exit paths associated to a classical stratified space (for example, conically stratified space or homotopically stratified space) always come with a conservative functor into the stratification poset. This is simply due to the fact that restricted to each stratum the construction produces the classical $\infty$-groupoid of paths. Furthermore, the fact that exit-paths define an $\infty$-category at all implies that, for such stratified spaces, the natural maps from the generalized homotopy links into the homotopy pullbacks of pairwise homotopy links
        \[
        \HolIP[{[p_0< \dots < p_n]}] (\tstr) \to \HolIP[p_0 < p_1](\tstr) \times^{H}_{X_{p_1}} \dots \times^{H}_{X_{p_{n-1}}} \HolIP[p_{n-1} < p_n](\tstr) 
        \]
    are weak homotopy equivalences. Hence, only certain diagrams $D \in \FunC (\sd(\pos)^\op, \iSpaces)$ --- those were the corresponding morphism is an isomorphism in the $\infty$-category of spaces --- may arise as the homotopy-link diagram of classical examples of stratified spaces. Diagrams that satisfy this property are called \define{d\'ecollages} (see \cite{haine2018homotopy,Exodromy}). The full $\infty$-subcategory of d\'ecollages, $\Dec$, turns out to be a (left) localization of $\FunC (\sd(\pos)^\op, \iSpaces)$ (i.e. a reflective subcategory). 
    \cite{haine2018homotopy} then defines his class of weak equivalences of stratified spaces over a fixed poset $\pos$ as the class of such stratified maps that map into isomorphisms under the composition
           \[
           \TopPN \xrightarrow{\HolIP[]}  \FunC (\sd(\pos)^\op, \iSpaces) \to \Dec.
           \]
    We call such weak equivalences (stratum-preserving) categorical equivalences, and denote the resulting homotopy theory by $\AltTopPC$. This makes $\AltTopPC$ a left localization of Douteau and Henriques' theory $\AltTopPD$, that is canonically equivalent to $\Dec$. One may then extend this class of weak equivalences to the case of varying poset $\StratN$, by requiring a weak equivalence to induce isomorphisms on stratifying posets and a stratum-preserving categorical equivalence after identifying the posets along the isomorphism. We call such stratified maps \define{poset-preserving categorical equivalences} (for reasons which will become apparent later) and denote the resulting homotopy theory by $\AltStratC$. 
    \item Finally, in \cite{nand2019simplicial}, Nand-Lal took the approach of transferring the weak equivalences along Lurie's functor of stratified singular simplices, which provides a model for the $\infty$-category of exit paths associated to a stratified space\footnote{Strictly speaking, Nand Lal works with a full subcategory of $\StratN$, but from a homotopy theoretic perspective this difference turns out to be irrelevant (see \cref{prop:eq_of_sur_strat})}:
        \[
    \SingS \colon \StratN \to \sSetN.
        \]
    More precisely, a stratified map $\tstr \to \tstr[Y]$ is called a \define{categorical equivalence} if the induced simplicial map $\SingS (\tstr) \to \SingS(\tstr[Y])$ is a Joyal equivalence (also called categorical equivalences; see \cite{HigherTopos}). One major difference compared to the previous two theories is that weak equivalences in this setting are not required to preserve the stratification poset. In this sense, the stratification posets are only essential to a stratified homotopy type insofar as they determine what paths are exit paths (what simplices are exit-simplices). We denote the induced homotopy theory by $\AltStratCR$.
\end{enumerate}
The first answer to \cref{con:stratified_homotopy_hypothesis} was given in \cite{haine2018homotopy}, where Haine used a result of Douteau (\cite[Thm. 3]{douteauEnTop}) to prove an equivalence of $\infty$-categories
    \[
    \AltTopPC
    \simeq \AbStrP,
    \]
where $\AbStrP$ is the $\infty$-category of conservative functors of $\infty$-categories with target $\pos$, so-called abstract stratified homotopy types. Thus, the equivalence provides an answer to \cref{con:stratified_homotopy_hypothesis} for the setting of a fixed poset.
The equivalence is constructed by passing through the equivalence of $\AltTopPC$ with $\Dec$, and then, in turn, constructing an equivalence of the latter with a model for $\AbStrP$ given in the language of complete Segal spaces. 
In this sense, the equivalence is not constructed directly in terms of Lurie's $\SingS$ construction (from a $1$-categorical perspective at least). A priori, it is unclear whether the weak equivalences defining $\AltTopPC$ can be characterized through $\SingS$. \\
\\
Here, we aim to provide an answer to the topological stratified homotopy hypothesis that is similarly tractable to the incarnation of the classical homotopy hypothesis in terms of \cref{PSHH:Thm_Kan_Quillen_Equ}. To begin with, this means we want to obtain tractable versions of the comparison functors between stratified spaces and $\infty$-categories, in terms of a presentation through the stratified singular simplicial set functor and its left adjoint given by stratified realization. Secondly, we want to obtain a better understanding of the homotopy theories of stratified spaces -- $\AltStratD$ ,$\AltStratC$ and $\AltStratCR$, and their fixed poset counterparts -- and how they interact with the $1$-category $\StratN$ as well as classical approaches to stratified homotopy theory, such as the one pursued in \cite{quinn1988homotopically}.
\begin{question}\label{question:explicit_results}
 More specifically, we aim to answer the following questions.
\begin{IntroQuestions}
    \item \label{iQ:Exitcon} Can the equivalence $\AltTopPC
    \simeq \AbStrP$ be presented through Lurie's stratified singular simplicial set (Exit-path) construction, and does $\SingS$ create stratum-preserving categorical equivalences, in the sense that a stratum-preserving map is a categorical equivalence, if and only if $\SingS(f)$ is a Joyal equivalence? If yes, do analogous results hold for the homotopy theories of stratified spaces with varying stratification posets, thereby presenting a global version of the topological stratified homotopy hypothesis?
    \item \label{iQ:LocatSHE} Without additional structure, $\infty$-categorical localizations of $1$-categories are generally difficult objects to study. What can we say about the homotopy theories of stratified spaces from a stratified topological perspective? For example, can we express $\AltStratCR$ as a subcategory of stratified spaces $\textnormal{\textbf{C}}$ localized at stratified homotopy equivalences, analogously to the situation of topological spaces and CW complexes?
    Such a category $\textnormal{\textbf{C}}$ should contain as many stratified spaces of classical interest as possible, to allow us to investigate the stratified homotopy-theoretic properties of geometrically interesting examples through the language of $\infty$-categories. 
    \item \label{iQ:Whatfibrant}
     In the same line of questioning as in the previous question: The category $\StratN$ admits a naive notion of mapping space, with points given by stratified maps and paths given by stratified homotopies. For what stratified spaces, $\tstr$ and $\str[Y]$, can we expect the mapping space $\AltStratCR(\tstr, \tstr[Y])$ to have the homotopy type of this naive mapping space, or similarly, which classical examples of stratified spaces are contained in $\textnormal{\textbf{C}}$?
    % , or more broadly, how the model structures interact with more classical approaches to stratified homotopy theory such as \cite{quinn1988homotopically}.
    \item \label{iQ:DifferenceBetweenTheories} What are the precise relationships between the several stratified homotopy theories introduced above and the more classical approach due to Quinn? In particular, how large are the differences when restricting to classical, geometrical examples of stratified spaces?
    \item \label{iQ:ComputationHocolim} How can $\infty$-categorical construction, such as colimits in a quasi-category (i.e., homotopy colimits), be interpreted in terms of the $1$-category $\StratN$? For instance, when is a pushout diagram of stratified spaces homotopy pushout in one of the homotopy theories of stratified spaces above?
     \item \label{iQ:ExponentialObjects} The sets of stratified maps can even be equipped with the structure of a stratified space, with each stratum corresponding to a map of the underlying posets (see \cite{nand2019simplicial}). In \cite{HughesPathSpaces}, such stratified (decomposed) mapping spaces and their relations with stratified notions of fibrations were studied. The central result states a stratified exit path version of the classical path-space fibration. Can we replicate this result in the stratified homotopy theories described above, using methods of modern abstract homotopy theory?
    \item \label{iQ:ConicallySmooth} What is the relationship of the homotopy theories described above with the homotopy theory of conically smooth stratified spaces discussed in \cite{AFRStratifiedHomotopyHypothesis}?
\end{IntroQuestions}
\end{question}
Classically, such questions concerning the relationship of a homotopy theory and a $1$-category from which it is obtained by localization are answered through the language of model categories. Thus, we may at least partially rephrase the questions above:
\begin{question}\label{question:mod_cat}
    Do the classes of weak equivalences defining the homotopy theories $\AltStratD$, $\AltStratC$ and $\AltStratCR$ (and their fixed, poset counterparts) extend to model structures on $\StratN$ ($\TopPN$)? Specifically, do model structures exist in which stratified spaces of classical interest -- for example, Whitney stratified spaces or, more recently, conically smooth stratified spaces -- are bifibrant? Supposing an affirmative answer, what are the properties of this model category, such as admitting a simplicial structure, cartesian closedness, cofibrant generation?
    Can one prove an answer to the stratified homotopy hypothesis in terms of Quillen equivalences through the adjunction of stratified realization and stratified singular simplicial sets?
\end{question}
In \cite{douteauEnTop,Henriques} model structures on $\StratN$ presenting $\AltStratD$ were defined. However, these model structures fail the criterion of having convenient bifibrant objects: Not even the stratified cone on a closed  manifold is a cofibrant object. In \cite{nand2019simplicial}, the author defined a model structure for the subcategory of such stratified spaces, $\str$, for which $\SingS (\str)$ is a quasi-category and used it to prove a Whitehead theorem for stratified spaces. The existence of a (semi-)model structure as above was left open as a conjecture (\cite[8.4.1]{nand2019simplicial}).
While this makes, for example, (appropriately stratified) piecewise linear pseudo-manifolds bifibrant, a full model category (in particular fibrant replacement) is needed to present the stratified homotopy hypothesis. Finally, in \cite{douteauwaas2021}, we showed that such model structures cannot exist (see \cref{prop:counter_example_mod} for a detailed proof adapted to the setting of this paper).
\subsection{Content of this article}
One main result of this paper is to prove that we can answer \cref{question:mod_cat} affirmatively if we generalize to (left) semi-model categories (see, for example, \cite{BarwickLeftRight,white2023left}). We note that, in practice, these usually turn out to be just as powerful as model categories (see, for example, \cite[Rem. 4.5]{white2023left}. In particular, the existence of left semi-model structures (just semi-model structures henceforth) allows us to obtain answers to \cref{question:explicit_results}.
Building on the work of \cite{haine2018homotopy}, \cite{douteauEnTop}, and \cite{douteau2021stratified}, we have laid the foundation to prove the existence of these semi-model structures in the two preceding papers: \cite{HoLinksWa} and \cite{ComModelWa}. The remaining task is to combine these results as follows:
\begin{enumerate}
    \item In \cite{douteau2021stratified} and \cite{haine2018homotopy} Douteau and Haine each defined simplicial counterparts to their stratified homotopy theories for the case of a fixed poset $\pos$. In \cite{ComModelWa}, we extended these models to model structures for stratified simplicial sets with varying stratification poset, the category of which we denote $\sStratN$,. Doing so, we obtain simplicial model categories $\sStratD$, $\sStratDR$, $\sStratC$ and $\sStratCR$. The latter two of these do, respectively, present the homotopy theories of abstract stratified homotopy types (with varying stratification poset) and of (small) layered $\infty$-categories, which are the $\infty$-categories in which every endomorphism is invertible.  
     We recall these results in \cref{subsec:comb_mod_hostrat}. 
    \item We want to construct model structures for $\AltStratC, \AltStratD, \AltStratCR$ and their fixed poset analogues by transferring the model structures on the simplicial side along the stratified singular simplicial set and realization adjunction. To this end, we need to prove that the functor $\SingS \colon \StratN \to \sStratN$ creates weak equivalences for the respective theories. This is proven in \cref{subsec:abstract_equ}. In fact, more than that, we show that the functor induces homotopy equivalences of categories with weak equivalences (\cref{thm:overview_over_all_hypothesis}). 
    \item Then, in \cref{sec:con_of_mod_cat}, we transfer the model structures from the combinatorial framework to the topological one, through a transfer lemma for semi-model categories (\cref{lem:transf_lem_spec}). Applying this lemma requires a deeper understanding of the interaction of the homotopy link functor with topological stratified constructions, such as pushouts along inclusions of cell complexes. This is achieved through \cite[Thm. \ref{hol:thm:hol_main_res_B}]{HoLinksWa}, which establishes the existence of  certain regular neighborhoods for multistrata interactions in stratified cell complexes. In many ways, these are the topological (as opposed to smooth) analogues of the Unzip construction of \cite{LocalStructOnStrat}:
    % takes in smooth stratified homotopy theory: 
    They provide the necessary glue to connect the geometry/topology of stratified spaces with their homotopy theory.
\end{enumerate}
We then combine these insights to obtain the following result, which -- using standard results on (semi-)model categories -- ultimately addresses \cref{iQ:LocatSHE,iQ:ComputationHocolim,iQ:ExponentialObjects} (see, particularly, \cref{subsec:answers}). For the sake of conciseness, we only cover the case of varying posets in this introduction.
 \begin{ctheorem}{A}[\cref{thm:overview_over_all_hypothesis,cor:semi_mod_transf,thm:ex_of_model_struct_strat_nonref,thm:ex_mod_struct_red,cor:cart_closed}]\label{thm:contents_transf_result}
          The category of stratified topological spaces $\StratN$ admits the structure of three distinct simplicial, cofibrantly generated, and cartesian closed left semi-model categories: $\StratD$, $\StratC$, and $\StratCR$. Weak equivalences are given, respectively, by the poset-preserving diagrammatic, poset-preserving categorical, and categorical equivalences.
          They are right-transferred along the adjunction
                        \[ 
                        \sReal{-} \colon \sStratN \rightleftharpoons \StratN \colon \SingS
                        \]
          from their respective counterpart on $\sStrat$. Moreover, with respect to these transferred structures, the adjunction defines simplicial Quillen equivalences (between the respective topological and simplicial counterparts) that create weak equivalences in both directions.
    \end{ctheorem}
    $\StratC$ is obtained from $\StratD$ in terms of a left Bousfield localization at stratified inner horn inclusions (providing a first step towards an answer to \cref{iQ:DifferenceBetweenTheories}). In turn, $\StratCR$ is obtained from $\StratC$ in terms of a right Bousfield localization (\cref{thm:overview_over_all_hypothesis}). \\
As a corollary of this result and \cite[Thm. \ref{comb:prop:Quillen_Equ_betw_ref_and_ord}]{ComModelWa}, we obtain the following version of the stratified homotopy hypothesis, providing answers to \cref{iQ:Exitcon}:
%%%
\tikzcdset{%
   harpoonlurd/.code={\pgfsetarrows{tikzcd right to-tikzcd right to}},
   harpoonldru/.code={\pgfsetarrows{tikzcd left to-tikzcd left to}},
   harpoonluru/.code={\pgfsetarrows{tikzcd right to-tikzcd left to}},
   harpoonldrd/.code={\pgfsetarrows{tikzcd left to-tikzcd right to}},
}
\makeatletter
\def\leftrightharpoon{%
  \@ifnextchar[{\@lrharp}{\@lrharp[]}%
}
\def\@lrharp[#1]{%
  \arrow[#1, leftharpoonup  ,yshift= 0.225ex]%
  \arrow[#1,rightharpoondown,yshift=-0.225ex]%
}
\makeatother
%%%%
\begin{ctheorem}{B}[\cref{thm:strat_ho_hy_quillen}]\label{thm:contents_strat_hohy_quillen}
     Mapping a simplicial set to its stratified realization (as in \cite{nand2019simplicial}) and conversely mapping a stratified space $\str$ to the underlying simplicial set of $\SingS \str$ induces a Quillen equivalence of (semi-)model categories 
    \[
    \sSetOrd \xrightleftharpoons{\simeq} {\StratCRN}
    \]
    that creates weak equivalences in both directions. \\
    Here, $\sSetOrd$ is the left Bousfield localization of the Joyal model structure on $\sSetN$ that presents layered $\infty$-categories.
\end{ctheorem}
To obtain a similarly convenient situation to the setup of the classical homotopy hypothesis, an answer to \cref{iQ:Whatfibrant} remains to be obtained. We do so in \cref{sec:cofibrants}.
It is an immediate consequence of the way that the model structures are constructed via transfer from stratified simplicial sets that any stratified space which admits a piecewise linear structure (or, more generally, a cell structure) that is compatible with the stratification is cofibrant in the semi-model categories presenting $\AltStratD$ and $\AltStratC$, namely $\StratD$ and $\StratC$. 
Hence, for example, by \cite{TriangulationsGoresky}, all Whitney stratified spaces are cofibrant.
Even more, we show \cref{prop:inductive_cofibrancy}, from which it follows that any stratified space (over a finite poset) whose strata are manifolds that additionally admit a stratified notion of mapping cylinder neighborhood is cofibrant. 
In particular, assuming the correctness of a result of \cite{AFRStratifiedHomotopyHypothesis}, this implies that conically smooth stratified spaces are cofibrant.
The semi-model category $\StratCR$ turns out to be a right Bousfield localization of $\StratC$, whose cofibrant objects are precisely the cofibrant objects in $\StratC$ that fulfill a weak version of the classical frontier condition and have connected strata (\cref{prop:ref_vs_frontier}). Conceptually, these are the stratified spaces (cofibrant in $\StratC$) in which the stratification poset structure arises entirely in terms of the topological relations of the strata.
\\
Faced with several different model structures and homotopy theories for stratified spaces, the obvious question about the precise relationship between these theories, in particular in application to geometric examples, arises. We have already illustrated above that the passage from $\AltStratC$ to $\AltStratCR$ essentially amounts to requiring that the poset structure is intrinsic to the topology of the space. The difference between $\AltStratC$ and $\AltStratD$ (that is, answering \cref{iQ:DifferenceBetweenTheories}) is more subtle.
As the cofibrant objects in the categorical and the diagrammatic semi-model categories agree, the difference between the resulting homotopy theories must lie in the conditions for fibrancy. 
It is a result of \cite{HigherAlgebra} that conically smooth stratified spaces are fibrant in $\StratCR$, and hence also in $\StratC$ and $\StratD$. In \cite{nand2019simplicial}, it was shown that the same holds for Quinn's homotopically stratified spaces (\cite{quinn1988homotopically}). This already covers most classically relevant examples of stratified spaces. Thus, it appears that at least in a geometric scenario there is not much difference between the two homotopy theories $\AltStratC$ and $\AltStratD$ at all. In fact, we make this result rigorous in terms of the following answer to \cref{iQ:Whatfibrant,iQ:DifferenceBetweenTheories}.
\begin{ctheorem}{C}[\cref{prop:equi_of_fibrancy,thm:all_theories_same}\label{mainthm:comparison}]
     Let $\tstr \in \StratN$ be a metrizable stratified space. Then the following conditions are equivalent:
    \begin{fibrancyProp}
        \item $\tstr$ is fibrant in $\StratCR$;
        \item $\tstr$ is fibrant in $\StratC$;
        \item $\tstr$ is fibrant in $\StratD$;
         \item For any pair of strata $[p<q]$, the starting point evaluation map $\HolIP[p<q](\tstr) \to \tstr_{p}$ is a Serre fibration.
        % \item\label{prop:equ_fibrancy_4} For any pair $p < q \in \ptstr$, and for any flag $\J = [p = \cdots = p_k < q = \cdots = q = p_{n_\J}]$ with $k \geq 1$, $\tstr$ has the horn filling property with respect to $\sReal{\Lambda^\J_k \hookrightarrow \Delta^\J}$;
    \end{fibrancyProp}
    In particular, when restricted to metrizable stratified spaces, the homotopy theories defined by $\StratC$ and $\StratD$ (in terms of simplicial categories of bifibrant objects) agree.
\end{ctheorem}
Similarly to the case of the classical homotopy hypothesis, it can be useful to restrict to a class of particularly convenient stratified spaces that mimic some of the properties of CW complexes. This is handled by \cref{cor:perf_strat_sp_are_essentially_surj}, which in particular states that the bifibrant stratified spaces in $\StratCR$ are precisely the retracts of the so-called CFF stratified spaces, i.e., the stratified spaces $\tstr \in \StratN$ fulfilling:
        \begin{enumerate}
            \item $\SingS{\tstr}$ is a quasi-category.
            \item $\tstr$ admits cell structure that is (in a sense which we will specify) compatible with the stratification.
            \item $\tstr$ has non-empty connected strata, and the structure of the stratification poset arises from the classical frontier condition.
        \end{enumerate}
It follows that we may rephrase \cref{thm:contents_strat_hohy_quillen} as the following statement:
\begin{ctheorem}{B'}
    Denote by $\CFFN$ the full subcategory of $\StratN$ given by CFF stratified spaces and let $H_s$ be the class of stratified homotopy equivalences in $\CFFN$.
    Lurie's exit path construction (\cite{HigherAlgebra}) induces an equivalence of $\infty$-categories
    \[
    \textnormal{Exit} \colon \CFFN[H_s^{-1}] \xrightarrow{\simeq} \iCatO
    \]
    where $\iCatO$ denotes the $\infty$-category of small layered $\infty$-categories.
\end{ctheorem} 
At the end of \cref{sec:cofibrants} in \cref{subsec:stratified_homotopy_link_fibrations}, we study the stratified homotopy link fibrations of \cite{HughesPathSpaces} from a model categorical perspective, thus providing further connections between our model categorical approach to stratified homotopy theory and more classical approaches. \\
\\
Finally, let us comment on the relationship of the homotopy theories investigated in this article with the homotopy theory of conically smooth stratified spaces investigated in \cite{AFRStratifiedHomotopyHypothesis} (\cref{iQ:ConicallySmooth}). Supposing the existence of stratified mapping cylinder neighborhoods (tubular neighborhoods), as asserted in \cite[Prop. 8.2.3]{LocalStructOnStrat}, it follows from \cref{prop:inductive_cofibrancy}, \cref{prop:ref_vs_frontier} and \cref{mainthm:comparison} that conically smooth stratified spaces are bifibrant in $\StratCR$.
Note that it follows from this result and \cref{thm:contents_strat_hohy_quillen} that the following two statements are equivalent:
\begin{enumerate}
    \item The exit path construction induces a fully faithful embedding from the $\infty$-category of conically smooth stratified spaces into $\infty$-categories (as asserted in \cite{AFRStratifiedHomotopyHypothesis}).
    \item For any pair of conically smooth stratified spaces $\tstr, \tstr[Y]$, the natural map 
    \[
    \Strat_{\mathcal{C}^{\infty}}(\tstr, \tstr[Y]) \to \Strat (\tstr, \tstr[Y])
    \]
    from the mapping space of conically smooth maps (defined in \cite{AFRStratifiedHomotopyHypothesis}) into the mapping space of continuous stratified maps is a homotopy equivalence of Kan complexes.
\end{enumerate}
In this sense, the smooth version of the stratified homotopy hypothesis as it is stated above is equivalent to a stratified smooth approximation theorem, as conjectured in \cite[Conjecture 1.5.1]{LocalStructOnStrat}.

%% file: 2Preliminaries.tex
\section{Preliminaries}
In this section, we introduce and recall some of the necessary language and notation, especially for $1$-categorical aspects of stratified spaces and higher homotopical frameworks.
\subsection{Models for \texorpdfstring{$(\infty,1)$-categories}{infinity-one-categories}}
    This paper is concerned with investigating several homotopy theories of stratified spaces. By a \textit{homotopy theory} we mean an $(\infty,1)$-category, not restricting to a specific model for $(\infty, 1)$-categories. For the sake of readability, we will always drop the $1$ and use the term $\infty$-category and homotopy theory synonymously. 
    For our purposes, it is going to be extremely useful to have several different models for $\infty$-categories available. Specifically, we are going to use the following models.
    \begin{notation}
        \begin{enumerate}
        \item The relative categories of \cite{BarwickRelative}, given by pairs $(\textnormal{\textbf{C}},W)$ of a $1$-category $\textnormal{\textbf{C}}$ with a wide subcategory $W \subset \textnormal{\textbf{C}}$.
        Names for relative categories will always begin with an italic letter, that is, they will be of the form $\textnormal{\textbf{\textit{N}ame}}$.
        \item Categories enriched over simplicial sets, also called simplicial categories (see \cite{DKsimLocCat,BergnerSimCat}): Names for simplicial categories will always be underlined, that is, written in the form $\underline{\textnormal{\textbf{Name}}}$.
        \item Quasi-categories in the sense of \cite{HigherTopos}: Names for quasi-categories will always be bold and begin with a calligraphic letter, i.e. they will be of the form $\mathcal{N}\textnormal{\textbf{ame}}$.\\
    \end{enumerate}
 If $\textbf{C}$ is an $\infty$-category (in particular a $1$-category), the mapping spaces between objects $X,Y \in \textbf{C}$ will be denoted by $\textbf{C}(X,Y)$. In the case of simplicial categories, these have an explicit model in the obvious way. For relative categories, they are given by the hammock localization of \cite{DwyerKanCalculating,BarWickKanCharacterization}. For quasi-categories, we use any of the equivalent models of mapping spaces of \cite{HigherTopos}.
    \end{notation}
    \begin{notation}
        To study $\infty$-categories through some underlying $1$-category which they are a localization of, we are going to make use of the language of (semi)model categories (see \cite{hirschhornModel} for model categories and \cite{BarwickLeftRight,white2023left} for a good overview of semi-model categories). We will usually denote ordinary categories in bold letters in the form $\textnormal{\textbf{Name}}$).
        The model structure will always be marked by adding some additional ornamentation - in the form $\underline{\textnormal{\textbf{Name}}}^{\mathfrak s}$ or $\textnormal{\textbf{Name}}^{\mathfrak s}$ - to the name of the underlying simplicial or $1$-category.
    \end{notation}
    \begin{notation}
        Functor categories, between two $1$-categorical, simplicial categories or quasi-categories $\cat[C], \cat[D]$ will always be denoted in the form $\FunC(\cat[C],\cat[D])$. This is to be understood in the sense that the type of categories inserted into $\FunC(-,-)$ specifies whether the resulting functor category is itself a $1$-category, quasi-category, or simplicial category. In case that the source is a $1$-category, and the target a quasi or simplicial category, we treat the $1$-category, respectively, as a quasi-category (via its nerve) or as a simplicial category with discrete mapping spaces.
        At times, we will also use exponential notation $\cat[D]^{\cat[C]}$ to refer to functor categories. 
    \end{notation}
    \begin{notation}\label{not:conventions_for_homotopy_theories}
        Often we are going to pass between different models of a homotopy theory. These passages will always follow the following ruleset for nomenclature:
        \begin{enumerate}
            \item Starting with a relative category $\textnormal{\textbf{\textit{N}ame}}$ or simplicial category $\underline{\textnormal{\textbf{Name}}}$, the underlying $1$-category will be denoted $\textnormal{\textbf{Name}}$. Similarly, if $\underline{\textnormal{\textbf{Name}}}^{\mathfrak s}$ is a simplicial model category, then we denote its underlying $1$-model category by $\textnormal{\textbf{Name}}^{\mathfrak s}$.
            \item Starting with a model category $\textnormal{\textbf{Name}}^{\mathfrak s}$, we denote by $\textnormal{\textbf{\textit{N}ame}}^{\mathfrak s}$ the relative category obtained by its underlying $1$-category, together with the wide subcategory of weak equivalences.
            \item Starting with a relative category $\textnormal{\textbf{\textit{N}ame}} = (\textnormal{\textbf{Name}}, W)$, we denote by $\mathcal{N}\textnormal{\textbf{ame}}$ the quasi-category $\textnormal{\textbf{Name}} [W^{-1}]$ obtained by taking the nerve of the underlying $1$-category of $\textnormal{\textbf{\textit{N}ame}}$, and then localizing at $W$. 
            \item In particular, following this language, if we start with a model category $\textnormal{\textbf{Name}}^{\mathfrak s}$, then the relative category $\textnormal{\textbf{\textit{N}ame}}^{\mathfrak s}$ and the quasi-category $\mathcal{N}\textnormal{\textbf{ame}}^{\mathfrak s}$ model the same $\infty$-category.
            \item If we start with a simplicial model category $\underline{\textnormal{\textbf{Name}}}^{\mathfrak s}$, then we denote by $\underline{\textnormal{\textbf{Name}}}^{\mathfrak s,o}$ the full simplicial subcategory of bifibrant objects. Following our notation, under the Quillen equivalence between simplicial categories and quasi-categories of \cite{BergnerInf}, $\underline{\textnormal{\textbf{Name}}}^{\mathfrak s,o}$ models the same $\infty$-category as $\mathcal{N}\textnormal{\textbf{ame}}^{\mathfrak s}$ (this follows from \cite[Prop. 4.8]{FunctionComplexesDwyerKan} together with \cite[Prop. 1.2.1]{HinichDwyerLoc}). 
        \end{enumerate}
        For the sake of completeness, one could of course also introduce notation to pass from relative to simplicial categories, etc., but we will not make use of this change of model in this article, and thus do not take this extra step.
    \end{notation}    
\subsection{Poset-stratified spaces}\label{subsec:top_notation}
Next, let us introduce the basic objects of study: topological spaces that are stratified over a partially ordered set (poset). We begin with notation for the world of posets.
\begin{notation} We are going to use the following terminology and notation for partially ordered sets, drawn partially from \cite{douSimp} and \cite{haine2018homotopy}:
\begin{itemize}
    \item  We denote by $\Pos$ the category of partially ordered sets, with morphisms given by order-preserving maps.
    \item  We denote by $\Delta$ the full subcategory of $\Pos$ given by the finite, nonempty, linearly ordered posets of the form $[n]:= \{0, \cdots,n \}$, for $n \in \mathbb{N}$.
    \item  Given $P \in \Pos$, we denote by $\catFlag$ the slice category $\Delta_{/ \pos}$. 
    That is, objects are given by arrows $[n] \to \pos$ in $\Pos$, $n \in \mathbb N$, and morphisms are given by commutative triangles.
    \item We denote by $\catRFlag$ the \define{subdivision of} $\pos$, given by the full subcategory of $\catFlag$ of such arrows $[n] \to \pos$, which are injective.
    \item The objects of $\catFlag$ are called \define{flags of} $\pos$. We represent them by strings $[p_0 \leq \cdots \leq p_n]$, of $p_i \in \pos$. 
    \item Objects of $\catRFlag$ are called regular \define{flags of} $\pos$. We represent them by strings $\standardFlag$, of $p_i \in \pos$.
\end{itemize}
\end{notation}
\begin{notation}\label{pshh:not:notation_top_spaces}
    $\TopN$ is going to denote either of the following three categories of topological spaces. \begin{enumerate}
        \item The category of \define{all topological spaces}, which we will also refer to as \define{general} topological spaces.
        \item The category of compactly generated topological spaces, i.e. such spaces which have the final topology with respect to compact Hausdorff spaces (see, for example, \cite{Rezk2017COMPACTLYGS}).
        \item The category of $\Delta$-generated topological spaces, i.e. such spaces which have the final topology with respect to realizations of simplices, or equivalently just with respect to the unit interval (compare \cite{duggerDelta,GaucherDelta}).
    \end{enumerate}
    We denote by $\sSet$ the simplicial category of simplicial sets, i.e. the category of set-valued presheaves on $\Delta^{\op}$, equipped with the canonical simplicial structure induced by the product (see \cite{HigherTopos} for all of the standard notation used for simplicial sets).
    We denote by $\real{-} \colon \sSetN \to \TopN$ the realization functor of simplicial sets and by $\Sing \colon \TopN \to \sSetN$ its right adjoint.
    $\TopN$ naturally carries the structure of a simplicial category, tensored and powered over $\sSetN$, induced by left Kan extension of the construction
    \[
    T \otimes \Delta^n := T \times \real{\Delta^n}.
    \]
    \end{notation}
    We denote the resulting simplicial category by $\Top$.
    Furthermore, we will always consider $\Top$ to be equipped with the Quillen model structure \cite{Quillen}, which makes
    $\real{-} \dashv \Sing$ a simplicial Quillen equivalence, between $\Top$ and $\sSet$ (with the latter equipped with the Kan-Quillen model structure), which creates all weak equivalences in both directions.
    \begin{remark}
    Note that one commonly only defines the simplicial structure for compactly or $\Delta$-generated spaces. However, this is mostly due to the fact that for general topological spaces $T$ and infinite simplicial sets $K$, the tensor $T\otimes K$ does not agree with the inner product $T \times \real{K}$. Instead, it is given by a colimit of products of $T$ with the simplices of $K$. Similarly, the power $T^{K}$ is not given by an internal mapping space (which does not even necessarily exist for arbitrary $K$) but by the limit of the mapping spaces with source given by simplices of $K$, which are equipped with the compact-open topology. 
    We want to emphasize that for the resulting homotopy theory the choice in underlying set-theoretic assumptions on topological spaces is immaterial.
    \end{remark}
    For the remainder of this subsection, we fix some category of topological spaces $\TopN$ as in \cref{pshh:not:notation_top_spaces}.
    \begin{notation}
    Having fixed a category of topological spaces $\TopN$, we then use the following notation for stratified topological spaces (all of these constructions already appear in \cite{douteauEnTop} among other places). The notation and language is mostly analogous with the language for the simplicial framework.
    \begin{itemize}
    \item We think of the $1$-category $\Pos$ as fully faithfully embedded in $\TopN$, via the Alexandrov topology functor $\mathcal{A} \colon \Pos \to \TopN$, equipping a poset $P$ with the topology where the closed sets are given by the downward closed sets. By abuse of notation, we usually just write $\pos$, for the Alexandrov space corresponding to the poset $\pos$ (compare \cite[Def. 2.2]{douteauwaas2021}). 
    \item For $\pos \in \Pos$, we denote by $\TopPN$ the slice category $\TopN_{/\pos}$. We treat $\TopPN$ as a simplicial category, denoted $\TopP$, with the structure inherited from $\Top$ (see \cite[Recol. 2.13]{douteauwaas2021}, for a detailed definition in the $\Delta$-generated case).
    \item Objects of $\TopP$ are called \define{$\pos$-stratified spaces}. They are given by a tuple $(\tstr,s\colon \utstr \to \pos)$. In the literature, a $P$-stratified space $(\tstr,s\colon \utstr \to \pos)$ is often simply referred to by its underlying space $\utstr$, omitting the so-called \define{stratification} $s \colon \utstr \to \pos$. We are not going to adopt this notation here and will generally use calligraphic letters for stratified spaces and stick to the notational convention $\tstr =(\utstr, \ststr)$ to refer to the underlying space and stratification.
    \item  Morphisms in $\TopPN$ are called \define{stratum-preserving maps.} 
   \item Given a map of posets $f\colon Q \to \pos$ and $\tstr \in \TopPN$, we denote by $f^*\tstr \in \TopPN[Q]$ the stratified space $\utstr \times_{\pos} Q \to Q$. We are mostly concerned with the case where $f$ is given by the inclusion of a singleton $\{p \}$, of a subset $\{q \sim p \mid q \in \pos \}$, for $p \in \pos$ and $\sim$ some relation on the partially ordered set $\pos$ (such as $\leq$), or more generally, a subposet $Q \subset \pos$. We then write $\tstr_{p}$ (or, respectively, $\tstr_{\sim p}$, $\tstr_Q$) instead of $f^*\tstr$. The spaces $\tstr_{p}$, for $p \in \pos$ are called the \define{strata of} $\tstr$.
   \item For $f\colon Q \to \pos$ in $\Pos$, we denote by $f_!$ the left adjoint to the simplicial functor $f^* \colon \TopP[Q] \to \TopP$. $f_!$ is given on objects by $(\ststr \colon \utstr  \to Q) \mapsto (f \circ \ststr \colon \utstr \to Q \to \pos)$.
    \item Let $\TopN^{[1]}$ be the category of arrows of topological spaces.
    We denote by $\StratN$ the category of all (poset-)stratified spaces, given by the full sub-category of $\TopN^{[1]}$ of such arrows $X \to \pos$, where $X \in \TopN$ and $\pos\in \Pos$ is a poset (equipped with the Alexandrov topology). In particular, every object of $\StratN$ is given by a $\pos$-stratified space, for some $\pos\in \Pos$, and a morphism $(X \to \pos) \to (Y \to \pos[Q]$) is given by a pair of morphisms $f\colon X \to Y$ and $g \colon \pos \to \pos[Q]$, where $f$ is a continuous map and $g$ can be seen as a map of posets, 
    making the obvious square commute (see also \cite[Def. 2.11]{douteauwaas2021}).
    Morphisms are called \define{stratified} maps. 
    \item Given $\tstr \in \StratN$, we are going to use the notations $\tstr=(\utstr, \ptstr, \ststr)$ to refer, respectively, to the underlying space, the poset, and the stratification and proceed analogously for morphisms.
    \item $\StratN$ is equipped with the structure of a simplicial category (tensored and powered over $\sSet$), $\Strat$,  with the simplicial structure induced by $\tstr \otimes \Delta^n = (\utstr \times \real{\Delta^n}  \to \utstr \to \ptstr)$. Note that since $\tstr \otimes -$ preserves colimits, this means that $\tstr \otimes \emptyset$ is not stratified over $\pos$, but instead over the empty poset.
    Simplicial homotopies, that is, homotopies with respect to the cylinder $- \otimes \Delta^1$ in $\Strat$, are called \define{stratified homotopies}. Simplicial homotopy equivalences are called \define{stratified homotopy equivalences}.
    \item The forgetful functor $\Strat \to \Top$, $\tstr \mapsto \utstr$, has a right adjoint. It is given by mapping $T \in \Top$ to the trivially stratified space $( T \to [0] )$. By abuse of notation, we will often write $T$ to refer to the trivially stratified space associated to a space $T$.
    \end{itemize}
    \end{notation}
    \begin{remark}\label{rem:bicomplete}
        Both $\StratN$ and $\TopPN$, for $P \in \Pos$, are bicomplete categories (see, for example, \cite{douteauEnTop}). Limits and colimits in $\TopPN$ are simply given by the limits and colimits in a slice category. Colimits in $\StratN$ are computed by taking the colimit both on the space and on the poset level. Limits in $\StratN$ are computed by taking the limit on the space and poset level, and then pulling back the map $\varprojlim_i X_i \to \varprojlim \mathcal{A}(P_i)$ along the natural comparison map $\mathcal{A}(\varprojlim_{i} P_i) \to \varprojlim_{i} \mathcal{A}(P_i)$, which is an isomorphism for finite diagrams.
    \end{remark}
    \begin{construction}\label{con:strat_map_space}
        If our choice of $\TopN$ is Cartesian closed (that is, if we are in the compactly generated or $\Delta$-generated case) then $\StratN$ is also a Cartesian closed category (see also \cite{nand2019simplicial}, for the slightly different setting of surjectively stratified spaces). 
        We use the notation $X^Y$ to refer to exponential in a Cartesian closed category.
        To be able to distinguish exponential objects in posets from those in topological spaces, we use the notation $\mathcal{A}(\pos)$ to denote the Alexandrov space associated to a poset $\pos$.
        Given a stratified space $\tstr \in \StratN$, the right adjoint $-^{\tstr}$ to the functor $- \times \tstr \colon \StratN \to \StratN$ is constructed as follows.  Given $\tstr[Z]$ in $\StratN$, consider a pullback square in $\TopN$
        \begin{diagram}
         {}\mathcal{A} ({\ptstr[Z]}^{\ptstr[X]} )  \times_{\mathcal{A}(\utstr[\ptstr[Z]])^{\utstr[X]}}  \utstr[Z]^{\utstr[X]}  \arrow[rr] \arrow[d]    & & \utstr[Z]^{\utstr[X]} \arrow[d] \\
          \mathcal{A} ({\ptstr[Z]}^{\ptstr[X]} )  \arrow[r] & \mathcal{A}(\ptstr[Z])^{\mathcal A(\ptstr[X])} \arrow[r] & \mathcal{A}(\ptstr[Z])^{\utstr[X]} \spacecomma
        \end{diagram}
        where $\ptstr[Z]^{\ptstr}$ is the poset obtained by equipping $\Pos ( \ptstr, \ptstr[Z])$ with the poset structure given by 
        \[
        f \leq g :\iff f(p) \leq g(p), \forall p \in \ptstr,\]
        which defines the exponential object $\ptstr[Z]^{\ptstr}$ in $\Pos$. The lower left horizontal map between the resulting space equipped with the Alexandrov topology into the mapping space $\mathcal{A}(\ptstr[Z])^{\mathcal{A}(\ptstr[X])}$ is always continuous, but generally only a homeomorphism if $\ptstr$ and $\ptstr[Z]$ are finite (\cite[Cor. 2.2.11]{may2016finite}). The stratified mapping space $\tstr[Z]^{\tstr}$ is defined by the left vertical of this pullback diagram \footnote{Note that the underlying set of $\mathcal{A} ({\ptstr[Z]}^{\ptstr[X]} ) \times_{\mathcal{A}(\utstr[\ptstr[Z]])^{\utstr[X]}}  \utstr[Z]^{\utstr[X]}$  is in natural bijection with $\StratN (\tstr, \tstr[Z])$, which shows that we can interpret this construction as a choice of topology on $\StratN (\tstr, \tstr[Z])$, which will generally be finer than the initial topology inherited from $\utstr[Z]^{\utstr[X]}$.}. This construction (and its functoriality in stratified maps) induces the right adjoint $-^{\tstr} \colon \StratN \to \StratN$ to $- \times \tstr$. In particular, we may treat $\StratN$ as a symmetric monoidal category (with monoidal structure induced by the product) enriched over itself.
       
    \end{construction}
\subsection{Stratified simplicial sets and stratified realization}
The approach to constructing and investigating homotopy theories of stratified spaces that we take in this article is to transfer homotopy theoretic structure from the combinatorial to the topological world. Let us quickly recall some notation and terminology concerning stratified simplicial sets, as introduced and investigated in \cite{douSimp,haine2018homotopy}. We have surveyed and expanded the results on homotopy theories of stratified simplicial sets in \cite{ComModelWa}, to which we refer for more details. 
\begin{notation}
    We use the following terminology and notation for (stratified) simplicial sets, drawn partially from \cite{douSimp} and \cite{haine2018homotopy}:
    \begin{itemize}
    \item When we treat $\sSet$ as a model category, this will generally be with respect to the Kan-Quillen model structure (see \cite{Quillen}), unless otherwise noted. When we use Joyal's model structure for quasi-categories (\cite{joyalNotes}) instead, we will denote this model category by $\sSetJoy$. 
    \item We think of $\Pos$ as being fully faithfully embedded in $\sSet$, via the nerve functor (compare \cite{haine2018homotopy}). By abuse of notation, we just write $\pos$, for the simplicial set given by the nerve of $\pos \in \Pos$. 
    \item For $\pos \in \Pos$, we denote by $\sSetPN$ the slice category $\sSetN_{/\pos}$, which is equivalently given by the category of set-valued presheaves on $\catFlag$. We treat $\sSetPN$ as a simplicial category, denoted $\sSetP$, with the structure inherited from $\sSet$ (see \cite[Recol. 2.21.]{douteauwaas2021}). Objects of $\sSetP$ are called \define{$\pos$-stratified simplicial sets}. 
     \item Most of the remaining language and notation we are going to use for stratified simplicial sets can be copied mutatis mutandis from the topological setting in \cref{subsec:top_notation}. See \cite[Subsec. \ref{comb:subsec:simp_lang_not}]{ComModelWa}, for a detailed list.
    \item The forgetful functor $\sStratN \to \sSetN$, $\str \mapsto \ustr$, which will be denoted $\forget$, has a right adjoint and a left adjoint. The left adjoint is given by left Kan extending the functor on simplices: $\Delta^n \mapsto \{ \Delta^n \xrightarrow{1_{\Delta^n}} \Delta^n = [n] \}$. We denote it by $\lstr \colon \sSetN \to \sStratN$. The right adjoint is given by mapping $K \in \sSetN$ to the trivially stratified simplicial set $\{ K \to [0] \}$. By abuse of notation, we will often write $K$ to refer to the trivially stratified simplicial set associated to a simplicial set $K$.
    \end{itemize}
    \end{notation}
    \begin{notation}
    We are going to need some additional notation for flags and stratified simplices.
    \begin{itemize}
    \item For a flag $\J = [p_0 \leq \cdots \leq p_n] \in \Delta_P$, we write $\Delta^\J$ for the image of $\J$ in $\sSetPN$ under the Yoneda embedding $\catFlag \hookrightarrow \sSetPN$. Equivalently, $\Delta^\J$ is given by the unique simplicial map $\Delta^{n} \to \pos$ mapping $i \mapsto p_i$. $\Delta^\J$ is called the \define{stratified simplex associated to} $\J$.
     \item Given a stratified simplex $\Delta^\J$, for $\J = [p_0 \leq \cdots \leq p_n]$, we write $\partial \Delta^\J$ for its \define{stratified boundary}, given by the composition $\partial \Delta^n \to \Delta^n \to \pos$.
    \item Furthermore, for $0 \leq k \leq n$, we write $\Lambda^\J_k \subset \Delta^\J$, for the stratified subsimplicial set given by the composition $\Lambda^n_k \to \Delta^n \to \pos$ (we use the horn notation as in \cite{HigherTopos}). The stratum-preserving map $\Lambda^\J_k \hookrightarrow \Delta^\I$ is called the \define{stratified horn inclusion associated to } $\J$ and $k$. The inclusion $\Lambda^\J_k \hookrightarrow \Delta^\I$ is called \define{admissible}, if $p_k = p_{k+1}$ or $p_{k} = p_{k-1}$. The inclusion $\Lambda^\J_k \hookrightarrow \Delta^\I$ is called \define{inner} if $0 <k <n$.
    % \item Using the fully faithful (and continuous) embedding $\catFlag \hookrightarrow \sSetPN$, we extend the base change notation for stratified simplicial sets to flags.  That is, for $f \colon Q \to \pos$ we write $f^*\J$ for the unique flag of $Q$ corresponding to $f^*(\Delta^\J)$. We use the same shorthand notation for subsets $Q \subset \pos$. For example, $\J_{\leq p}$ is the flag obtained from $\J$, by removing all entries not lesser equal to $p$.
    \item It will also be convenient to have a concise notation for the images of simplices, horns, and boundaries under $\lstr \colon \sSetN \to \sStratN$. These are denoted by replacing the exponent $n \in \mathbb N$, by the poset $[n]$. That is, we write $ \stratSim := \lstr(\Delta^n)$, $ \stratBound := \lstr(\partial \Delta^n)=$, $\stratHorn := \lstr(\Lambda^n_k)$, for $0 \leq k \leq n$. 
\end{itemize}
    \end{notation}
     \begin{recollection}\cite{douSimp}
        For fixed $P \in \Pos$, the two simplicial categories $\TopP$ and $\sSetP$ are connected through a realization and functor of singular simplices type of adjunction, denoted 
        \begin{align*}
            \sReal{-} \colon \sSetP \rightleftharpoons \TopP \colon \SingS .
        \end{align*}
        The left adjoint is constructed by mapping a stratified simplex $\Delta^\J \to \pos$, with $\J = [ p_0 \leq \dots  ,\leq p_n]$, to the stratified space
        \begin{align*}
            \real{\Delta^n} &\to \pos \\
            x &\mapsto \sup\{p_i \in \J \mid x_i > 0\},
        \end{align*}
        where we consider $\real{\Delta^n}$ as embedded in $\mathbb R^{n+1} \cong \mathbb R^{\J}$. 
        If we consider $\sSetPN$ as the category of set-valued presheaves on $\Delta_P$, then by the logic of a nerve and realization functor, $\SingS \tstr$ is hence given the stratified simplicial set
        \[
        \SingS\tstr(\J) = \TopPN(\sReal{\Delta^\J}, \tstr)
        \]
        with the obvious structure morphisms. 
        % Equivalently, we may describe $\SingS \tstr$ via the left vertical in the pullback diagram
        % \begin{diagram}
        %     K \arrow[d] \arrow[r, hook]& \Sing \utstr \arrow[d] \\
        %     P \arrow[r, hook] & \Sing P \spacecomma
        % \end{diagram}
        % where $\pos \hookrightarrow \Sing\pos $ is the adjoint map to the terminal map $\sReal{P} \to \pos$. 
        The adjunction $\sReal{-} \dashv \SingS$ is simplicial.
    \end{recollection}
    \begin{recollection}[{\cite[Recol. 2.23]{douteauwaas2021}}]
        The two functors $\sReal{-}$ and $\SingS$ are compatible with functoriality in the poset, in the sense that for any morphism $f \colon P \to Q$ there are natural isomorphisms 
        \begin{align*}
          \sReal{-}  f_! &\cong f_!  \sReal{-}\\
          \SingS f^* &\cong f^*  \SingS.
        \end{align*}
        It follows that the adjunctions
        \begin{align*}
            \sReal{-} \colon \sSetP \rightleftharpoons \TopP \colon \SingS,
        \end{align*}
        extend to a global adjunction
        \begin{align*}
            \sReal{-} \colon \sStrat \rightleftharpoons \Strat \colon \SingS,
        \end{align*}
        denoted the same, by a slight abuse of notation. Specifically, the maps of simplicial mapping spaces are given by applying the fixed poset versions of $\sReal{-} \dashv \SingS$ component-wise, under the identifications
        \[
        \Strat (\tstr[S], \tstr[T]) \cong \bigsqcup_{f \in \Pos(\ptstr[S],\ptstr[T])} \TopP[{\ptstr[S]}](\tstr[S],f^*\tstr[T]) \cong \bigsqcup_{f \in \Pos(\ptstr[S],\ptstr[T])} \TopP[{\ptstr[T]}](f_!\tstr[S],\tstr[T]);
        \]
        \[
          \sStrat (\str, \str[Y]) \cong \bigsqcup_{f \in \Pos(\pstr,\pstr[Y])} \sSetP[{\pstr}](\str,f^*\str[Y]) \cong \bigsqcup_{f \in \Pos(\pstr,\pstr[Y])} \sSetP[{\ptstr[Y]}](f_!\str,\str[Y]) .
        \]
    \end{recollection}

%% file: 3StratifiedHomotopyTheories.tex
    \section{Homotopy theories of stratified objects}
    Equipping $\StratN$ with a simplicial structure and powers (even an enrichment over Kan complexes) allows us to study stratified spaces from a homotopy-theoretic perspective. It follows from \cite[\S 2.5,2.6]{DwyerKanEqu} that the $\infty$-category defined by $\Strat$ is equivalently given by localizing $\StratN$ at stratified homotopy equivalences $H_s$\footnote{We make now use of this statement, so we omit a proof, but this formally follows from the fact that the powering of $\Strat$ under $\sSetN$ can be used to define an equivalence of relative categories between $(\StratN,H_s)$ and the flattening of $\Strat$.}. It is, however, a general paradigm in abstract homotopy theory that it is often more fruitful to consider a class of weak equivalences larger than the homotopy equivalences with respect to the simplicial structure (for example, weak homotopy equivalences of topological spaces, instead of homotopy equivalences).
    We may encode this information in terms of a relative category given by equipping $\StratN$ with the additional data of a class of weak equivalences $W$. The $\infty$-categorical localization $\StratN[W^{-1}]$ then defines a homotopy theory of stratified spaces.
    In the presence of sufficient extra structure, this process of localizing a larger class of maps is often equivalent to restricting the simplicial category to a subclass of particularly nice (i.e., bifibrant) objects (formally, this is \cite[Prop. 1.1.10]{goerss2005moduli} or, more classically, \cite[§7]{FunctionComplexesDwyerKan}). For example, if one is interested in studying topological manifolds, by Whitehead's theorem (and the fact that topological manifolds have the homotopy type of CW complexes), it is often perfectly sufficient to do so in the setting of weak homotopy equivalences. 
    Following this perspective from classical homotopy theory, the question of a good notion of weak equivalence of stratified spaces and hence after a convenient homotopy theory of (poset) stratified spaces arose. 
    \subsection{Homotopy theories of stratified topological spaces}\label{subsec:top_homotopy_theories}
    The question after a convenient homotopy theory of stratified spaces was already investigated in \cite{quinn1988homotopically} by Quinn, who took the approach of restricting the class of objects, rather than increasing the class of weak equivalences. Quinn followed the insight that in geometric scenarios a stratified space $\tstr[W]$, with finitely many strata and minimal stratum $p$, can be decomposed into a diagram 
    \[
    \tstr[W'] \hookleftarrow \tstr[L] \to \tstr[W]_p 
    \]
    with $\tstr \to \tstr[W]_p$ a fiber bundle with fiber stratified over $\pos_{>p}$, and $\tstr[W'] \hookleftarrow \tstr[L]$ a stratum-preserving boundary inclusion over $\pos_{>p}$. If we inductively repeat this process with $\tstr[W]_{>p}$ and $\tstr[L]$, we ultimately end up with a diagram of spaces indexed over the regular flags $\pos$, that is, over $\sd(\pos) ^{\op}$. In a less geometric scenario, we may not have access to the geometric link $\tstr[L]$ and its further decomposition into smaller pieces. However, we may still consider a homotopy-theoretic analogue.    \begin{recollection}\label{recol:pshh_homotopy_link}
                Let $\tstr \in \TopPN$, and $\I \in \sd(\pos)$. The $\I$-th homotopy link of $\tstr$, $\HolIP(\tstr) \in \TopN$, is the topological space obtained by equipping $\TopPN( \sReal{\Delta^\I}, \tstr)$ with the subspace topology inherited from the compact-open topology (or the respective Kelleyfication thereof, to end up in compactly generated or $\Delta$-generated spaces). If $\I = \{p\}$, then $\HolIP(\tstr) = \tstr_p$ and if $\I = [p<q]$, then $\HolIP( \tstr)$ is the space of paths which start in $\tstr_p$ and immediately exit into the $q$-stratum. This construction is functorial, both in $\I$ and $\tstr$, inducing a right-adjoint, simplicial functor
                    \[
                    \HolIP[] \colon \TopP \to \FunC( \sd(\pos)^{\op}, \Top).
                    \]
                It will often be more convenient to present $\HolIP \tstr$ as a simplicial set. By abuse of notation, we will also denote by $\HolIP[]$ the composition 
                \[
                    \HolIP[] \colon \TopP \to \FunC( \sd(\pos)^{\op}, \Top) \xrightarrow{\Sing_*} \FunC( \sd(\pos)^{\op}, \sSet).
                    \]
                If $\tstr[W]$ is a Whitney stratified space with two strata, then the diagram $\HolIP[](\tstr[W])$ is weakly equivalent to the geometric decomposition diagram illustrated above. 
    \end{recollection}
    Quinn studied a class of (metrizable) stratified spaces, for which the natural evaluation maps $\HolIP[p<q]{\tstr} \to \tstr_p$ are Hurewicz fibrations and which additionally fulfill a cofibrancy condition for inclusions of strata (see \cite{quinn1988homotopically} for details), so-called homotopically stratified spaces (also called homotopically stratified sets). In this type of framework he showed, for example, a stratified analogue of the s-cobordism theorem. In \cite{miller2013}, Miller proved a Whitehead-style theorem for homotopically stratified spaces: A stratum-preserving map between the latter is a stratified homotopy equivalence if and only if it induces homotopy equivalences on all pairwise homotopy links and strata. 
    Inspired by this, both Douteau and Henriques (see \cite{douteauEnTop,Henriques}) independently suggested the following class of weak equivalences for stratified spaces.
    \begin{recollection}\label{rec:douTheoryFixed}\cite{douteauEnTop,douteauwaas2021}
        We call a stratum-preserving map $\tstr \to \tstr[Y] \in \TopPN$ a \define{diagrammatic equivalence}, if and only if it induces weak homotopy equivalences on all generalized homotopy links $\HolIP$, $\I \in \sd(\pos)$. 
        We denote by $\RelTopPD$, the relative category obtained by equipping $\TopPN$ with the class of diagrammatic equivalences. It follows from \cite[Thm. 3]{douteauEnTop} that the homotopy link functor 
        \[
        \HolIP[] \colon \TopP \to \FunC (\sd(\pos)^{\op}, \Top) 
        \]
        induces an equivalence of quasi-categories
        \[
        \AltTopPD \simeq  \FunC (\sd(\pos)^{\op}, \TopN)[\{\text{pointwise weak equivalences}\}^{-1}] \simeq  \FunC (\sd(\pos)^{\op},\iSpaces).
        \]
    \end{recollection}
    This result can be extremely useful insofar, as it allows one to study stratified spaces in terms of presheaves on a fairly simple category, which is a homotopy theoretical setting that is well understood. At first glance, if one takes the perspective of Miller's Whitehead theorem for homotopically stratified sets, it may seem surprising that generalized homotopy links of flags $\I$ containing more than two elements are part of the data detecting weak equivalences. Roughly speaking, higher homotopy links cannot be ignored, as even though they are not necessary to detect stratified homotopy equivalences between two homotopically stratified spaces, they nevertheless provide obstructions for such maps to exist at all. Compare this to the situation of equivalences of categories: Whether a functor is an equivalence of categories may be detected in terms of objects and hom-sets, but whether a functor exists at all requires us to consider higher-dimensional data. Namely, one needs to take into account the composition laws in the categories involved.\\
    \\
    There is, however, another way of homotopically approximating the framework of homotopically stratified spaces. It follows from \cite[A.5]{HigherAlgebra} or \cite[Thm. 8.1.2.6.]{nand2019simplicial} that the homotopy link diagrams arising from classical examples of stratified spaces usually have the following property:
    \begin{recollection}\label{rec:Haine_theory_fixed}\cite{haine2018homotopy,ComModelWa}
        A diagram $D \in \FunC (\sd(\pos)^{\op},\iSpaces)$ is called a \define{d\'ecollage}, if for every regular flag $\I = [p_0 < \cdots < p_n] \in \sd(\pos)$, the canonical morphism
        \[
        D({\I}) \to D{(p_0 < p_1)} \times_{D({p_1})} \cdots \times_{D({p_{n-1}})} D{(p_{n-1} < p_n)}.
        \]
        is a weak homotopy equivalence (that is, an isomorphism in the $\infty$-category of spaces). If $D$ is presented by a commutative diagram in topological spaces or simplicial sets, this condition is equivalent to the natural map
        \[
        D({\I}) \to D{(p_0 < p_1)} \times^{H}_{D({p_1})} \cdots \times^{H}_{D({p_{n-1}})} D{(p_{n-1} < p_n)}. 
        \]
        into the iterative homotopy pullback being a weak homotopy equivalence.
        In particular, if the homotopy link diagram $\HolIP[] (\tstr)$ of a stratified space $\tstr$ is a d\'ecollage, then, roughly speaking, exit paths in $\tstr$ admit compositions which are unique up to higher coherence.
       The homotopy theory of such d\'ecollages can be constructed in terms of a left Bousfield localization of $\SDiag :=\FunC (\sd(\pos)^{\op}, \sSet)$, equipped with the injective model structure, by localizing at such weak equivalences, which are local with respect to d\'ecollages (see \cite[Subsec. \ref{comb:subsec:decollages}]{ComModelWa}). In other words, $f\colon E \to E'$ is a weak equivalence in the model structure that presents d\'ecollages, if and only if for every injectively fibrant diagram $D \in \SDiag$ that has the d\'ecollage property the induced map of mapping spaces
        \[
            \SDiag (E',D) \to \SDiag (E,D)
        \]
        is a weak homotopy equivalence. We denote this left-Bousfield localization of the injective model structure on $\SDiag$ by $\SDiag^\dec$.
        We say that a stratum-preserving map $f\colon \tstr \to \tstr[Y] \in \TopPN$ is a \define{categorical equivalence}, if and only if the induced morphism $\HolIP[] (\tstr ) \to \HolIP[] (\tstr[Y])$ is a weak equivalence in $\SDiag^\dec$, and denote the corresponding relative category by $\TopPCN$. Categorical weak equivalences are precisely the weak equivalences of stratified spaces suggested by Haine in \cite{haine2018homotopy}. It follows by construction and \cite[Thm. 3]{douteauEnTop} that the induced functor of quasi-categories
        \[
        \HolIP[] \colon \AltTopPC \to \AltDiag^{\dec}
        \]
        is an equivalence. That is, if we localize categorical equivalences in $\TopPN$, we obtain the homotopy theory of d\'ecollages. \cite[Thm. 1.1.7]{Exodromy, haine2018homotopy} states that $\AltDiag^{\dec}$ is in turn equivalent to the homotopy theory of quasi-categories with a conservative functor in $\pos$, so-called \define{abstract stratified homotopy types}. This already shows that $\AltTopPC$ fulfills a version of a stratified homotopy hypothesis (\cite{haine2018homotopy}).  
    \end{recollection}
    However, in \cite{haine2018homotopy} it was not yet known whether the equivalence between abstract stratified homotopy types and $\AltTopPC$ could be constructed on the nose through the stratified singular simplicial set construction. In fact, it was not known whether any categorical equivalence even induced a categorical equivalence (Joyal-equivalence) on stratified singular simplicial sets. We are going to prove that this is indeed the case (see \cref{thm:overview_over_all_hypothesis}). Before we do so, let us generalize from the setting of stratum-preserving maps to flexible posets. There are two apparent ways to generalize to the case of a flexible poset. The first essentially amounts to first requiring an isomorphism on the poset level, and then a weak equivalence under the resulting of homotopy theories over fixed posets:
    \begin{recollection}
        Let $f \colon \tstr \to \tstr[Y] \in \StratN$. We say $f$ is a \define{poset-preserving diagrammatic (categorical) equivalence} if the underlying map of posets $P(f) \colon \ptstr \to \ptstr[Y]$ is an isomorphism and the induced stratum-preserving map $P(f)_! \tstr \to \tstr[Y]$ is a diagrammatic (categorical) equivalence. We denote the resulting relative categories on $\StratN$ by $\RelStratD$ and $\RelStratC$. This is how weak equivalence in $\StratN$ are defined, respectively, in \cite{douSimp,haine2018homotopy}. 
    \end{recollection}
    This definition of homotopy theories on $\StratN$ has the advantage that by making use of the Grothendieck bifibration $\StratN \to \Pos$, most questions about the global homotopy theories may be reduced to questions of the fibers. From a conceptual point of view, however, it has the side effect that the resulting homotopy theory contains a lot of highly pathological stratified spaces. Namely, both $\AltStratD$ and $\AltStratC$ contain a fully faithful copy of $\Pos$, given by empty stratified spaces. If we are looking to obtain a theory which embeds fully faithfully into the homotopy theory of (small) quasi-categories $\iCat$, i.e., fulfills a homotopy hypothesis closer to the classical one, then the stratification poset should generally not be a homotopy invariant, but merely an additional piece of data used to specify the allowable paths in a stratified space. This can be achieved by following the other possible approach to generalize the definitions in \cref{rec:douTheoryFixed,rec:Haine_theory_fixed}. Namely, one can transfer weak equivalences along a functor defined on stratified maps:
    \begin{construction}
          The extended homotopy link $\AltHolIP[](\tstr)$ of a stratified space $\tstr \in \StratN$ is the bisimplicial set given by 
        \[
        n \mapsto \Strat (\sReal{\Delta^{[n]}}, \tstr),
        \]
        with the obvious functoriality in $n$. Note that 
        \begin{align*}
            \Strat (\sReal{\Delta^{[n]}}, \tstr) &= \bigsqcup_{\I \in N(P)_n} \TopP[\ptstr]( \sReal{\Delta^\I}, \tstr)  \\ 
            &= \bigsqcup_{\I \in \sd \ptstr} \HolIP(\tstr) \sqcup \bigsqcup_{\J \in N(\ptstr)_n,  \textnormal{degenerate}} \TopP[\ptstr]( \sReal{\Delta^\J}, \tstr)  .
        \end{align*}
        In particular, $\AltHolIP[](\tstr)$ contains the data of all homotopy links of $\tstr$. Furthermore, whenever $\J$ is a flag that degenerates from a regular flag $\I$, then the degeneracy map induces weak equivalences $\HolIP(\tstr) \to \TopP[\ptstr](\sReal{\Delta^\J},\tstr)$. In this sense, if we remove homotopically redundant data, then $\AltHolIP[](\tstr)$ stores exactly the data of all homotopy links of $\tstr$. At the same time, allowing for degenerate simplices in homotopy links makes $\AltHolIP[]$ functorial in stratified maps, not only stratum-preserving ones. 
    \end{construction}
    \begin{definition}
        A stratified map $f \colon \tstr \to \tstr[Y]$ is called a \define{diagrammatic equivalence}, if it induces weak homotopy equivalences of simplicial sets on all extended homotopy links $\AltHolIP[n](\tstr) \to \AltHolIP[n](\tstr[Y])$, for $n \in \mathbb{N}$. We denote by $\RelStratDR$ the relative category given by $\StratN$ together with diagrammatic equivalences. 
    \end{definition}
    Finally, let us consider the categorical analogue of $\RelStratDR$.
    \begin{recollection}
        In \cite{nand2019simplicial}, Nand-Lal defined a stratified map $f \colon \tstr \to \tstr[Y]$ to be a \define{categorical equivalence}, if and only if the underlying simplicial map of 
        \[
        \SingS(f) \colon \SingS (\tstr) \to \SingS (\tstr[Y])
        \]
        is a categorical equivalence of simplicial sets (also called Joyal equivalences). We denote by $\RelStratCR$ the relative category given by $\StratN$ together with categorical equivalences. 
    \end{recollection}
    As we have already illustrated in the introduction, defining a homotopy theory of stratified spaces in terms of an $\infty$-categorical localization, especially in terms of maps which are themselves obtained through a localization, comes with an immediate series of questions. To just name a few: Are the stratum-preserving categorical equivalences detected by $\SingS$? Do we have any intrinsic description of the mapping spaces in the $\infty$-categorical localizations? More precisely, how do they relate to the classical mapping spaces of stratified maps given by the simplicial structure on $\Strat$? In the same vein of questioning, if we restrict ourselves to classical examples of stratified space, do we obtain the same homotopy theory as is investigated in \cite{quinn1988homotopically}? And if this is the case, how come \cite{quinn1988homotopically} has no need for higher homotopy links in his definition of a good class of stratified spaces? The classical way of making localizations of $1$-categories more tractable is, of course, the theory of model categories. In \cite{douteauEnTop}, the relative categories $\RelTopPD$ and $\RelStratD$, were extended to model structures. However, these model structures have only limited expressive power when it comes to investigating classical examples of stratified spaces, since these are almost never bifibrant (for example, no stratified cone on a non-empty manifold is ever bifibrant). As we have already illustrated in the introduction, the goal of this work is to construct (semi-)model structures for the homotopy theories defined in this section that make as many of such classical examples bifibrant.
    \subsection{Combinatorial models for stratified homotopy theory}\label{subsec:comb_mod_hostrat}
    One general strategy of constructing model categories is to first construct a model structure in a framework where this can easily be done, i.e., a category of presheaves or a similarly convenient setting, and then to transfer the theory along some adjunction. In this section, we will discuss simplicial analogues to the homotopy theories defined in the last section. We have investigated these theories in detail in \cite{ComModelWa}, connecting the theories of \cite{haine2018homotopy,douSimp} and constructing combinatorial models for $\AltStratCR$ and $\AltStratDR$.
    Let us begin with the case of $\AltTopPD$.
    \begin{recollection}[\cite{douSimp,douteauwaas2021}]\label{pshh:rec:Dou-Hen-mod}
    The \define{Douteau-Henriques model structure} on $\sSetP$, defined first in \cite{douSimp}, is the Cisinski model structure (see \cite[Thm. 2.4.19]{HigherCatCisinki}) induced by the simplicial cylinder $X \mapsto X \otimes \Delta^1$, with the empty set of anodyne extensions.
    % That is, it is the model-structure with the minimal amount of weak equivalences, in which stratified simplicial homotopy equivalences are weak equivalences and all monomorphisms are cofibrations.
    This defines a combinatorial, cofibrant, simplicial model structure on $\sSetP$, which may be characterized as the minimal model structure (in the sense of the smallest possible class of weak equivalences) in which the cofibrations are the monomorphisms and all stratified simplicial homotopy equivalences are weak equivalences. 
    We denote the resulting simplicial model category $\sSetPD$.
    $\sSetPD$ is cofibrantly generated by the classes of stratified boundary inclusions and admissible horn inclusions. It follows purely abstractly that weak equivalences between stratified simplicial sets $\tstr, \tstr[Y]$ that have the horn filling property with respect to admissible horn inclusions are precisely such stratum-preserving simplicial maps $\str \to \str[Y]$, for which the induced simplicial map of the simplicial homotopy links
    \[
    \HolIPS(\str) := \sSetP(\Delta^\I, \str) \to \sSetP( \Delta^\I, \str[Y]) =:  \HolIPS(\str[Y])
    \]
    is a weak homotopy equivalence, for all $\I \in \sd(\pos)$. More surprisingly, in \cite{douteauwaas2021} it was shown that this detection criterion holds for all $\tstr, \tstr[Y]$ and that no fibrancy assumptions are necessary.
    % We call such maps, inducing weak equivalences on all homotopy links, \define{diagrammatric equivalence}.
    Furthermore, mapping $\tstr \in \sStratN$ to the simplicial presheaf on $\sd(\pos)$ given by $\I \mapsto \Strat (\Delta^\I, \tstr$) induces a Quillen equivalence \[
\Diag \xrightleftharpoons{\simeq }\sSetPDN \colon \HolIPS[],
\]
where the left-hand side denotes the category of simplicial presheaves, $\FunC (\sd(\pos)^{\op}, \sSetN)$ equipped with the injective (or projective) model structure. Even more, this Quillen equivalence creates weak equivalences in both directions (see \cite[Thm. 1.3]{douteauwaas2021}).
In particular, $\sSetPD$ presents the $\infty$-category of space-valued diagrams indexed over $\sd(\pos)^{\op}$.
\end{recollection}
\begin{recollection}\cite{haine2018homotopy}\label{pshh:recol:haine_mod_cat}
    The \define{Joyal-Kan model structure} on $\sSetP$ is obtained by left Bousfield localizing $\sSetPD$ at the class of stratified inner horn inclusions. 
    % Weak equivalences in this model structure will be called \define{Joyal-Kan equivalences}. 
    Fibrant objects are precisely the stratified simplicial sets $\str$, for which the underlying simplicial set $\ustr$ is a quasi-category and $\sstr\colon \ustr  \to \pstr$ is a conservative functor. 
    It follows that $\sSetPC$ presents the $\infty$-category of conservative functors from a quasi-category into $\pos$, also called \define{abstract stratified homotopy types} over $\pos$. We can explicitly present the equivalence between d\'ecollages and abstract stratified homotopy types of \cite{haine2018homotopy} in terms of a Quillen equivalence. In \cite[Subsec. \ref{comb:subsec:decollages}]{ComModelWa}, we have shown that the adjunction
        \[
        \SDiag^{\dec} \xrightleftharpoons{}
        \sSetPC 
        \colon \HolIPS[]
        \]
    induces a Quillen equivalence between $\sSetPC$ and the model structure for d\'ecollages, which creates weak equivalences in both directions. This lifts the equivalence between abstract stratified homotopy types and d\'ecollages already proven in \cite{Exodromy} to the level of model categories.
\end{recollection}
\begin{recollection}\cite{ComModelWa}
Both the Douteau-Henriques as well as the Joyal-Kan model structures admit a global analogue on $\sStrat$ (named accordingly) obtained by gluing the model structures together along the Grothendieck bifibration $\sStratN \to \pos$, using \cite[Thm 4.4]{CagneMellies}. We denote by $\sStratD$ and $\sStratC$ the simplicial model categories with underlying category $\sStrat$, defined by applying \cite[Thm 4.4]{CagneMellies} to the forgetful functor \[\pos \colon \sStrat \to \Pos,\] with the fiberwise model structures given by $\sSetPD$ and $\sSetPC$, for $P \in \Pos$, respectively. 
In more detail, a stratified simplicial map $f \colon \tstr \to \tstr[Y]$ is a cofibration, if and only if the induced map $\pos (f)_! \tstr \to \tstr[Y]$ is a monomorphism, 
and weak equivalences are precisely such maps for which $\pos(f)$ is an isomorphism, and $\pos (f)_! \tstr \to \tstr[Y]$ is a weak equivalence in $\sSetPD[{\pos_{\str[Y]}}]$ ($\sSetPC[{\pos_{\str[Y]}}]
$).  $\sStratD$ and $\sStratC$ are cofibrant combinatorial model categories. A set of generating cofibrations for both model categories is given by the set of stratified boundary inclusions $\{ \stratBound  \hookrightarrow \stratSim \mid n \in \mathbb N\}$, together with the two morphisms
    \begin{diagram}
        \emptyset \arrow[r]  \arrow[d] & \emptyset \arrow[d] & & \emptyset  \arrow[r]  \arrow[d] & \emptyset \arrow[d] \\
        \emptyset \arrow[r] & {[0]}                           \spacecomma& &{[0] \sqcup [0]} \arrow[r, hook] & {[1]} \spaceperiod
    \end{diagram}
Furthermore, in $\sStratD$ acyclic cofibrations are generated by the stratified admissible horn inclusions 
 \begin{diagram}
        \Lambda_k^n \arrow[rr] \arrow[rd] & &\Delta^n \arrow[ld] \\
        & {[m]} &,
    \end{diagram}
    for $n,m \in \mathbb N$. 
% The model structure on $\sStratD$ is called the \define{Douteau-Henriques model structure} on $\sStrat$. The model structure on $\sStratC$, is called the \define{Joyal-Kan model structure} on $\sStrat$.
   % Weak equivalences in these model categories are called \define{diagrammatic equivalences} and \define{categorical equivalences} respectively.
\end{recollection}
Generally, the poset of a stratified simplicial set may contain a lot of redundant elements and relations that do not reflect in the simplicial set itself. 
For example, $\sStratC$ does not present a subcategory of $\iCat$, but instead an $\infty$-category of categories together with a conservative functor into a poset, and contains a fully faithful copy of the category of posets, given by mapping a poset $\pos$ to the functor $\emptyset \to \pos$. If we are looking to obtain a version of a stratified homotopy hypothesis, referring only to layered $\infty$-categories without a choice of additional data, a further (right) localization is necessary. This uses the fact that any layered $\infty$-category naturally comes with a functor to the poset of its isomorphism classes, with relations induced by morphisms between the latter.
\begin{recollection}[{\cite[Subsec. \ref{comb:subsec:refining_strat_sset}]{ComModelWa}}]\label{rec:ref_mod_structs}
    Given a stratified simplicial set $\tstr \in \StratN$, its refined poset $\rpstr$ is the poset generated with elements the vertices of $\tstr$, and a relation $x \leq y$, if and only if there is a path of $1$-simplices 
        \[ x = x_0 \leftrightarrow x_1 \leftrightarrow x_2 \leftrightarrow \cdots \leftrightarrow x_n = y \] 
    where only simplices that are contained within a stratum of $\str$ are allowed to point in direction of $x$. In other words, $\rpstr$ is the poset whose elements are the path component of strata of $\tstr$ and whose relations are generated by exit paths in $\str$.
    The stratification map of $\str$, $\sstr \colon X \to \ptstr$, factors through $\rpstr$, inducing a new stratification of $X$. The resulting stratified simplicial set is called \define{the refinement of $\tstr$} and denoted $\str^{\ared}$, and comes with a natural stratified simplicial map $\str^{\ared} \to \str$. Stratified simplicial sets for which this map is an isomorphism are called \textit{refined}.
    We denote by $\sStratDR$ ($\sStratCR$) the two (simplicial) right Bousfield localizations of respectively $\sStratD$ and $\sStratC$ obtained by localizing the refinement maps (the existence of these localizations was proven in  \cite[Thm. \ref{comb:prop:ex_red_struct}]{ComModelWa}). They are respectively called the \define{diagrammatic and categorical} model structure on $\sStrat$.
    % Weak equivalences in $\sStratDR$ are called \define{refined diagrammatic equivalences.}
    % Weak equivalences in $\sStratCR$ are called \define{refined categorical equivalences.}
    Both form a combinatorial model category, with cofibrations generated by the class of stratified boundary inclusions $\{ \stratBound  \hookrightarrow \stratSim \mid n \in \mathbb N\}$, together with the boundary inclusion $\partial \Delta^{[1]} \to \Delta^1$, into the trivially stratified simplex.
    The fibrant objects of $\sStratCR$, i.e., the refined abstract stratified homotopy types are precisely what \cite{Exodromy} call the \define{$0$-connected abstract stratified homotopy types} - such abstract stratified homotopy types which have $0$-connected strata and pairwise homotopy links. Forgetting the underlying stratification poset induces an equivalence 
    \[
    \AbStr^{\ared} \simeq \iCatO
    \]
    between the the homotopy theory of refined abstract stratified homotopy types $\AbStr^{\ared}$ and the homotopy theory of layered $\infty$-categories $\iCatO$ (see \cite{Exodromy} and \cite[Prop. \ref{comb:prop:sStraC_pres_Astrat} and Prop. \ref{comb:prop:Quillen_Equ_betw_ref_and_ord}]{ComModelWa}).
\end{recollection}
    \subsection{Equivalences of homotopy theories of stratified objects}\label{subsec:abstract_equ}
    Now that we have introduced both topological as well as simplicial homotopy theories of stratified spaces, the obvious question concerning the relation between the two arises. More precisely, what are the homotopy-theoretic properties of the adjunctions $\sReal{-} \dashv \SingS$?
    We prove that each homotopy theory of stratified topological spaces is, in fact, equivalent to its simplicial counterpart, with the equivalence induced through the adjunctions $\sReal{-} \dashv \SingS$.
     \begin{theorem}[Partially already found in \cite{douteauwaas2021,haine2018homotopy}]\label{thm:overview_over_all_hypothesis}
        Let $\textnormal{\textit{\textbf{R}}}$ be any of the relative categories of topological stratified spaces of \cref{subsec:top_homotopy_theories} and let $\mathrm{s}\textnormal{\textit{\textbf{R}}}$ be its simplicial counterpart from \cref{subsec:comb_mod_hostrat}. Then, the adjunction $\sReal{-} \dashv \SingS$ preserves all weak equivalences in both directions (i.e. induces functors of relative categories), and has unit and counit a weak equivalence (in other words, it is a strict homotopy equivalence of relative categories in the sense of \cite{BarwickRelative}).
        In particular, the adjunctions induce equivalences of quasi-categories
        \begin{align*}
         \AltsSetPD  &\simeq \AltTopPD;  \quad & \quad \AltsSetPC &\simeq \AltTopPC \spacecomma 
           \end{align*}
              for each $P \in \Pos$, as well as equivalences of quasi-categories
           \begin{align*}
           \AltsStratD  &\simeq \AltStratD ;  \quad & \quad  \AltsStratC &\simeq \AltStratC  ; \\
            \AltsStratDR  &\simeq \AltStratDR ;  \quad & \quad \AltsStratCR &\simeq \AltStratCR   .
        \end{align*}
     
    \end{theorem}
    Before we give a proof, let us remark on the history of these results.
    \begin{remark}
        The existence of an equivalence of quasi-categories $\AltsStratD \simeq \AltStratD$ was first shown in \cite{douteau2021stratified}, in terms of a Quillen equivalence. \cite{haine2018homotopy} used the equivalence between abstract stratified homotopy types and d\'ecollages, $\AltsSetPC \simeq \AltDiag^{\dec}$, proven in \cite{Exodromy}, as well as \cite[Thm. 3]{douteauEnTop}, to subsume that $\AltsSetPC$ was equivalent to a left localization of $\AltTopPD$. If one carefully follows the argument in \cite{haine2018homotopy}, this localization turns out to be precisely $\AltTopPC$. However, a priori, the resulting equivalence $\AltsSetPC \simeq \AltTopPC  $ is not induced by the adjunction $\sReal{-} \dashv \SingS$, but rather as a composition of equivalences $\AltTopPC \xrightarrow{\HolIP[]} \AltDiag^{\dec} \xrightarrow{\simeq} \AltsSetPC$. At this point, it was not yet known that $\SingS \colon \TopPN \to \sSetPN$ maps categorical equivalences to weak equivalences in the Joyal-Kan model structure. The strict homotopy equivalence of relative categories version for the cases $\AltStratD$ and $\AltTopPD$ was first shown in \cite[Sec. 5]{douteauwaas2021}, albeit it is stated there in slightly different language. As we will see, all other cases can be derived from this one quite readily.
    \end{remark}
    Let us first show that the nomenclature for categorical equivalences in $\TopPN$ is justified:
    \begin{lemma}\label{lem:char_of_cat_equ}
        Let $f\colon \tstr \to \tstr[Y] \in \TopPN$ be a stratum-preserving map. Then the following are equivalent:
        \begin{enumerate}
            \item $f$ is a categorical equivalence;
            \item $\SingS (f)$ is an equivalence in the Joyal-Kan model structure;
            \item The underlying simplicial map of $\SingS(f)$ is a categorical equivalence.
        \end{enumerate}
    \end{lemma}
    \begin{proof}
        There is a natural isomorphism of functors of $1$-categories $\HolIPS \circ \SingS \cong \HolIP$. By definition, $f$ is a categorical equivalence if and only if $\HolIP[](f)$ is an equivalence in the model structure for d\'ecollages. Hence, to prove that equivalence of the first two conditions holds, it suffices to see that $\HolIPS[]$ creates weak equivalences, which is part of the statement of  \cite[Thm. \ref{comb:prop:equ_decol_haine}]{ComModelWa}. The final equivalence is a consequence of \cite[Thm. 0.2.2]{haine2018homotopy}.
    \end{proof}
    As an immediate corollary, we obtain:
     \begin{corollary}\label{cor:char_of_cat_equ_glob}
          Let $f\colon \tstr \to \tstr[Y] \in \StratN$ be a stratified map. Then the following are equivalent:
        \begin{enumerate}
            \item $f$ is a poset-preserving categorical equivalence;
            \item $\SingS (f)$ is an equivalence in the poset-preserving Joyal-Kan model structure;
            \item $f$ is a categorical equivalence and induces an isomorphism on posets.
        \end{enumerate}
     \end{corollary}
    Furthermore, the following equivalence for diagrammatic equivalences holds:
    \begin{lemma}\label{lem:char_of_diag_equ_glob}
        A stratified map $f \colon \tstr \to \tstr[Y]$ in $\StratN$ is a poset-preserving diagrammatic equivalence if and only if $\SingS(f)$ is a weak equivalence in $\sStratDN$.
    \end{lemma}
    \begin{proof}
        This follows from \cite[Prop. \ref{comb:prop:equ_char_of_ref_diag_equ}]{ComModelWa} together with the natural isomorphism of $\AltHolIP[]$ with the composition of its simplicial counterpart with $\SingS$.
    \end{proof}
    Furthermore, we are going to make use of the following lemma:
    \begin{lemma}\label{lem:homo_equ_of_loc}
        Given two categories with weak equivalences $( \cat[C], W_{ \cat[C]})$, $( \cat[D], W_{ \cat[D]})$ suppose that a pair of functors
        \[
        L \colon   \cat[C] \rightleftharpoons  \cat[D] \colon R,
        \]\
        together with natural transformations $\varepsilon \colon 1_{ \cat[C]} \to RL$ and $\eta \colon LR \to 1_{ \cat[D]} $ define a strict homotopy equivalence of relative categories. 
        Let $W'_{ \cat[C]}$ be any wide subcategory such that $W_{\cat[C]} \subset W'_{\cat[C]}$,
        and such that $({\cat[C]}, W'_{\cat[C]})$ is again a category with weak equivalences. Denote $W'_{{\cat[D]}} = R^{-1}(W'_{ \cat[C]})$. Then the quadruple $(L,R,\varepsilon,\eta)$ also defines a strict homotopy equivalence between $( \cat[C], W'_{ \cat[C]})$ and $( \cat[D], W'_{ \cat[D]})$.
    \end{lemma}
   \begin{proof}
     Clearly, $W_{{\cat[D]}} \subset W'_{ \cat[D]}$ and $R(W'_{ \cat[D]}) \subset W'_{ \cat[C]}$.
      Really, the only thing one needs to verify is that $L$ remains a functor of relative categories, i.e., maps $W'_{\cat[C]}$ into $W'_{\cat[D]}$. So let $w \colon X \to Y \in W'_{\cat[C]}$. Consider the commutative diagram 
      % https://q.uiver.app/#q=WzAsNCxbMCwwLCJYIl0sWzAsMSwiWSJdLFsxLDAsIkxSKFgpIl0sWzEsMSwiTFIoWSkiXSxbMiwzLCJMUih3KSJdLFsxLDMsIlxcdmFyZXBzaWxvbl9ZIiwyXSxbMCwxLCJ3IiwyXSxbMCwyLCJcXHZhcmVwc2lsb25fWCJdXQ==
        \begin{diagram}
	X & {RL(X)} \\
	Y & {RL(Y)} \spaceperiod
	\arrow["{\varepsilon_X}", from=1-1, to=1-2]
	\arrow["w"', from=1-1, to=2-1]
	\arrow["{LR(w)}", from=1-2, to=2-2]
	\arrow["{\varepsilon_Y}"', from=2-1, to=2-2]
        \end{diagram}
        By the two-out-of-three property for weak equivalences, and the assumption that $\varepsilon_X, \varepsilon_Y \in W_{\cat[C]} \subset W_{\cat[C]}'$, it follows that $w$ being in $W'_{\cat[C]}$ implies that $RL(w) \in W'_{\cat[C]}$. As $R$ creates weak equivalences, this is the case if and only if $L(w)$ is in $W'_{\cat[D]}$. This finishes the proof.
   \end{proof}
   \begin{proof}[Proof of \cref{thm:overview_over_all_hypothesis}]
        Note that all model structures on the simplicial side are obtained in terms of Bousfield localizations from $\sSetPDN$ or $\sStratDN$ (see (see \cite[Prop. \ref{comb:prop:loc} and \ref{comb:prop:diag_of_bousfield_loc_sglob}]{ComModelWa}).
       It follows from \cref{lem:homo_equ_of_loc} together with \cref{lem:char_of_cat_equ,cor:char_of_cat_equ_glob,lem:char_of_diag_equ_glob} that it suffices to show the cases of $\RelTopPD$ and $\RelStratD$. These were proven in the proof of \cite[Thm. 5.1]{douteauwaas2021}. Specifically, see \cite[Lem 5.2 and 5.3]{douteauwaas2021} and their proofs.
    \end{proof}
    % As an immediate corollary of \cref{thm:overview_over_all_hypothesis,cor:equ_inj_mod_struct}, together with the Quillen equivalence between $\sSet$ and $\Top$ one obtains the following statement, which will play a decisive role in the proof of \cref{thm:cell_diag_are_homotopy_pushout}.
    % \begin{corollary}\label{cor:equ_of_top_diag}
    %     The topological homotopy link functor
    %     induces a homotopy equivalence of relative categories
    %     \[
    %     \HolIP[] : \RelTopPD \xrightarrow{\simeq} (\DiagTop, W_{\mathrm{inj}}),
    %     \]
    %     with $W_{\mathrm{inj}}$ the class of pointwise weak homotopy equivalenes. 
    % \end{corollary}

%% file: 4Modelstructures.tex
\section{(Semi-)model categories of stratified spaces}\label{sec:con_of_mod_cat}
The difficulty with the results in \cref{thm:overview_over_all_hypothesis} is their interpretability.  Without some additional structure on the side of $\StratN$, the only real information that we have is that the homotopy theories we defined are equivalent to the combinatorial world we transferred them from. 
% This might make the result seem almost tautotolical (even though the proofs in \cite{douteauwaas2021} show this is not the case). 
We thus encounter the difficulty of interpreting what these results actually say about classical examples of stratified spaces. That is, we need to determine how the localizations of \cref{thm:overview_over_all_hypothesis} relate to stratified maps, stratified homotopies or more generally the stratified simplicial mapping spaces of $\Strat$.
The usual approach to relating a localization to a simplicial category is the language of simplicial model categories. In \cite{douteauEnTop}, such simplicial model structures were suggested for $\AltStratD$ and $\AltTopPD$. The difficulty with these model structures, however, is that classically relevant examples such as (for example) Whitney stratified spaces are not bifibrant in the latter (in fact, no stratified cone on a manifold is cofibrant). Hence, the expressive power of these model structures when it comes to investigating mapping spaces for such spaces is limited. In \cite[A.5]{HigherAlgebra} or \cite[Thm. 8.1.2.6.]{nand2019simplicial}, it was shown that for essentially all stratified spaces of classical relevance $\tstr$ the stratified singular simplicial set $\SingS(\tstr)$ is fibrant in $\sStratCR$, and hence in any of the model structures for stratified simplicial sets. Furthermore, at least Whitney stratified spaces or PL-pseudo manifolds admit triangulations compatible with their stratifications (\cite{TriangulationsGoresky}), and are thus stratum-preserving homeomorphic to the realization of a stratified simplicial set. 
This suggests the approach of transferring the simplicial model structures on stratified simplicial sets to the topological setting. Constructing such a model structure in the case of $\AltStratCR$ was one of the stated goals and conjectures of \cite{nand2019simplicial}, which was ultimately left as a conjecture. In fact, it turns out that the model structures in the combinatorial setting do not transfer to the topological framework, which was first shown in \cite[A.5]{douteauwaas2021}.
We repeat the result here, in slightly different form and with a more detailed proof, as it illustrates well the difficulties which may arise in the topological but not in the combinatorial world.  
\begin{proposition}[{\cite[Prop. A.1, Rem. A.4, Rem. A.5]{douteauwaas2021}}]\label{prop:counter_example_mod}
There does not exist a model structure on $\StratN$, for which all poset-preserving diagrammatic equivalences are weak equivalences, and the inclusion 
$\sReal{\Lambda^{[3]}_1 \hookrightarrow \Delta^{[3]}}$
is a cofibration. The analogous result for $\TopPN$, for $P$ non-discrete, also holds.
\end{proposition}
\begin{proof}
    Consider the stratum-preserving map over $\{p < q\}$ illustrated in the following picture:
    \begin{center}
        \begin{tikzpicture}
        \coordinate(A) at (0, {sqrt(3)});
        \coordinate (B) at (-1, 0);
        \coordinate (C) at (1, 0);
        \coordinate (G) at (barycentric cs:A=1,B=1,C=1);
        \draw[color = blue, fill = blue, opacity = 0.6] (A) -- (B) -- (C) -- cycle;
        \coordinate (B') at (-2, 0);
        \coordinate (C') at (2, 0);
        \node[above] at (A){a};
        \node[below] at (B){b};
        \node[below] at (B'){b'};
        \node[below] at (C){c\textcolor{white}{'}};
        \node[below] at (C'){c'};
        \node[below] at (G){a'};
    \draw[color = blue, fill = blue, opacity = 0.6] (B') -- (B) -- (A) -- cycle;
    \draw[color = blue, fill = blue, opacity = 0.6]  (C') -- (C) -- (A) -- cycle;
    \foreach \x in {1, ..., 5}
        {\pgfmathsetmacro\wA{2*(2^\x -1)}  
        \coordinate (E1) at (barycentric cs:B=1,B'=1,A=\wA); 
            \node[fill=red, circle, inner sep = 0.4pt] at (E1){};
            }
    \foreach \x in {1, ..., 5}
        {\pgfmathsetmacro\wA{2*(2^\x -1)}  
        \coordinate (E1) at (barycentric cs:C=1,C'=1,A=\wA); 
            \node[fill=red, circle, inner sep = 0.4pt] at (E1){};
            }
    % \coordinate (E2) at (barycentric cs:B=1,B'=1,A=6);
    % \coordinate (E3) at (barycentric cs:B=1,B'=1,A=14);
    % \coordinate (E4) at (barycentric cs:B=1,B'=1,A=30);
    % \node[fill=red, circle, inner sep = 0.4pt] at (E1){};
    % \node[fill=red, circle, inner sep = 0.4pt] at (E2){};
    % \node[fill=red, circle, inner sep = 0.4pt] at (E3){};
    % \node[fill=red, circle, inner sep = 0.4pt] at (E4){};
    \draw[color = red] (G) -- (A);
    \draw[color=blue] (B) -- (G);
    \draw[color=blue] (C) -- (G);
    \node[fill=red, circle, inner sep = 0.6pt] at (A){};
    \node[fill=red, circle, inner sep = 0.6pt] at (G){};
    \begin{scope}[shift={(5,0)}]
        \coordinate(A) at (0, {sqrt(3)});
        \coordinate (B) at (-1, 0);
        \coordinate (C) at (1, 0);
        \draw[color = blue, fill = blue, opacity = 0.6] (A) -- (B) -- (C) -- cycle;
        \coordinate (B') at (-2, 0);
        \coordinate (C') at (2, 0);
         \node[above] at (A){a};
        \node[below] at (B){b};
        \node[below] at (B'){b'};
        \node[below] at (C){c\textcolor{white}{'}};
        \node[below] at (C'){c'};
    \draw[color = blue, fill = blue, opacity = 0.6] (B') -- (B) -- (A) -- cycle;
    \draw[color = blue, fill = blue, opacity = 0.6]  (C') -- (C) -- (A) -- cycle;
    \coordinate (G) at (barycentric cs:A=1,B=1,C=1);
    \foreach \x in {1, ..., 4}
        {\pgfmathsetmacro\wA{2*(2^\x -1)}  
        \coordinate (E\x) at (barycentric cs:B=1,B'=1,A=\wA); 
        \coordinate (L\x) at (barycentric cs:B=0.5,B'=1.5,A=\wA); 
            \node[fill=red, circle, inner sep = 0.4pt] at (E\x){};
            }
    \foreach \x in {1, ..., 4}
        {\pgfmathsetmacro\wA{2*(2^\x -1)}  
        \coordinate (F\x) at (barycentric cs:C=1,C'=1,A=\wA); 
            \node[fill=red, circle, inner sep = 0.4pt] at (F\x){};
        \coordinate (R\x) at (barycentric cs:C=0.5,C'=1.5,A=\wA); 
            }
        \draw[color = white] plot[smooth] coordinates {(B) (L1) (R2) (L3) (R4) ($(A)-(0, 0.05)$) (A)};
            \foreach \x in {1, ..., 4}
        { 
            \node[fill=red, circle, inner sep = 0.4pt] at (E\x){};
            }
    \foreach \x in {1, ..., 4}
        {  
            \node[fill=red, circle, inner sep = 0.4pt] at (F\x){};
            }
\end{scope}
    \draw[->] (2,1) -- (3,1) node[midway, above] {$r$};

    % \coordinate (E2) at (barycentric cs:B=1,B'=1,A=6);
    % \coordinate (E3) at (barycentric cs:B=1,B'=1,A=14);
    % \coordinate (E4) at (barycentric cs:B=1,B'=1,A=30);
    % \node[fill=red, circle, inner sep = 0.4pt] at (E1){};
    % \node[fill=red, circle, inner sep = 0.4pt] at (E2){};
    % \node[fill=red, circle, inner sep = 0.4pt] at (E3){};
    % \node[fill=red, circle, inner sep = 0.4pt] at (E4){};
    % \draw[color = red] (G) -- (A);
    % \draw[color=blue] (B) -- (G);
    % \draw[color=blue] (C) -- (G);
     \node[fill=red, circle, inner sep = 0.6pt] at (A){};
    % \node[fill=red, circle, inner sep = 0.6pt] at (G){};
    \end{tikzpicture}   
    \end{center}
    The picture shows two stratifications over $\{p<q\}$ of the realization of simplicial complexes (ignoring the white path from $a$ to $b$ for now), embedded in $\mathbb R^2$. Denote, from right to left, the corresponding stratified spaces by $\tstr$ and $\tstr[Y]$. The $p$-stratum is marked in red. In $\tstr$ the $p$-stratum consists of a vertical line and two sequences, given in barycentric coordinates with respect to $a,b,b'$ and $a,c,c'$ by $(1-\frac{1}{2^n}, \frac{1}{2^{n+1}},\frac{1}{2^{n+1}})$. In $\tstr[Y]$, the vertical line is replaced by a single point.
    Denote these sequences, respectively, by $b_n$ and $c_n$. The map $r$ is given by convexly extending the map which keeps all vertices besides $a'$ the same and maps the latter to $a$. This map is not a categorical equivalence. Indeed, consider the induced map of simplicial sets
    \[
    \SingS(\tstr) \to \SingS(\tstr[Y])
    \]
    and the induced map of hom-sets 
    \[
    \ho  \SingS(\tstr) (a,b) \xrightarrow{r_*} \ho  \SingS(\tstr[Y]) (a,b)
    \]
    in the associated homotopy categories. First, note that $\SingS(\tstr[Y])$ is, in fact, a quasi-category. This can be seen, for example, from the fact that every non-trivially stratified horn $\sReal{\Lambda^\J} \to \tstr$, $\J = [p \leq p \leq \cdots \leq q ]$ is constant on the $p$-stratum. It follows that it suffices to show that $\tstr[Y]$ admits extensions with respect to $ \sReal{\Lambda^\J/\Delta^{\J_p}} \hookrightarrow \sReal{\Delta^\J / \Delta^{\J_p}}$, which always admits a stratified retraction.
    Consequently, elements of $\SingS(\tstr[Y]) (a,b)$ are given by stratified homotopy classes (rel boundary) of exit paths from $a$ to $b$.
    It turns out that $r_*$ is not surjective: Consider a path $\gamma$ from $a$ to $b$ that (described as starting from $b$) ascends monotonously in height and passes to the left of $b_1$, to the right of $c_2$, to the left of $b_3$ etc., as illustrated in the right picture. The class of $\gamma$ does not lie in the image of $r_*$. Roughly speaking, this is due to the fact that no exit path in $\tstr$ can oscillate around the points of both sequences $(b_n)$ and $(c_n)$. For a formal proof, see \cref{appendix:proof_of_nonexistence_of_oath}. 
    Next, let $\J = [p \leq p < q \leq q]$ and $\J' = [p < q \leq q]$ consider the following diagram of pushout squares
    \begin{diagram}\label{diag:pushouts_in_mod_struct_counterexample}
        \sReal{\Lambda_1^{[3]}} \arrow[r,  hook] \arrow[d] &  \sReal{\Delta^{[3]}} \arrow[d] \arrow[r, "\sReal{s^1}"]& \sReal{\Delta^{[2]}} \arrow[d] \\ 
        \sReal{\Lambda^{\J}_1} \arrow[d, hook] \arrow[r, hook] & \sReal{\Delta^{\J}} \arrow[d, hook] \arrow[r] & \sReal{\Delta^{\I}} \arrow[d, hook] \\
        \tstr \arrow[r]& \tstr[Z] \arrow[r, "r'"] & \tstr[Y] \spacecomma
    \end{diagram}
    with the upper verticals induces by the identity on spaces and the flags $\J$ and $\J'$, the lower left vertical given by affinely extending the map of vertices $ e_0 \mapsto a'$, $e_1 \mapsto a$ and $e_2 \mapsto b$ and $e_3 \mapsto c$, and with $s^1$ the degeneracy map collapsing the edge $\{1,2\}$ to $\{1\}$.
    The bottom composition is precisely $r$, by construction of the latter. Now, the middle left horizontal is the realization of a diagrammatic equivalence and hence a diagrammatic equivalence. It follows that if a model structure as claimed in the statement of the proposition existed, then the left lower horizontal would be an acyclic cofibration.
    Now, suppose that the lower right vertical $r'$ was also a categorical equivalence. Then it would follow that $r$ is also a categorical equivalence, in contradiction to what we have just shown. In fact, the lower left vertical is even a stratum-preserving homotopy equivalence, and hence a categorical equivalence. To see this, note that $r'$ admits a section, induced by the pushout
    \begin{diagram}
        \sReal{\Delta^{\I}} \arrow[d] \arrow[r, hook, "\sReal{d^0}"] & \sReal{\Delta^\J} \arrow[d] \\
        \tstr[Y] \arrow[r, "s'"]& \tstr[Z] \spaceperiod
    \end{diagram}
    $\sReal{d_0}$ is even a stratum-preserving strong deformation retract (that is, admits a deformation retraction which is also stratum-preserving). It follows that $s'$ is also a stratum-preserving deformation retraction, and in particular a stratum-preserving homotopy equivalence. Since $r'$ is a retraction of $s'$, it follows that $r'$ is a stratified homotopy equivalence. For the case of a fixed poset $\pos$, construct $\tstr$ and $\tstr[Y]$ over any trivial relation $p <q$ in $\pos$, and simply repeat the argument omitting the first line in \cref{diag:pushouts_in_mod_struct_counterexample}.
\end{proof}
\subsection{A transfer lemma for cofibrantly generated semi-model structures}\label{subsec:transf_lemma}
% As we have illustrated in \cref{subsec:intro_to_strat_ho_hy}, if we want to understand the precise implications of the several equivalences of homotopy theories in \cref{thm:overview_over_all_hypothesis}, as well as the stratified homotopy hypotheses \cref{thm:hohy_haine_cat,thm:hohy_nand_cat}, better - in particular when it comes to classical examples of stratified spaces such as Whitney stratified spaces - it would be preferable to extend the relative categories describing the homotopy theories of stratified spaces to appropriate model categories, which make classical examples of stratified spaces cofibrant-fibrant.
% The obvious approach is to not only transfer the homotopy theories along the adjunctions $\sReal{-} \dashv \SingS$, but also the combinatorial model structures. This, however, turns out to be impossible in all cases, but the trivial one of $\TopP$ for $P$ a discrete poset.
If one takes a precise look at the counterexample in \cref{prop:counter_example_mod}, one will notice that it involves a pushout of an acyclic cofibration along a stratified map whose target is somewhat pathologically stratified. 
In particular, we cannot expect the target to be cofibrant. Hence, what one may nevertheless hope to obtain is a left semi-model structure. 
Semi-model structures provide a weaker version of model structures in which one can only expect pushouts of acyclic cofibrations to remain acyclic if all objects involved are cofibrant. 
Despite having slightly weaker axioms, for almost all intents and purposes, semi-model categories (cofibrantly generated ones to be more precise) are just as useful as regular model categories. Essentially, every theorem about model categories admits an analogue in the world of semi-model categories, as long as one takes care that sources of morphisms often need to be assumed to be cofibrant (see, for example, \cite{white2023left} for an overview of some results). We will not go through the ordeal of reproducing every necessary result from the model category world in the semi-model category world here. Instead, when necessary, we will often reference a proof in the model-categorical setting and explain what needs to be adapted in the setting of semi-model categories. Since we only ever make use of the cofibrantly generated scenario in this paper, and the latter tends to be significantly more well-behaved, we are going to use the following definition. Recall first the notion of weak factorization system (for example, from \cite[Ch. 11]{RiehlCategoricalHomotopyTheory})
\begin{definition}\label{def:left_semi}
    Let $\cat[M]$ be a bicomplete category. A cofibrantly generated left semi-model structure on $\cat[M]$ is the data of three classes $(C,W,F)$, called respectively cofibrations, weak equivalences, and fibrations, such that:
    \begin{enumerate}
        \item $W$ contains all isomorphisms and is closed under two-out-of-three and retracts;
        \item $(C, W \cap F)$ is a (functorial) weak factorization system;
        \item $F$ consists exactly of those morphisms which have the right lifting property with respect to all morphisms in $W \cap C$ with cofibrant source;
        \item Every morphism with cofibrant source factors (functorially) into a morphism in $W \cap C$, followed by a fibration. 
    \end{enumerate}
   Furthermore, we assume that there exist sets of cofibrations $I$ and of acyclic cofibrations with cofibrant source $J$, such that $F$ is the class of morphisms that have the right lifting property with respect to $J$, and $F \cap W$ is precisely the class of morphisms that have the right lifting property with respect to $I$. 
   A (cofibrantly generated, left) semi-model category is the data of a complete and cocomplete category $\textnormal{\textbf{C}}$, together with a (cofibrantly generated) model structure on $\textnormal{\textbf{C}}$. Since we are only concerned with left semi-model categories in this work, we will just speak of semi-model categories when we mean left semi-model category.
\end{definition}
It is a routine exercise to show that since we assumed cofibrant generation, our definition is equivalent to \cite[Def. 2.1]{white2023left}. Note that this definition is slightly stronger than the one used, for example, in \cite[12.1.1]{FresseOperads}. We use the following definition of simpliciality:
\begin{definition}\label{def:tranferred_model_structure}
    Let $\textbf{M}$ be a bi-complete simplicial category, admitting powers and tensors, with the underlying $1$-category equipped with a semi-model structure. We say that $\textbf{M}$ is a \define{simplicial semi-model} category, if for every cofibration $i \colon A \hookrightarrow B$ in the Kan-Quillen model structure on simplicial sets (that is, for every monomorphism), and for every fibration $f \colon X \to Y$ in $\textbf{M}$, the induced morphism
    \[
    X^B \to X^A \times_{Y^A} Y^B
    \]
    is a fibration, and furthermore a weak equivalence if $i$ or $f$ is a weak equivalence.
\end{definition}
\begin{remark}
   Note that this is a slightly different definition of simplicial semi-model category than the one in \cite[Def. 2.3]{white2023left} or \cite[Def. 1.1.8]{goerss2005moduli}. 
Our definition is a priori stronger than the ones in \cite{white2023left,goerss2005moduli}. The latter only implies that when $i$ is a cofibration and $f$ a fibration, then $Y^B \to X^B \times_{X^A} Y^A$ has the right lifting property with respect to all acyclic cofibrations with cofibrant source. However, since in our cases all relevant semi-model categories are generated by (acyclic) cofibrations with cofibrant source, for the semi-model categories we are concerned with, all of these definitions agree.
\end{remark}
\begin{definition}
Let $\textnormal{\textbf{D}}$ be a bicomplete category.
  Consider an adjunction
  \[
  L\colon  \textnormal{\textbf{C}} \rightleftharpoons \textnormal{\textbf{D}} \colon R 
  \]
and assume that \textnormal{\textbf{C}} is a model category. We say that the \define{right transferred semi-model structure on $\textnormal{\textbf{D}}$ along $R$ exists}, if the following classes form a semi-model structure on $\textnormal{\textbf{C}}$:
\begin{enumerate}
    \item Fibrations (weak equivalences) in $\textnormal{\textbf{D}}$ are precisely such $f \in \textnormal{\textbf{D}}$, for which $R(f)$ is a fibration (weak equivalence).
    \item Cofibrations in $\textnormal{\textbf{D}}$ are those morphisms that have the left lifting property with respect to all acyclic fibrations (i.e. morphisms which are both a weak equivalence and a fibration). 
\end{enumerate}
    We then call this structure the \define{right transferred semi-model structure} on $\textnormal{\textbf{D}}$ (with respect to $R$).
\end{definition}
Given a class of morphisms $I$ in a bicomplete category $\cat[C]$, we denote by $\Cell(I)$ the class of relative cell complexes with respect to $I$ (see \cite[Def. 10.5.8]{hirschhornModel}).
Similarly to the scenario of regular model categories (see, for example, \cite{nlab:transferred_model_structure} for a proof), one obtains:
\begin{lemma}\label{lem:ex_of_semi_gen}
    In the situation of \cref{def:tranferred_model_structure}, suppose that $\textnormal{\textbf{C}}$ is cofibrantly generated by cofibrations $I$ and acyclic cofibrations with cofibrant source $J$. Suppose furthermore that:
        \begin{enumerate}
        \item Every morphism in $\textnormal{\textbf{D}}$ factors into an element of $\Cell L(I)$, followed by an acyclic fibration.
        \item Every morphism in $\textnormal{\textbf{D}}$ with cofibrant source factors into an element of $\Cell L(J)$, followed by a fibration.
        \item Every element of $\Cell L(J)$ with cofibrant source is a weak equivalence.
    \end{enumerate}
    Then the right transferred model structure on $\textnormal{\textbf{D}}$ exists and is cofibrantly generated by $L(I)$ and $L(J)$.
\end{lemma}
Let us now consider the following transfer lemma for semi-model structures. We note that this result (omitting the simplicial part) can also be obtained as a special case of \cite[Thm. 2.2.2]{white2018bousfield} (see also \cite[Thm. 3.3]{fresse2010props}). We provide a proof here, purely in order for the reader not to have to translate back and forth between the slightly different setup in \cite{white2018bousfield}.
\begin{lemma}\label{lem:transf_lem_spec}
     Let $\textnormal{\textbf{D}}$ be a complete and cocomplete category.
  Consider an adjunction
  \[
  L\colon  \textnormal{\textbf{C}} \rightleftharpoons \textnormal{\textbf{D}} \colon R 
  \]
and assume that $\textnormal{\textbf{C}}$ is a cofibrantly generated model category, with cofibrations generated by a set $I$ and acyclic cofibrations generated by a set of morphisms with cofibrant source $J$, arrows in which have $\kappa$-small source, for some cardinal $\kappa$. Assume furthermore the following:
    \begin{TransfLem}
        \item {\label{lem:transf_ass_1}} $R$ preserves transfinite compositions of morphisms in $\Cell{L(I)}$ and $\Cell L(J)$ respectively;
        \item \label{lem:transf_ass_2} Weak equivalences in $\textnormal{\textbf{C}}$ are stable under transfinite composition;
        \item \label{lem:transf_ass_3} For any a pushout diagram in $\textnormal{\textbf{D}}$
            \begin{diagram}
                L(A) \arrow[r, "L(j)"] \arrow[d] & L(B) \arrow[d]\\
                X \arrow[r, "j'"] & Y
            \end{diagram}
        with $j \in J$ and $X \in \Cell(L(I))$ an absolute cell complex, the map $j'$ is a weak equivalence. 
    \end{TransfLem}
     Then the right transferred semi-model structure on $\textnormal{\textbf{D}}$ exists, and it is furthermore cofibrantly generated by $L(I)$ and $L(J)$. Even more, if $L \dashv R$ is an adjunction of simplicial categories, with $\cat[D]$ admitting a simplicial tensoring and powering, and the model structure on $\textnormal{\textbf{C}}$ is simplicial, then so is the induced semi-model structure on $\textnormal{\textbf{D}}$.
\end{lemma}
\begin{proof}
   It follows from \cref{lem:transf_ass_1} and the assumption on small sources that the classes $L(I)$ and $L(J)$ permit the small object argument. In fact, it follows under the adjunction $L \dashv R$ that the sources of $L(I)$ and $L(J)$ are $\kappa$ small, with respect to $\Cell L(I)$ and $\Cell L(J)$. Consequently, the only remaining requirement of \cref{lem:ex_of_semi_gen} is that every element $j' \in \Cell{L(J)}$ with cofibrant source is a weak equivalence.
  It follows by \cref{lem:transf_ass_1,lem:transf_ass_2} that we only need to consider the case of a pushout diagram 
    \begin{diagram}
                L(A) \arrow[r, "L(j)"] \arrow[d] & L(B) \arrow[d]\\
                X \arrow[r, "j'"] & Y
    \end{diagram}
   with $j \in J$ and $X$ cofibrant. Furthermore, since $X$ is cofibrant, it is the retract of an absolute $L(I)$ cell complex $X'$, through morphisms $X \xhookrightarrow{i} X'$ and $X' \xrightarrow{r} X$, with $r \circ i = 1_X$.
   We may extend the square in question to a diagram of pushout squares
  % https://q.uiver.app/#q=WzAsOCxbMCwwLCJMKEEpIl0sWzEsMCwiTChCKSJdLFswLDEsIlgiXSxbMCwyLCJYJyJdLFswLDMsIlgiXSxbMSwxLCJZIl0sWzEsMiwiWSciXSxbMSwzLCJZIl0sWzAsMV0sWzAsMl0sWzIsNCwiMV9YIiwyLHsiY3VydmUiOjJ9XSxbMiwzLCJpIl0sWzMsNCwiciJdLFsyLDUsImonIl0sWzEsNV0sWzUsNl0sWzMsNiwiaicnIl0sWzQsNywiaiciXSxbNiw3XSxbNSw3LCIxX1kiLDAseyJjdXJ2ZSI6LTJ9XV0=
\begin{diagram}
	{L(A)} & {L(B)} \\
	X & Y \\
	{X'} & {Y'} \\
	X & Y \spaceperiod
	\arrow[from=1-1, to=1-2]
	\arrow[from=1-1, to=2-1]
	\arrow[from=1-2, to=2-2]
	\arrow["{j'}", from=2-1, to=2-2]
	\arrow["i", from=2-1, to=3-1]
	\arrow["{1_X}"', curve={height=12pt}, from=2-1, to=4-1]
	\arrow[from=2-2, to=3-2]
	\arrow["{1_Y}", curve={height=-12pt}, from=2-2, to=4-2]
	\arrow["{j''}", from=3-1, to=3-2]
	\arrow["r", from=3-1, to=4-1]
	\arrow[from=3-2, to=4-2]
	\arrow["{j'}", from=4-1, to=4-2]
\end{diagram}
    In particular $j'$ is a retract of $j''$, and it suffices to show that the latter is a weak equivalence.
    By the composability of pushout squares, it follows that we have thus reduced to the case where $X$ is an absolute $L(I)$-cell complex. This case is precisely the content of \cref{lem:transf_ass_3}. Simpliciality follows immediately by using that 
   \[
   R( X^B \to X^A \times_{Y^A} Y^B) \cong \big ( R(X)^B\to R(X)^A \times_{R(Y)^A} R(Y)^B \big ),
   \]
   and the explicit definition of (acyclic) fibrations in $\textnormal{\textbf{D}}$ in terms of $R$.
\end{proof}
\subsection{Semi-model structures on \texorpdfstring{$\Strat$}{Strat}}\label{subsec:proof_ex_top_modstruct}
To transfer the model structures from the combinatorial world to $\TopP$ and $\Strat$, we now need to verify the assumptions of \cref{lem:transf_lem_spec}. Most of these will turn out to be fairly standard abstract arguments. The difficult part is showing \cref{lem:transf_ass_3}. This requires an in-depth study of the behavior of homotopy links of absolute cell complexes, with respect to realizations of stratified boundary inclusions, which was performed in \cite{HoLinksWa}.
\begin{definition}
    Denote $I:=\{ \sReal{ \partial \Delta^\J \hookrightarrow \Delta^\J} \mid \J \in \Delta_P \} \subset \TopPN$. Absolute cell complexes with respect to $I$ are called \define{cellularly stratified spaces.} 
    If $\tstr[A]$ and $\tstr$ are cellularly stratified spaces, then we call a stratum-preserving map $f \colon \tstr[A] \hookrightarrow \tstr$ which is a relative $I$-cell complex an \define{inclusion of cellularly stratified spaces}. 
\end{definition}
The following result is the decisive argument to apply \cref{lem:transf_lem_spec} to transfer the model structure on combinatorially stratified objects to the stratified framework. 
It is a corollary of \cite[Thm. \ref{hol:thm:hol_main_result}]{HoLinksWa}, which is one of the main result of \cite{HoLinksWa}.
\begin{lemma}\label{lem:pushout_of_cell}
    Consider a pushout diagram in $\StratN$, 
    \begin{diagram}
        \sReal{\tstr[A]} \arrow[r, "\sReal{i}"] \arrow[d]&  \sReal{\tstr[B]} \arrow[d]\\
         \tstr[X] \arrow[r, "i'"]& {\tstr[Y]} 
    \end{diagram}
    with $\tstr[X]$ cellularly stratified. Suppose that $i$ is an acyclic cofibration in $\sStratDN$. Then, $i'$ is a weak equivalence in $\RelStratD$, i.e., an isomorphism in $\AltStratD$. The analogous result holds for any of the model categories of stratified simplicial sets and their corresponding topological $\infty$-category in \cref{not:all_top_strat_cat}.
\end{lemma}
Before we give a proof, note that \cref{lem:pushout_of_cell} can be taken as a statement about certain pushout diagrams of stratified spaces being homotopy pushout, i.e. pushout in the associated $\infty$-category obtained by inverting weak equivalences. This already shows that the existence of semi-model structures is closely related with \cref{iQ:ComputationHocolim}.
\begin{proof}
    Let us begin with the case over a fixed poset $P$, and all morphisms on posets given by the identity. 
    We claim that the diagram in the statement of the proposition is a pushout diagram in the quasi-categories $\AltTopPD$ and $\AltTopPC$. We already know that $\AltTopPC$ is a left-localization of $\AltTopPD$ (by \cref{thm:overview_over_all_hypothesis} and the fact that $\sSetPC$ is a left Bousfield localization of $\sSetPD$, together with \cite[Prop. 5.2.4.6]{HigherTopos}) and in particular that the canonical functor $\AltTopPD \to \AltTopPC$ is left adjoint. Therefore, by \cite[5.2.3.5]{HigherTopos}, which states that left-adjoint functors preserve colimits, it suffices to show the case of $\AltTopPD$. Recall from \cref{pshh:rec:Dou-Hen-mod} that taking homotopy links induces equivalences between $\AltTopPD$ and $\FunC ( \sd(\pos)^{\op}, \iSpaces)$. Thus, it suffices to show that the diagram becomes pushout after applying homotopy links. Now, by \cite[Cor. 5.1.2.3]{HigherTopos} a diagram in a functor quasi-category is pushout if and only if it is so after evaluating at every point. Thus, it suffices to show that the diagram in $\iSpaces$
    \begin{diagram}\label{diag:important_pushout_hol}
        \HolIP\sReal{\tstr[A]} \arrow[r ] \arrow[d]&   \HolIP \sReal{\tstr[B]} \arrow[d]\\
          \HolIP \tstr[X] \arrow[r]&  \HolIP {\tstr[Y]} 
    \end{diagram}
    is pushout. Now, the quasi-category $\iSpaces$ is given by the homotopy coherent nerve of $\Top^{o}$. \cite[Thm. \ref{hol:thm:hol_main_res_B}]{HoLinksWa} states that the diagram in $\Top$ corresponding to \cref{diag:important_pushout_hol} is a homotopy pushout square. Finally, by \cite[Thm. 4.2.4.1.]{HigherTopos}, this is equivalent to \cref{diag:important_pushout_hol} being a pushout diagram in $\iSpaces$. \\
    We have established that the diagram in the statement of the lemma is pushout in $\AltTopPD$ and $\AltTopPC$. Consequently, by \cref{thm:overview_over_all_hypothesis}, so is its image under $\SingS$ in $\AltsSetPD$ ($\AltsSetPC$). By assumption, $i$ defines an isomorphism in $\AltsSetPD$ ($\AltsSetPC$). It follows, again by \cref{thm:overview_over_all_hypothesis}, that $\SingS (\sReal{i})$ also defines an isomorphism in $\AltsSetPD$ ($\AltsSetPC$). Consequently, since pushouts preserve isomorphisms, it follows that $\SingS(i')$ is an isomorphism in the quasi-category $\AltsSetPD$ ($\AltsSetPC$), and thus a weak equivalence in $\sSetPDN$ ($\sSetPCN$) (by \cite[Thm. 8.3.10.]{hirschhornModel} and the fact that the homotopy category functor from quasi-categories to 1-categories commutes with localization). It follows that $i$ is a weak equivalence, as was to be shown. \\
    % To see that $i'$ is a diagrammatic equivalence, we need to see that for any regular flag $\I \in \sd(\pos)$ the induced map $\HolIP(i')$ is a weak homotopy equivalence.
    % Consider the induced diagram in $\Fun ( \sd(\pos), \TopN)$ 
    % \begin{diagram}
    %     \HolIP[] \sReal{\tstr[A]} \arrow[r] & \HolIP[] \\
    %     \HolIP \tstr & 
    % \end{diagram}
    % \cref{thm:cell_diag_are_homotopy_pushout} states that such a diagram is homotopy pushout in $\AltTopPD$, i.e. pushout in the $\infty$-categorical sense. 
    % % Furthermore, by \cref{thm:overview_over_all_hypothesis}, $\sReal{i}$ is a diagrammatic equivalence. \\
    % Hence, the result follows, since a morphism is a weak equivalence in $\RelTopPD$, if any only if it is an isomorphism in $\AltTopPD$, and isomorphisms are stable under pushout. 
    % Note, that the same argument applies to $\AltTopPC$, since by \cref{prop:loc}, the functor of $\infty$-categories $\AltTopPD \to \AltTopPC$ is left adjoint, and hence preserves pushout diagrams. 
    Next, let us cover the poset-preserving cases in $\sStratN$. Note that $i$ induces an isomorphism on posets.
    By pushing forward with $\pstr[A] \to \ptstr[X]$, and using the fact that this preserves realizations of acyclic cofibrations (by \cref{prop:model_global_mod_struct_are_glued}), we may assume that the whole diagram lies in $\TopPN[{\ptstr[X]}]$. Hence, this case follows from the previous cases.
    Finally, for the non-poset preserving cases, note that the model structure $\sStratDRN$ ($\sStratCRN$) on $\sStratN$ has the same acyclic cofibrations as $\sStratDN$ ($\sStratCN$) as it is obtained via right Bousfield localization. Hence, this case follows from the case $\sStratDN$ ($\sStratCN$).
\end{proof}
We may now prove the following.
\begin{corollary}\label{cor:semi_mod_transf}
    Let $\underline{\textnormal{\textbf{cS}}}$ be any of the simplicial model categories of stratified simplicial sets listed in \cref{subsec:comb_mod_hostrat}. Let $\underline{\textnormal{\textbf{S}}}$ be correspondingly the simplicial category of topological stratified spaces $\TopP$ or $\Strat$. Then, along the adjunction
    \[
    \sReal{-} \colon \underline{\textnormal{\textbf{cS}}} \rightleftharpoons \underline{\textnormal{\textbf{S}}} \colon \SingS ,
    \]
    the simplicial model structure on $\underline{\textnormal{\textbf{cS}}}$ right-transfers to a cofibrantly generated simplicial semi-model structure on $\underline{\textnormal{\textbf{S}}}$.
\end{corollary}
\begin{proof}
    We verify the assumptions of \cref{lem:transf_lem_spec}. Note first that $\underline{\textnormal{\textbf{cS}}}$
    is a cofibrantly generated simplicial model category (with acyclic generators with cofibrant source) as we have recalled in \cref{subsec:comb_mod_hostrat}.
    \cref{lem:transf_ass_1} follows from the classical fact that any compactum in a relative cell complex only intersects finitely many cells (see, for example, \cite[Prop. 10.7.4]{hirschhornModel}). \cref{lem:transf_ass_2} was verified in \cite[Prop. \ref{comb:prop:we_stable_colim},\ref{comb:prop:strat_we_stable_colim} and \ref{comb:prop:we_stable_under_colim_ref}]{ComModelWa}. 
    Finally, \cref{lem:transf_ass_3} follows from \cref{lem:pushout_of_cell}.
\end{proof}
\begin{remark}\label{rem:extending_nand_struct}
    We note that the analogous proof of \cref{cor:semi_mod_transf} also applies to the category of surjectively stratified spaces, as investigated in \cite{nand2019simplicial}, $\Strat^s$, and one may replace $\sStrat$ by $\sSetN$ with the Joyal-model structure (omitting the simpliciality statement). 
    Indeed, \cref{lem:pushout_of_cell} also applies to surjectively stratified spaces and realizations of Joyal acyclic cofibrations, using \cref{prop:eq_of_sur_strat} and that surjectively stratified spaces form a coreflective subcategory of $\Strat$.
    In particular, this also gives an affirmative answer to the semi-model structure conjectured in \cite[Subsec. 8.4.1]{nand2019simplicial}. As a consequence of \cref{prop:eq_of_sur_strat}, the resulting semi-model category is Quillen equivalent to the one on all stratified spaces $\Strat$ transferred from $\sStratCR$.
    We also note that the alternative method suggested for the proof in \cite{nand2019simplicial} cannot succeed, as \cite[Conjecture 2]{nand2019simplicial} is false (realizations of fibrant simplicial sets are generally not fibrant). We will provide a counterexample in future work.
\end{remark}
\begin{notation}\label{not:all_top_strat_cat}
    We denote by $\TopPD$ ($\TopPC$), $\StratD$ ($\StratC$) and $\StratDR$ ($\StratCR$), respectively, the corresponding simplicial semi-model categories whose existence is guaranteed by \cref{cor:semi_mod_transf}. We, respectively, call the corresponding semi-model structures the \define{diagrammatic (categorical), poset-preserving diagrammatic (poset-preserving categorical) and diagrammatic (categorical)} model structures.
\end{notation}
We use the remainder of the section to make \cref{cor:semi_mod_transf} more explicit by restricting it to its special cases. From \cref{pshh:rec:Dou-Hen-mod,pshh:recol:haine_mod_cat,cor:semi_mod_transf} we obtain:
\begin{theorem}\label{thm:ex_of_mod_struct_local}
    Let $P \in \Pos$. The simplicial category $\TopP$ admits the structures of cofibrantly generated simplicial semi-model categories, $\TopPD$ and $\TopPC$ - called the diagrammatic and categorical model structure - defined by the following classes:
        \begin{enumerate}
            \item Cofibrations are generated by the set of stratified boundary inclusions \[ \{ \sReal{\partial \Delta^\J \hookrightarrow \Delta^\J} \mid \J \in \Delta_P\}.\] 
            \item Weak equivalences are the stratum-preserving diagrammatic equivalences or respectively the stratum-preserving categorical equivalences of $P$-stratified spaces.
            \item Fibrations are the stratum-preserving maps that have the right lifting property with respect to all acyclic cofibrations with cofibrant source. In $\TopPD$ they are equivalently characterized by having the right lifting property with respect to all realizations of admissible stratified horn inclusions $\sReal{\Lambda^\J_k \hookrightarrow \Delta^\J}$, $\J \in \catFlag$.
        \end{enumerate}
\end{theorem}
% \Lukas{edit appropriately}
% \begin{proposition}\label{prop:cof_and_fib_top_over_p} Cofibrant and fibrant objects in $\TopPD$ and $\TopPC$, as well as fibrations between them, can be characterized as follows.
% \begin{itemize}
%     \item The cofibrant objects in $\TopPD$ and $\TopPC$ are precisely the retracts of celullarly stratified spaces.
%     \item The fibrant objects of $\TopPD$ ($\TopPC$) are precisely such stratified spaces which have the horn filling property with respect to realizations of admissible stratified horn inclusions (inner stratified horn inclusions).
%     % \item The fibrant objects of $\TopPC$ are precisely such stratified spaces which have the horn filling property with respect to realizations of inner stratified horn inclusions.
%     \item Fibrations in $\TopPD$ are equivalently characterized by having the right lifting property with respect to all realizations of admissible stratified horn inclusions.
%     \item The fibrations between fibrant objects in $\TopPC$ are precisely the stratified maps, which have the right lifting property with respect to to realizations of inner stratified horn inclusions.
% \end{itemize}
% \end{proposition}
Next, let us consider the global version of $\TopP$, allowing for varying stratification posets.
From \cref{pshh:rec:Dou-Hen-mod,cor:semi_mod_transf} we obtain the following two results:
\begin{theorem}\label{thm:ex_of_model_struct_strat_nonref}
     The simplicial category $\Strat$ admits the structures of cofibrantly generated simplicial semi-model categories, $\StratD$ and $\StratC$  - called the poset-preserving diagrammatic and poset-preserving categorical model structure - with the following classes:
        \begin{enumerate}
            \item Cofibrations are generated by the set of stratified boundary inclusions \[ \{ \sReal{\stratBound \hookrightarrow \stratSim}  \mid n \in \mathbb N\},\] together with the maps of empty stratified spaces over $[0] \sqcup [0] \hookrightarrow \{0<1\}$ and $\emptyset \hookrightarrow [0]$.
            \item Weak equivalences are the poset-preserving diagrammatic equivalences or, respectively, the poset-preserving categorical equivalences of stratified spaces.
            \item Fibrations are the stratum-preserving maps that have the right lifting property with respect to all acyclic cofibrations with cofibrant source. 
            In $\StratD$ they are equivalently characterized by having the right lifting property with respect to all realizations of admissible stratified horn inclusions $\sReal{\Lambda^\J_k \hookrightarrow \Delta^\J}$, $\J \in \Delta_{[m]}$, for $m \in \mathbb N$.
        \end{enumerate}
 \end{theorem}
 Similarly, to the case of stratified simplicial sets, the global diagrammatic and categorical model structures can essentially be understood entirely from the ones over a fixed poset and the functoriality of base changes along the poset. This is due to the fact that they may again be interpreted as being pieced together from the local pieces along a cartesian bifibration:
 \begin{proposition}\label{prop:model_global_mod_struct_are_glued}
      The model structure on $\StratDN$ ( $\StratCN$ ) is the unique semi-model structure on $\StratN$ pieced together (in the sense of \cite{CagneMellies}, using notation from there)\footnote{ \cite{CagneMellies} covers the case of model structures. For a generalization to semi-model structures see \cite{batanin2023model}.} from the diagrammatic (categorical) model structures on the fibers $\TopPN$ of the Grothendieck bifibration
      \[
      \pos \colon \StratN \to \Pos.
      \]
      In particular, a stratified map $f \colon \tstr \to \tstr[Y]$ is a fibration if and only if the canonical map $\baseBack \colon  \tstr \to \pos(f)^*{\tstr[Y]}$
      is a fibration, and $f$ is a cofibration if and only if $\baseForw \colon \pos(f)_!\tstr \to \tstr[Y]$ is a cofibration.
 \end{proposition}
 \begin{proof}
    We cover the case of $\StratDN$, the case of $\StratCN$ is completely analogous.
    Consider the forgetful functor $\pos \colon \StratN \to \Pos$, mapping a stratified space $\tstr$ to its stratification poset $\ptstr$. It is a Grothendieck bifibration (\cite{CagneMellies}), with the left action given by $f \mapsto f_!$ and the right action given by $g \mapsto g^*$. If we equip $\Pos$ with the trivial model structure, with all morphisms bifibrations and weak equivalences given by isomorphisms, then $\pos \colon \StratN \to \Pos$ tautologically forms a Quillen bifibration in the sense of \cite{CagneMellies} (replacing model categories by cofibrantly generated semi-model categories).
    We obtain two commutative diagrams (one for each horizontal) with diagonals Quillen bifibrations
    \begin{diagram}
        \StratDN \arrow[rr, -left to  , yshift= 0.225ex, "\simeq"] \arrow[rr, left to-,yshift=-0.225ex] \arrow[rd] & & \arrow[ld] \sStratDN \\
        & \Pos & \spacecomma
    \end{diagram}
    where the right adjoint horizontal preserves the right action, and the left adjoint the left actions. The fibers of the right hand vertical are the model categories $\sSetPDN$ (see the construction of the global model categories in  \cite[Subsec.\ref{comb:subsec:from_local_to_global}]{ComModelWa}). Since the semi-model structure on $\StratDN$ is transferred along the horizontal adjunction, it follows that its restriction to the fibers of the left vertical $\TopPN$, are precisely the transfers of the model structure on $\sSetPDN$ to $\TopPN$, which is $\TopPDN$. Using the obvious cofibrantly generated semi-model categorical version of \cite[Prop. 3.4,3.5]{CagneMellies} (see also \cite{batanin2023model}), which is proven identically, the characterization of cofibrations and fibrations follows.
 \end{proof}
% \begin{proposition}\label{prop:cof_fib_nonref_top}
%     The global analogue of \cref{prop:cof_and_fib_top_over_p} holds for $\StratD$ ($\StratC$) mutatis mutandis. For the characterizations of (fibrations and) fibrant objects, the result holds with the realizations of the following restricted classes of horn inclusions:
%     \begin{itemize}
%         \item The set of admissible horn inclusions $\Lambda_k^\J \hookrightarrow \Delta^\J$ over the poset $[m]$, where $m$ ranges over $\mathbb N$, for $\sStratD$.
%         \item The set of inner horn inclusions $\stratHorn \hookrightarrow \stratSim$, for $\sStratC$.
%     \end{itemize}
% \end{proposition}
% Note, that there is no need to also consider admissible horn inclusions for checking fibrancy in $\StratC$, as the strata of $\SingS \tstr$ are automatically Kan-complexes.
Finally, let us state the explicit results for the refined framework.
We also are going to need the following construction.
\begin{construction}\label{con:unstrat_bound_inc_is_retract}
    Consider the pushout diagram in $\StratN$
    \begin{diagram}
       \sReal{\stratBound[1]} \arrow[r, hook] \arrow[d, hook, "t" ]& \sReal{\stratSim[1]} \arrow[d, hook]\\
       \sReal{\stratSim[1]} \arrow[r, hook]& S^1 \spacecomma
    \end{diagram}
    where the vertical $t$ is given by the opposite inclusion of the boundary points of the stratified interval. The pushout is a trivially stratified $S^1$. In particular, the diagonal of this diagram - which is the inclusion of two points lying over separate strata into a trivially stratified $S^1$ - is a cofibration in the refined model structures. 
    The realization of the stratified boundary inclusion $\stratBound[1] \hookrightarrow \Delta^1$, where $\Delta^1$ is trivially stratified, is a retract of the diagonal of this diagram. 
\end{construction}
The following are a consequence of \cite[Thm. \ref{comb:prop:ex_red_struct}]{ComModelWa}, which we recalled in \cref{rec:ref_mod_structs}, and \cref{cor:semi_mod_transf}. Note that, by \cref{con:unstrat_bound_inc_is_retract}, the topological framework needs one less generator of cofibrations than the combinatorial framework.
\begin{theorem}\label{thm:ex_mod_struct_red}
     The simplicial category $\Strat$ admits the structures of cofibrantly generated simplicial semi-model categories, $\StratDR$ and $\StratCR$  - called the diagrammatic and categorical model structures - with the following classes:
        \begin{enumerate}
            \item Cofibrations are generated by the set of stratified boundary inclusions \[ \{ \sReal{\stratBound \hookrightarrow \stratSim}  \mid n \in \mathbb N\}.\]
            \item Weak equivalences are the diagrammatic equivalences or respectively the categorical equivalences of stratified spaces.
            \item Fibrations are the stratum-preserving maps that have the right lifting property with respect to all acyclic cofibrations with cofibrant source. 
            In $\StratDR$, they are equivalently characterized by having the right lifting property with respect to all realizations of admissible stratified horn inclusions $\sReal{\Lambda^\J_k \hookrightarrow \Delta^\J}$, $\J \in \Delta_{[m]}$, for all $m \in \mathbb N$.
        \end{enumerate}
\end{theorem}
% \begin{proposition}\label{prop:Grothendieck_bifib}
%     \Lukas{Add a proposition on being the induced model structure with respect to a Grothendieck-bifibration}
% \end{proposition}
Finally, we may summarize the connections between the various semi-model categories of stratified objects in the following result, which follows by \cite[Prop. \ref{comb:prop:loc},\ref{comb:prop:c_is_left_bous_of_d},\ref{comb:prop:c_is_left_bous_of_d_ref},\ref{comb:prop:ex_red_struct},\ref{comb:prop:diag_of_bousfield_loc_sglob}]{ComModelWa} and \cref{prop:ref_com_sings,thm:ex_mod_struct_red,thm:overview_over_all_hypothesis}. Here, by a Quillen equivalence of (cofibrantly generated left) semi-model categories, we mean an adjunction fulfilling any (or equivalently all) of the classical characterizations for model categories (see, for example, \cite{nlab:quillen_equivalence}). In the definition of a Quillen adjunction of cofibrantly generated left semi-model categories (see \cite{nlab:quillen_adjunction}), one only needs to take care that the left adjoint generally only is required to preserve acyclic cofibrations with cofibrant source.
\begin{theorem}\label{thm:huge_overview}
    The various simplicial semi-model categories discussed in this section fit into a diagram
    \begin{diagram}
        \sSetPD  \arrow[rd, -left to  ,yshift= 0.16ex, xshift= 0.16ex, crossing over] \arrow[rd, left to-,yshift=-0.16ex,  xshift= -0.16ex]
        \arrow[dd, hook] \arrow[rr, -left to  , yshift= 0.225ex, "\simeq"] \arrow[rr, left to-,yshift=-0.225ex]&      &\TopPD \arrow[dd, hook]\arrow[rd, -left to  ,yshift= 0.16ex, xshift= 0.16ex, crossing over] \arrow[rd, left to-,yshift=-0.16ex,  xshift= -0.16ex] & \\
              &\sSetPC \arrow[rr, -left to  ,yshift= 0.225ex, "\simeq" xshift= 10pt, crossing over] \arrow[rr, left to-,yshift=-0.225ex] &      & \TopPC \arrow[dd, hook]\\
        \sStratD \arrow[dd, -left to  , xshift= 0.225ex] \arrow[dd, left to-,xshift=-0.225ex] \arrow[rd, -left to  ,yshift= 0.16ex, xshift= 0.16ex, crossing over] \arrow[rd, left to-,yshift=-0.16ex,  xshift= -0.16ex] \arrow[rr, -left to  ,yshift= 0.225ex,  "\simeq" xshift= 20pt] \arrow[rr, left to-,yshift=-0.225ex]&      &\StratD \arrow[dd, -left to  , xshift= 0.225ex, "\simeq"] \arrow[dd, left to-,xshift=-0.225ex]  \arrow[rd, -left to  ,yshift= 0.16ex, xshift= 0.16ex, crossing over] \arrow[rd, left to-,yshift=-0.16ex,  xshift= -0.16ex] & \\
              &\sStratC \arrow[rr, -left to  ,yshift= 0.225ex, "\simeq" xshift= - 10pt, crossing over] \arrow[rr, left to-,yshift=-0.225ex] &      & \StratC \arrow[dd, -left to  , xshift= 0.225ex] \arrow[dd, left to-,xshift=-0.225ex] \\
        \sStratDR \arrow[rd, -left to  ,yshift= 0.16ex, xshift= 0.16ex, crossing over] \arrow[rd, left to-,yshift=-0.16ex,  xshift= -0.16ex] \arrow[rr, -left to  ,yshift= 0.225ex,  "\simeq" xshift= 20pt] \arrow[rr, left to-,yshift=-0.225ex]&      &\StratDR \arrow[rd, -left to  ,yshift= 0.16ex, xshift= 0.16ex, crossing over] \arrow[rd, left to-,yshift=-0.16ex,  xshift= -0.16ex] & \\
             &\sStratCR \arrow[rr, -left to  ,yshift= 0.225ex, "\simeq" xshift=-10pt] \arrow[rr, left to-,yshift=-0.225ex]&      & \StratCR 
        \arrow[from = 2-2, to = 4-2 , crossing over, hook]
        \arrow[from = 4-2, to = 6-2 , crossing over, -left to, xshift= 0.225ex]
          \arrow[from = 4-2, to = 6-2 , left to-, xshift=-0.225ex]
    \end{diagram}
    of simplicial functors, with all horizontals induced by the adjunction $\sReal{-} \dashv \SingS$, and all remaining functors given by inclusions or the identity. This diagram has the following properties. 
    \begin{enumerate}
        \item The arrows pointing in two directions are simplicial Quillen adjunctions, where the upper (left) arrow is the left part of the adjunction.
        \item The left adjunction part (and, respectively, the right adjunction part) together with the inclusions form a commutative diagram.
        \item All horizontals are Quillen equivalences, which create weak equivalences in both directions.
        \item The upper vertical maps are given by including the fibers under the Quillen bifibrations to $\Pos$.
        \item The diagonals pointing downward are given by left Bousfield localization at inner horn inclusions.
        \item The lower verticals that point downward are given by right Bousfield localization at the refinement maps $(-)^{\ared} \to 1$ (see \cref{subsec:refined_strat_spaces}, for the topological case).
    \end{enumerate}
\end{theorem}
% \subsection{Some Notation To Try Out}
% \begin{enumerate}
%     \item $\textnormal{Strat}$
%     \item $\textnormal{Strat}_P$
%     \item $\textbf{Strat}$
%     \item $\textbf{Strat}_P$
%     \item $\underline{\textnormal{Strat}}$
%     \item $\underline{\textbf{Strat}}$
%     \item $\textnormal{\textbf{\textit{S}trat}}$
%     \item $\mathcal{S}\textnormal{trat}$
%         \item $\mathcal{S}\textnormal{\textbf{trat}}$
% \end{enumerate}
\subsection{Topological stratified mapping spaces}
Give two layered $\infty$-categories $X$ and $Y$, the $\infty$-category of functors $Y^X$ is itself layered. Indeed, it follows from the fact that isomorphisms of functors are detected pointwise that it even suffices for $Y$ to be layered.
Together with \cref{thm:overview_over_all_hypothesis}, this already shows that at least the homotopy theory $\AltStratCR$ admits a notion of internal mapping space on the $\infty$-categorical level.
Similarly, in the world of topological stratified spaces (more specifically homotopically stratified spaces) \cite{HughesPathSpaces} equipped the space of stratified maps with a natural decomposition (which generally may not be a stratification) and investigated the lifting properties of such mapping spaces. \\
If we are looking to bring these two notions of internal (i.e. stratified) mapping spaces - one on the $\infty$-categorical and one on the $1$-categorical level together - the language of (semi-)model categories provides the notion of a cartesian closed (semi-)model category (see, for example, \cite{RezkCartesian} for the case of model categories).
\begin{definition}\label{def:cart_closed_semi}
    Let $\cat[M]$ be a semi-model category the underlying $1$-category of which is cartesian closed. We say that $\cat[M]$ is a \define{cartesian closed semi-model category}, if for every cofibration $i \colon A \to B$ and every fibration $f \colon X \to Y$, the induced morphism of exponential objects
        \[
        X^B \to Y^B \times_{Y^A} X^A
        \]
    is a fibration and furthermore is a weak equivalence, if additionally $f$ is a weak equivalence, or $i$ is a weak equivalence with cofibrant source. 
    \end{definition}
    In \cite[Subsec. \ref{comb:subsec:strat_mapping_spaces}]{ComModelWa}, we proved that all of the model structures for stratified homotopy theory defined on $\sStratN$ (see \cref{subsec:comb_mod_hostrat}) are, in fact, cartesian closed. We may use the following lemma to transfer this result to the topological setting.
    \begin{lemma}
        Let $\cat[M]$ be a cartesian semi-model category and $\textnormal{\textbf{N}}$ a cartesian closed, bicomplete category. Suppose that 
        there is an adjunction
        \[
        L \colon \cat[M] \rightleftharpoons \textnormal{\textbf{N}} \colon R,
        \]
        such that $L$ preserves finite products and the right transferred semi-model structure on $\textnormal{\textbf{N}}$ exists. Then, together with this model structure, $\textnormal{\textbf{N}}$ is a cartesian closed semi-model category.
    \end{lemma}
    \begin{proof}
    Note first that if $L$ preserves products, then there are canonical isomorphisms 
    \[
    R(X^{L(A})) \cong R(X)^A
    \]
    for $A \in \cat[M]$ and $X \in \textnormal{\textbf{N}}$. Now, let us verify the hypothesis in the definition of a cartesian closed semi-model category.
    Suppose that we are given a cofibration $i \colon A \to B$ and a fibration $f \colon X \to Y$ in the transferred semi-model structure on $\textnormal{\textbf{N}}$.
    In the context of \cref{def:cart_closed_semi} in $\textnormal{\textbf{N}}$, let us first reduce to the case where $A \to B$ is of the form $L(A' \to B')$, for $A \to B$ a cofibration (acyclic cofibration with cofibrant source). For the general case, we need to show that $R(X^B \to Y^B \times_{Y^A} X^A)$ is a fibration in $\cat[M]$, or equivalently, that every lifting problem 
    \begin{diagram}
        L(A') \arrow[r] \arrow[d] &  X^B \arrow[d]\\ 
        L(B') \arrow[r] \arrow[ru, dashed] & Y^B \times_{Y^A} X^A
    \end{diagram}
    with $A' \to B'$ an (acyclic) cofibration (with cofibrant source) has a solution. Every such lifting problem is, in turn, equivalent to a lifting problem of the form
    \begin{diagram}
        A \arrow[r] \arrow[d] &  X^{L(B')} \arrow[d]\\ 
        B \arrow[r] \arrow[ru, dashed] & Y^{L(B')} \times_{Y^{L(A')}} X^{L(A')}.
    \end{diagram}
    Hence, it suffices to see that the right-hand map is an (acyclic) fibration, which implies the reduction. Finally, under the natural isomorphism we mentioned at the beginning of the proof, using the fact that right adjoint functors preserve fiber products, the map \[ R(X^{L(B')}  \to  R(Y^{L(B')} \times_{Y^{L(A')}} X^{L(A')}) \] is equivalently given by the induced morphism \[ R(X)^{B'} \to R(Y)^{B'} \times_{R(Y)^{A'}} R(X)^{A'},\] which is an (acyclic) fibration in $\cat[M]$ by the assumption that the latter is a cartesian semi-model category.
    \end{proof}
    \begin{corollary}\label{cor:cart_closed}
                Whenever $\TopN$ is a cartesian closed category (i.e., in the compactly generated or $\Delta$-generated case), the semi-model categories $\StratDN, \StratDRN, \StratCN, \StratCRN$ are cartesian closed.
    \end{corollary}
    From this corollary, one may deduce that whenever $\tstr$ is a fibrant stratified space and $\tstr[A]$ a cofibrant stratified space, in any of the semi-model structures on $\StratN$, then the stratified mapping space $\tstr^{\tstr[A]}$ presents the $\infty$-categorical internal mapping space in the corresponding $\infty$-category. In \cref{subsec:frontier_cond_and_refinement}, we discuss the relationship with the stratified mapping spaces of \cite{HughesPathSpaces}.
\subsection{Applying the semi-model structures of stratified spaces}\label{subsec:answers}
\cref{cor:semi_mod_transf} opens up the study of homotopy theories of stratified spaces to all of the tools available for cofibrantly generated simplicial semi-model categories. Notably, these essentially include all of the tools available for model categories, as long as one is willing to assume cofibrancy in the appropriate places (see, for example, \cite{BarwickLeftRight,white2023left} for a good overview of the literature). In this section, we are going to name a few that stand out, as they provide answers to some questions that have been open for a while (see \cref{iQ:Exitcon,iQ:LocatSHE,iQ:ComputationHocolim} in the introduction). Let us begin by noting that the situation is even better than what one can usually expect from a semi-model category. Namely, the simplicial structure gives us simplicial resolutions of fibrant objects (see \cite{FunctionComplexesDwyerKan}), without having to cofibrantly replace first. 
\begin{corollary}
    Let $\underline{\textnormal {\textbf{S}}}$ be any of the simplicial semi-model categories for stratified spaces of \cref{not:all_top_strat_cat}. Then, for any $ \str \in \underline{\textnormal {\textbf{S}}}$ which is fibrant, the functor 
    \begin{align*}
          \Delta^{\op} \to \StratN \\
        [n] \mapsto \tstr^{\Delta^{n}}
    \end{align*}
    defines a simplicial resolution of $\tstr$.
\end{corollary}
In particular, even non-cofibrant objects that are fibrant admit path-objects, and right and left homotopy classes of morphisms with a cofibrant source and fibrant target agree.
We may now compute $\infty$-categorical mapping spaces of stratified spaces in terms of the simplicial structures. 
\begin{corollary}\label{cor:computing_mapping_spaces}
    Let $\underline{\textnormal{\textbf{S}}}$ be any of the simplicial semi-model categories for stratified spaces defined in \cref{subsec:proof_ex_top_modstruct}. Let $W$ be the class of weak equivalences in $\cat[S]$ and denote $\cat[S][W^{-1}]$ the quasi-categorical localization of $\cat[S]$ at $W$. Let $\str[A] \in \cat[S]$ be cofibrant and $\str[X]$ be fibrant. Then there is a natural (zigzag) of weak equivalences between mapping spaces
    \[
    \cat[S][W^{-1}](\str[A], \str[X]) \simeq \underline{\textnormal{\textbf{S}}}(\str[A], \str[X]).
    \]
    In particular, the canonical map bijection
    \[
     [\str[A], \str[X]]_s \to \pi_0 ( \cat[S][W^{-1}](\str[A], \str[X]))
    \]
    is a bijection, where the right-hand expression denotes homotopy classes of stratified (stratum-preserving when the poset is fixed) maps with respect to stratified homotopies.
\end{corollary}
\begin{proof}
    This follows from \cite[Prop. 1.1.10]{goerss2005moduli} which is the left semi-model category version of the argument in \cite[§7]{FunctionComplexesDwyerKan}, together with \cite[1.2]{HinichDwyerLoc}.
\end{proof}
Note that \cref{cor:computing_mapping_spaces} applies to all piecewise linear pseudomanifolds or, more generally, conically stratified spaces, (equipped with the refined stratification, when necessary) and provides, in particular, a strengthened version of \cite[Theorem 5.7]{douteauwaas2021}.\\
From a global perspective, we obtain the following result, which states that the homotopy theory of stratified spaces obtained by inverting weak equivalences may equivalently be described in terms of the simplicial categories of bifibrant stratified spaces, or equivalently a localization of bifibrant stratified spaces at stratified homotopy equivalences, answering \cref{iQ:LocatSHE}.
\begin{corollary}\label{cor:equ_char_of_infty_cat}
   In the setting of \cref{cor:computing_mapping_spaces}, denote by $\underline{\cat[S]}^{o}$ the full simplicial subcategory of bifibrant objects. Denote by $H_s$ the class of stratified homotopy equivalences between bifibrant stratified spaces. 
   There are canonical equivalences of $\infty$-categories 
            \[
            \underline{\cat[S]}^{o} \simeq \cat[S]^o[H_s^{-1}] \simeq \cat[S][W^{-1}]
            \]
    where equivalences between simplicial and quasi-categories are to be understood in terms of the Quillen equivalence between quasi-categories and simplicial categories of \cite{BergnerSimCat}.
\end{corollary}
\begin{proof}
    This is the semi-model category version of \cite[Prop. 4.8]{FunctionComplexesDwyerKan} the proof of which is identical. We note that the first equivalence also follows from the existence of cofibrant replacements and the powering of the stratified categories over $\sSetN$, using a similar flattening argument as in \cite[2.5]{DwyerKanEqu}. Therefore, the full power of a semi-model category is not necessary to obtain this result. The second equivalence, however, needs both fibrant replacement and cofibrant replacement, as well as the Whitehead theorem in a semi-model category.
\end{proof}
Furthermore, we obtain the following two versions of the stratified homotopy hypothesis, providing an answer to \cref{iQ:Exitcon}.
\begin{theorem}\label{thm:strat_ho_hy_quillen}
Mapping a simplicial set to $K$ to the stratified realization of its stratified simplicial set $\lstr{(K)}$, and conversely mapping a stratified space $\str$ to the underlying simplicial set of $\SingS \str$, $\forget (\SingS \str)$, induces a Quillen equivalence of (semi-)model categories 
\[
\sSetOrd \xrightleftharpoons{\simeq} {\StratCRN},
\]
that creates weak equivalences in both directions. \\
The left-hand model category $\sSetOrd$ is the left Bousfield localization of $\sSetJoy$ which presents the full reflective subcategory $\iCatO \subset \iCat$ of small $\infty$-categories in which every endomorphism is an isomorphism.
\end{theorem}
\begin{proof}
    This is the combination of \cref{cor:semi_mod_transf,thm:ex_mod_struct_red} and \cite[Thm. \ref{comb:prop:Quillen_Equ_betw_ref_and_ord}]{ComModelWa}. The only thing we need to pay extra attention to is the statement that both functors create weak equivalences.
    In the case of the right adjoint, note that the forgetful functor \[
    \sStratCRN \to \sSetOrd
    \]
    does not create weak equivalences. However, it creates weak equivalences between such stratified simplicial sets whose strata are Kan-complexes (\cite[Prop. \ref{comb:prop:char_of_ref_cat_equ}]{ComModelWa}). Since $\SingS$ has image in this category, it follows that the composition of $\SingS$ with the functor forgetting the stratification creates weak equivalences. Conversely, note that $\sReal{-} \circ \lstr$ is left Quillen, with source a cofibrant model category. Since the adjunction above is a Quillen equivalence (as the composition of two Quillen equivalences), the derived unit of adjunction is a weak equivalence.
    For $X \in \sSetN$,
    It is given by 
    \[
    X \to \forget ({\SingS\sReal{\lstr(X)}}) \to \forget (\SingS \tstr[Y]) 
    \]
    where $\sReal{\lstr(X)} \to \tstr[Y]$ is a fibrant replacement of $\sReal{\lstr(X)}$. Note, however, that as we have already shown that $\forget \circ \SingS$ preserves weak equivalences, it follows by two-out-of three that the ordinary unit of adjunction is also a weak equivalence.
    It follows, again by two-out-of three, that $f\colon X \to Y$ in $\sSetOrd$ is a weak equivalence, if and only if $\forget({\SingS\sReal{\lstr(f)}})$ is a weak equivalence.
    Finally, since $\forget \circ \SingS$ creates weak equivalences, this is equivalent to $\sReal{-} \circ \lstr$ creating weak equivalences.
\end{proof}
A quasi-categorical counterpart of this result, using \cref{cor:equ_char_of_infty_cat} is:
\begin{corollary}
\label{thm:equ_bifib_and_layered_qc}
    Denote by $\StratN^{o}$ the full subcategory of $\StratCRN$ of bifibrant stratified spaces. Denote by $H_s$ the wide subcategory of stratified homotopy equivalences. 
     % The homotopy coherent nerve of $\Strat^{o}$ is naturally weakly equivalent to
    % the quasi-category given by the localization $\StratN[H_s^{-1}]$.\\
    $\SingS$ (i.e., Lurie's exit-path construction of \cite{HigherAlgebra}) induces an equivalence of quasi-categories
    \[
    \textnormal{Exit} \colon \StratN^o[H_s^{-1}] \xrightarrow{\simeq} \iCatO,
    \]
     where $\iCatO$ denotes the quasi-category of small quasi-categories, in which every endomorphism is an isomorphism. 
\end{corollary}
Finally, let us comment on \cref{iQ:ComputationHocolim}, concerning the relationship between homotopy colimits and $1$-categorical colimits of stratified spaces. 
The main application for the methods of regular neighborhoods in stratified cell complexes we developed in \cite{HoLinksWa} was to prove that certain pushout diagrams of stratified spaces are homotopy pushout. Note that, a posteriori, this also follows from the existence of semi-model structures making inclusions of cellularly stratified spaces cofibrations. In fact, we can now use the full machinery for the computation of homotopy colimits in a simplicial (semi-)model category (see, for example, \cite{hirschhornModel}) to compute the latter.
As a particular consequence, we obtain the following result.
\begin{corollary}
    Suppose we are given a pushout diagram of in $\TopPN$
    \begin{diagram}\label{diag:cor_homotopy_pushout}
        {\tstr[A]} \arrow[r,"c", hook] \arrow[d]& {\tstr[B]} \arrow[d] \\
    {\tstr[X]} \arrow[r, hook]  & {{\tstr[Y]}},
    \end{diagram}
    all of which are cofibrant in $\TopPDN$ or equivalently $\TopPCN$ (that is, retracts of cellularly stratified spaces, by \cref{prop:char_cofibrants_wo_ref}), with the upper vertical given by a cofibration. Then this diagram descends to a pushout diagram in $\AltTopPD$ and $\AltTopPC$. \\
    In particular, this holds when all spaces involved are cellularly stratified spaces and $c$ is an inclusion of cellularly stratified spaces. The analogous claim for $\StratN$ and $\StratDN$, $\StratCN$ holds.
\end{corollary}
Clearly, \cref{lem:pushout_of_cell} is implied by this result. \cref{lem:pushout_of_cell} was used to prove the existence of semi-model structures. Now, however, we can see that the existence of semi-model structures is essentially equivalent to diagrams of stratified cell complexes as in \cref{diag:cor_homotopy_pushout} being homotopy pushout, which illustrates the importance of \cref{iQ:ComputationHocolim}.
% \begin{proof}
%     $\iCatO$ is equivalent to the quasi-category obtained by localizing the relative category $s\textnormal{\textbf{\textit{S}et}}^{\mathfrak{O}}$. By \cref{thm:hohy_nand_cat}, $\SingS$ induces an equivalence of relative categories of $\RelStratCR \xrightarrow{\simeq} s\textnormal{\textbf{\textit{S}et}}^{\mathfrak{O}}$. Consider the inclusion of relative categories $(\CFFN, H) \hookrightarrow \RelStratCR$. This is a homotopy equivalence of relative categories. Indeed, every CFF-stratified space is cofibrant-fibrant by \cref{cor:perf_strat_sp_are_essentially_surj}, and hence the stratified homotopy equivalences between CFF-stratified spaces are precisely the refined categorical equivalences. A homotopy inverse to this inclusion is then given, by first cofibrantly replacing, and then fibrantly replacing through the small object argument, applied applied to the classes $\{\sReal{\stratBound \hookrightarrow \stratSim} \mid n\in \mathbb N\}$ and $\{\sReal{\stratHorn \hookrightarrow \stratSim} \mid 0<k <n\}$. Hence, combinining this equivalence with \cref{thm:strat_ho_hy_quillen}, we obtain an equivalence of relative categories
%     \[
%     (\CFFN, H) \xrightarrow{\simeq } s\textnormal{\textbf{\textit{S}et}}^{\mathfrak{O}}.
%     \]The equivalence of quasi-categories in the statement of the theorem is obtained by passing to localizations of quasi-categories, at the respective class of weak equivalences.
% \end{proof}

%% file: 5BifibrantObjects.tex
\section{(Co)fibrant stratified spaces}\label{sec:cofibrants}
In the previous section, we have established the existence of semi-model structures of stratified spaces, which were transferred from the combinatorial setting of stratified simplicial sets. These were constructed with the specific goal of connecting the various homotopy theories of stratified spaces with their topology and geometry. If we are looking to deepen our understanding of this connection, the obvious question at hand is how classical examples of stratified spaces interact with these model structures. In particular, which classical examples of stratified spaces are bifibrant (\cref{iQ:Whatfibrant}). For example, \cite{miller2013,douteauEnTop,nand2019simplicial} all showed Whitehead theorems for certain stratified spaces, which characterize the stratified homotopy equivalences between them in terms of a priori weaker conditions. Every simplicial semi-model category comes with its own Whitehead theorem, stating that the homotopy equivalences with respect to the simplicial cylinder between the bifibrant objects are precisely the weak equivalences.
The case of $\sStratCRN$ was already proven through a direct proof in \cite{nand2019simplicial}. Summarizing old and proving new criteria for (co)fibrancy and studying how these properties relate to more classical properties of stratified spaces is the main  content of this section. Finally, at the end of this section, we provide a short investigation of stratified versions of the homotopy link as investigated in \cite{hughes1999stratifications} in our model categorical framework.
\subsection{Fibrant stratified spaces and Quinn's homotopically stratified sets}
Let us begin by studying the class of fibrant objects in the model structures constructed in \cref{subsec:proof_ex_top_modstruct}. 
The question of what classical examples are fibrant in $\StratCRN$ was already investigated in \cite{nand2019simplicial} and \cite{HigherAlgebra}, with the core results that both Quinn's homotopically stratified spaces and Siebenmann's conically stratified spaces (see later in this section) are such that the associated stratified singular simplicial set is a quasi-category. Hence, it turns out that most examples of stratified spaces of classical geometrical interest fall into the class of fibrant objects in the semi-model structures we defined. The central new contribution from our side in this section is that we prove that, as long as one restricts to metrizable spaces, the class of fibrant objects is in fact independent of the choice of model structure we presented here (see \cref{prop:equi_of_fibrancy}). 
A consequence of this is \cref{thm:all_theories_same}, which may be taken as the statement that in any geometric scenario there really is no difference between the diagrammatic homotopy theories for stratified spaces and their categorical counterparts at all. \\
\\
Before we begin, note that we really only need to study fibrancy in the scenario of a fixed poset, due to the following fact.
\begin{lemma}\label{lem:equ_fibrancy_cond}
Let $\tstr \in \StratN$. Then the following statements are equivalent:
    \begin{enumerate}
        \item $\tstr$ has the horn filling property over $\ptstr$ with respect to all realizations of admissible stratified horn inclusions. That is, it admits fillers
            \begin{diagram}
                \sReal{\Lambda^\J_k} \arrow[r] \arrow[d]& \tstr \\
                \sReal{\Delta^\J} \arrow[ru, dashed] & 
            \end{diagram}
        in $\TopPN[\ptstr]$, for all admissible pairs $\J \in \catFlag$, $0 \leq k \leq n_{\J}$.
        \item $\tstr$ is fibrant in all of the semi-model categories $\TopPDN[\ptstr], \StratDN,\StratDRN$.
        \item $\tstr$ is fibrant in one of the semi-model categories $\TopPDN[\ptstr], \StratDN,\StratDRN$.
    \end{enumerate}
Furthermore, the following analogous statements about the categorical semi-model are equivalent:
    \begin{enumerate}
        \item $\tstr$ has the horn filling property over $\ptstr$ with respect to all realizations of admissible stratified horn inclusions, and all inner stratified horn inclusions.
        %     \begin{diagram}
        %         \sReal{\Lambda^\J_k}  \arrow[r] \arrow[d] & \tstr \\
        %         \sReal{\Delta^\J} \arrow[ru, dashed] & 
        %     \end{diagram}
        % in $\TopPN[\ptstr]$, for all admissible pairs $\J \in \catFlag$, $0 \leq k \leq n_{\J}$ and whenever $0 < k < n_{\J}$;
        \item $\tstr$ has the horn filling property over $\ptstr$ with respect to all realizations of inner stratified horn inclusions.
        \item $\tstr$ is fibrant in all of the semi-model categories $\TopPCN[\ptstr], \StratCN, \StratCRN$.
        \item $\tstr$ is fibrant in one of the semi-model categories $\TopPCN[\ptstr], \StratCN, \StratCRN$.
    \end{enumerate}

\end{lemma}
\begin{proof}
The diagrammatic case is immediate by the construction of the left model structures via left-transfer and \cref{prop:model_global_mod_struct_are_glued} together with the fact that $\StratDRN$ is obtained from $\StratDN$ via a right Bousfield localization \cref{thm:huge_overview}. For the categorical statement, the proof is almost identical; special attention only needs to be paid to the implication from the second to the first property. Note that for any stratified space $\tstr$ the strata of the stratified-simplicial set $\SingS{\tstr}$ are Kan-complexes. Hence, the implication follows from \cite[Prop. 2.2.3]{haine2018homotopy}.
\end{proof}
\begin{definition}
    We call a stratified topological space $\tstr$ that satisfies any of the equivalent conditions of \cref{lem:equ_fibrancy_cond} concerning the diagrammatic (categorical) semi-model structures \define{diagrammatically} (\define{categorically}) fibrant.
\end{definition}
Clearly, every categorically fibrant stratified space is diagrammatically fibrant. The converse must necessarily be false. If it were true, then by \cref{cor:equ_char_of_infty_cat,thm:overview_over_all_hypothesis} together with \cite[Thm. 1.1.7]{haine2018homotopy} and \cite[Thm. 3]{douteauEnTop}, it would imply that the inclusion from the $\infty$-category of d\'ecollages into $\FunC (\sd(\pos)^{\op}, \iSpaces)$ is an equivalence of $
\infty$-categories\footnote{It is easy to write down an example of a diagram which is not a d\'ecollage. For example, take the constant diagram over $(\sd [2])^{\op}$ and replace the entry at $[0 <1<2]$ by the empty set.}.
However, it is surprisingly hard to provide a geometrical example for a stratified space that is diagrammatically fibrant and not categorically fibrant (depending on your definition of geometric, it is not possible at all). Before we investigate this question further, let us first recall results on which classical examples of stratified space we can expect to be fibrant.
\begin{remark}
    It was first shown in \cite[A.5]{HigherAlgebra} that all conically stratified spaces - roughly a stratified space that locally has the structure $U \times C(\tstr[L])$, where $C(\tstr[L])$ is the stratified cone of a stratified space $\tstr[L]$ - are categorically fibrant.
    In particular, this implies that all topological pseudomanifolds (see, for example, \cite{banagl2007topological}), as well as all conically smooth stratified spaces (see \cite{LocalStructOnStrat}) and hence also all Whitney stratified spaces (\cite{nocera_whitney_2023} or classically \cite{thom1969ensembles}), are categorically and hence also diagrammatically fibrant. 
\end{remark}
It turns out, however, that fibrancy can already be obtained under significantly less geometric conditions. One of the crucial insights contributing to the foundations of stratified homotopy theory was Quinn's observation that a powerful homotopy theoretical setup for stratified spaces can already be achieved by only posing requirements on pairs of strata (\cite{quinn1988homotopically}). Quinn defined a notion of homotopically stratified set (with slightly more restrictive conditions on the stratification poset) by requiring $\tstr \in \StratN$ to be metrizable, and requiring that for any two-element flag $[p < q] \in \catFlag$ it holds that:
\begin{enumerate}
    \item \label{prop:QuinStrat1} $\tstr_{p} \hookrightarrow \tstr_{p<q}$ is tame, i.e. admits a nearly strict neighborhood deformation retraction (see \cite{quinn1988homotopically}).
    \item \label{prop:QuinStrat2}$\HolIP[p<q]{\tstr} \to \tstr_{p}$ is a Hurewicz fibration.
\end{enumerate}
Together, these assumptions provided an excellent homotopy theoretic framework for classification results concerning manifold stratified spaces (see \cite{quinn1988homotopically,weinberger1994topological}).
Note that both of these assumptions are more in line with the Hurewicz approach (see \cite{strom_homotopy_1972}), than with the Serre-Quillen approach to homotopy theory (\cite{Quillen}). Since the approach we pursue here is combinatorial (i.e. following Serre and Quillen), we should expect cofibrancy in our semi-model categories to be a stronger condition than condition (\ref{prop:QuinStrat1}) above. Indeed, it is a consequence of \cite[Prop. \ref{hol:prop:aspire_for_finite_cell}]{HoLinksWa} that any locally compact cofibrant stratified space (in any of the model structures of \cref{subsec:proof_ex_top_modstruct}) satisfies the tameness condition above.
Conversely, in the two strata case it is not hard to see that $\HolIP[p<q]{\tstr} \to \tstr_{p}$ being a Serre fibration is equivalent to diagrammatic and equivalently categorical fibrancy (see the proof of \cref{prop:equi_of_fibrancy} below). 
Surprisingly, it follows from \cite[Thm. 4.9]{Miller2009PopathsAH} that every homotopically stratified set is, in fact, categorically fibrant (\cite[Prop. 8.1.2.6]{nand2019simplicial}), without any restrictions on the stratification poset. 
Note that being fibrant is manifestly a statement about the interaction of more than two strata, which makes it so surprising that it can be achieved by only making requirements on two-strata interaction.
Even more, this result does not leave much space for diagrammatically fibrant spaces that are not categorically fibrant at all (see \cref{iQ:DifferenceBetweenTheories}). If we assume pairwise tameness and metrizability, then it would follow that any counterexample would involve $\HolIP[p<q]{\tstr} \to \tstr_{p}$ being a Serre fibration, but not a Hurewicz fibration. It turns out that no cofibrancy assumptions are necessary whatsoever, and in fact in any geometric setting diagrammatic and categorical fibrancy are equivalent conditions, which are furthermore equivalent to the Serre fibration analogue of Quinn's pairwise homotopy link condition. Furthermore, this result does not require any additional tameness assumptions.
\begin{proposition}\label{prop:equi_of_fibrancy}
    Let $\tstr \in \StratN$ be a metrizable stratified space. Then the following conditions are equivalent:
    \begin{fibrancyProp}
        \item \label{prop:equ_fibrancy_1} $\tstr$ is categorically fibrant.
        \item \label{prop:equ_fibrancy_2} $\tstr$ is diagrammatically fibrant.
         \item \label{prop:equ_fibrancy_3} For any pair $[p<q] \in \catFlag[\ptstr]$, the starting point evaluation map $\HolIP[p<q](\tstr) \to \tstr_{p}$ is a Serre fibration.
        \item \label{prop:equ_fibrancy_4} For any pair $[p<q] \in \catFlag[\ptstr]$, and for any flag $\J = [p = \cdots = p_k < q = \cdots = q = p_{n_\J}]$ with $k \geq 1$, $\tstr$ has the horn filling property with respect to $\sReal{\Lambda^\J_k \hookrightarrow \Delta^\J}$.
    \end{fibrancyProp}
\end{proposition}
    Before we give a proof, we note the implications of this result. By \cref{thm:huge_overview}, the categorical setting is always obtained as a left Bousfield localization of the diagrammatic one. In particular, the two settings have the same cofibrations; furthermore, we obtain:
    \begin{corollary}\label{cor:equ_of_fibrations}
        Let $\tstr,\tstr[Y] \in \StratN$ be metrizable stratified spaces, $\tstr$ diagrammatically fibrant and $f \colon \tstr[Y] \to \tstr$ a stratified map. Then there are equivalences:
        \begin{enumerate}
            \item $f$ is a fibration in $\TopPDN$ $\iff$ $f$ is a fibration in $\TopPCN$;
            \item $f$ is a fibration in $\StratDN$ $\iff$ $f$ is a fibration in $\StratCN$;
            \item $f$ is a fibration in $\StratDRN$ $\iff$ $f$ is a fibration in $\StratCRN$.
        \end{enumerate}
    \end{corollary}
    \begin{proof}
        This is immediate from \cref{prop:equi_of_fibrancy}, together with \cref{thm:huge_overview}, and the characterization of fibrations between fibrant objects in a left Bousfield localizations (the semi-model category version of \cite[Prop. 3.3.16.]{hirschhornModel}, see \cite[Rem. 4.5]{white2018bousfield}.)
    \end{proof}
        We may hence interpret \cref{prop:equi_of_fibrancy} as stating that in any geometrical framework the diagrammatic approach is really the same as the categorical approach, obtaining an answer to \cref{iQ:DifferenceBetweenTheories}.
        More rigorously, we can phrase this insight as follows.
        \begin{theorem}\label{thm:all_theories_same}
            Let $\underline{\textnormal{\textbf{D}}}$ be any of the diagrammatic simplicial semi-model categories of \cref{subsec:proof_ex_top_modstruct} and let $\underline{\textnormal{\textbf{C}}}$ be its categorical pendant. Denote by $\underline{\textnormal{\textbf{D}}}^{f}_{m}$ and $\underline{\textnormal{\textbf{C}}}^{f}_{m}$ the respective restrictions of the semi-model structures to the full simplicial subcategory of fibrant, metrizable stratified spaces. Then equality \[ \underline{\textnormal{\textbf{D}}}^{f}_{m} = \underline{\textnormal{\textbf{C}}}^{f}_{m} \]
            holds (on the nose). In particular, if we denote by $\underline{\textnormal{\textbf{D}}}^{o}_{m}$ and $\underline{\textnormal{\textbf{C}}}^{o}_{m}$ the corresponding simplicial categories of bifibrant metrizable stratified spaces, then 
            \[
            \underline{\textnormal{\textbf{D}}}^{o}_{m} = \underline{\textnormal{\textbf{C}}}^{o}_{m}
            \]
            and hence the full sub-$\infty$-categories of bifibrant metrizable stratified spaces agree.
        \end{theorem}
        \begin{proof}
            Equality on the object level follows by \cref{prop:equi_of_fibrancy}. Equality of cofibrations holds even without additional assumptions. Equality of fibrations was shown in \cref{cor:equ_of_fibrations}. It remains to verify the equality of weak equivalences. We prove the stratum-preserving case; the others are analogous: Let $f \colon \tstr \to \tstr[Y]$ be a stratum-preserving map in $\underline{\textnormal{\textbf{D}}}^{f}_{m}$. We may factorize $f$ as a cofibration and an acyclic fibration
            \[
            \tstr \to \tstr[Z] \to \tstr[Y],
            \] 
            with $\tstr[Z]$ not necessarily metrizable.
            Note that since cofibrations in both structures agree, so do acyclic fibrations. 
            By \cref{prop:equi_of_fibrancy}, it follows that $\tstr[Z]$ is also categorically fibrant (this did not use any metrizability assumptions on $\tstr[Z]$).
            $f$ is a weak equivalence if and only if the first of these two maps is a weak equivalence. However, the latter is a map between categorically fibrant spaces. In particular, by the semi-model category version of the Whitehead theorem for Bousfield localizations (see \cite[Thm 3.2.13]{hirschhornModel} for the model category version)\footnote{\cite[Rem. 4.5]{white2023left} states the semi-model category case for bifibrant objects in the localized model structure. Note, however, that we may always assume bifibrancy, by replacing both source and target cofibrantly first, which is done through acyclic fibrations, which are the same in both model structures. Alternatively, we can just use the full model categorical result and pass to the simplicial setting using the fact that all semi-model structures are transferred.}, it is a diagrammatic equivalence if and only if it is a categorical equivalence.
        \end{proof}
        In other words, for geometric examples, the two homotopy theories agree. \\
        Before we provide a proof of \cref{prop:equi_of_fibrancy}, we need the following elementary set-theoretic topological lemma, the proof of which is provided in \cref{appendix:elementary_lem}. It is the only place where the metrizability of $\tstr$ comes into play. In the following, $D^{n+1}$ denotes the euclidean unit disk of dimension $n+1$, for $n \in \mathbb{N}$, $S^{n}$ denotes its boundary, and $\mathring{D}^{n+1} =D^{n+1} \setminus S^n$ its interior.
         \begin{lemma}\label{lem:extension_lemma}
        Let $n \geq 0$, $X$ a metrizable topological space and suppose that we are given a solid commutative diagram of (general) topological spaces of the following form:
        \begin{diagram}
                        {}                                                                 & D^{n+1} \times \{ 0 \} \cup S^{n} \times [0,1] \arrow[rd, "f"] \arrow[d, hook] &  \\
           D^{n+1} \times \{ 0 \} \cup S^{n} \times [0,1) \arrow[ru, hook] \arrow[rd, hook]  & D^{n+1} \times [0,1] \arrow[r, dashed, "\tilde{f}"] & X \\
                                                                                    & D^{n+1} \times [0,1) \arrow[u, hook] \arrow[ru, "f'"'] & \spaceperiod             
        \end{diagram}
         Then there exists $\tilde f \colon   D^{n+1} \times [0,1] \to X$ making the upper triangle commute, and furthermore such that the inclusion \[\tilde{f}(\mathring D^{n+1} \times [0,1] ) \subset f'(\mathring D^{n+1} \times [0,1) )\]
    holds.
    \end{lemma}
        As a consequence of \cref{lem:extension_lemma} we obtain:
\begin{lemma}\label{lem:metric_reduction_lemma}
    Let $\tstr \in \TopPN$ be a metrizable stratified space and $\J = [p_0 \leq \cdots \leq p_n]$ be a flag in $P$. Let $0 < k < n$ and denote by $\I$ the regular flag containing all $p_i \in \J$ with $i \geq k$. If $\tstr_{\I}$ is categorically fibrant, then $\tstr$ has the horn filling property with respect to $\sReal{\Lambda^\J_k \hookrightarrow \Delta^\J}$.
\end{lemma}
\begin{proof}
Let $\J_k$ be the flag obtained by removing $p_k$ from $\J$, corresponding to the $k$-th face of $\Delta^\J$.
We may identify $\sReal{\Delta^\J}$ with an (appropriately stratified) join $\sReal{\Delta^{[p_k]}} \star \sReal{\Delta^{\J_k}}$.
This identification induces join coordinates $[s,x]$, $s \in [0,1]$, $x \in \sReal{\Delta^{\J_k}}$ on $\sReal{\Delta^\J}$, with the point $[0,x]$ corresponding to the unique element of $\sReal{\Delta^{[p_k]}}$\footnote{Beware that at other points in the text, we have parametrized joins the wother way around.}. 
Under this identification, $\sReal{\Lambda^\J_k}$ corresponds to the join $\sReal{\Delta^{p_k}} \star \sReal{\partial \Delta^{\J_k}}$. 
For $0 \leq \alpha \leq \alpha ' <1$, we write
\begin{align*}
   S^{\leq \alpha'}_{\geq \alpha} = \{ [x,s] \in \sReal{\Delta^{[p_k]}} \star \sReal{\Delta^{\J_k}} \mid \alpha \leq s \leq \alpha' \},
\end{align*}
and use analogous notation replacing $\leq $ by $<$ and $\geq$ by $>$.
We may then decompose $\sReal{\Delta^\J}$ into 
\begin{align*}
  \sReal{\Delta^\J} = S^{\leq \frac{1}{2}}_{\geq 0} \cup S^{\leq 1}_{\geq \frac{1}{2}}.
\end{align*}
It is immediate by the definition of the stratified simplex $\sReal{\Delta^\J}$, that the inclusion 
\[
\real{\Lambda^\J_k} \cap S^{\leq \frac{1}{2}}_{\geq 0} \hookrightarrow S^{\leq \frac{1}{2}}_{\geq 0}
\]
is stratum-preserving homeomorphic to the inclusion
\[
\sReal{\Lambda^{\J'}_k} \hookrightarrow \sReal{\Delta^{\J'}},
\]
where $\Delta^{\J'}$ is obtained by replacing every entry $p_i$ in $\J$ with $i < k$ by $p_k$. In particular, by the assumption on $\tstr_{\I}$, $\tstr$ admits a filler with respect to the latter inclusion. Thus, it suffices to show that $\tstr$ admits fillers with respect to 
\[
A := S^{\leq 1}_{\geq \frac{1}{2}} \cap \real{\Lambda^\J_k} \cup S^{\leq \frac{1}{2}}_{\geq \frac{1}{2}} \hookrightarrow S^{\leq 1}_{\geq \frac{1}{2}} =: B.
\]
Choosing any homeomorphism $\real{\Delta^n} \cong D^{m+1}$, where $m = n-1$, and using the affine order preserving homeomorphism $[\frac{1}{2},1] \cong [0,1]$ we may identify the latter inclusion with a non-stratified inclusion
\[
D^{m+1} \times \{0 \} \cup S^{m} \times [0,1] \hookrightarrow D^{m+1} \times [0 , 1 ].
\]

Since $p_k \leq p_n$, under this identification $\mathring D^{m+1} \times [0,1]$ is entirely contained in the $p_n$-stratum of $\sReal{\Delta^\J}$.
Thus, given a stratum-preserving map $f \colon A \to \tstr$, an extension of $f$ to $B$ is the same data as an extension of $f$ to $D^{m+1} \times [0,1]$ mapping $\mathring D^{m+1} \times [0,1]$ to $\tstr_{p_n}$. We may now apply \cref{lem:extension_lemma}, from which it follows that it suffices to obtain an extension of $f_{D^{m+1} \times \{0 \} \cup S^{m} \times [0,1)}$ to $D^{m+1} \times [0,1) \cong S^{<1}_{\frac{1}{2}}$.
Now write
\[
S^{<1}_{\frac{1}{2}} = \bigcup_{l \geq 1}  S^{\leq 1- \frac{1}{2^{l+1}}}_{\geq 1-\frac{1}{2^{l}}},
\]
which defines a locally finite, closed cover of $S^{<1}_{\frac{1}{2}}$. Inductively, we may hence reduce to solving extension problems with respect to
\[
S^{\leq 1- \frac{1}{2^{l}}}_{\geq 1-\frac{1}{2}^{l}} \cup \big ( \sReal{\Lambda^\J_k} \cap  S^{\leq 1- \frac{1}{2^{l+1}}}_{\geq 1-\frac{1}{2^{l}}}) \hookrightarrow S^{\leq 1- \frac{1}{2^{l+1}}}_{\geq 1-\frac{1}{2^{l}}}.
\]
Finally, if we denote by $\J'_k$ the $k$-th face of $\J'$, then the latter inclusion is stratum-preserving homeomorphic to the simplicial box product
\[
\sReal{\Delta^{\J'_k}} \otimes \Delta^{0} \cup_{\sReal{\partial \Delta^{\J'_k}} \otimes \Delta^{0}} \sReal{\partial \Delta^{\J'_k}} \otimes \Delta^1 \hookrightarrow \sReal{\Delta^{\J'_k}} \otimes \Delta^{1},
\]
for the inclusion of $\Delta^0$ into $\Delta^1$ at $0$. Since $\TopPC[\I]$ is a simplicial semi-model category, the latter defines an acyclic cofibration in $\TopPC[\I]$. Since $\tstr_{\I}$ was assumed to be fibrant, the existence of extensions with respect to such an inclusion follows.
\end{proof}
We may now complete the following proof.
\begin{proof}[Proof of \cref{prop:equi_of_fibrancy}]
    Clearly, \cref{prop:equ_fibrancy_1} implies \cref{prop:equ_fibrancy_2}. To see that \cref{prop:equ_fibrancy_2} implies \cref{prop:equ_fibrancy_3}, note that any lifting problem
    \begin{diagram}
     {}   \real{ \Lambda^n} \arrow[d, hook] \arrow[r] & \HolIP[p<q](\tstr) \arrow[d] \\
        {} \real{\Delta^n} \arrow[r] \arrow[ru, dashed] & \tstr_p 
    \end{diagram}
    is equivalent to a lifting problem (over $\ptstr$)
    \begin{diagram}
       {} \sReal{ \Delta^{[p]}} \otimes \Delta^n \cup_{\sReal{\Delta^{[p]}} \otimes  \Lambda^n} \sReal{\Delta^{[p<q]}} \otimes \Lambda^n  \arrow[d, hook] \arrow[r] & \tstr \\
        {} \sReal{\Delta^{[p<q]}} \otimes \Delta^n \arrow[ru, dashed] & \spaceperiod
    \end{diagram}
    Hence, \cref{prop:equ_fibrancy_3} follows from the simpliciality of the semi-model category $\TopPD[\ptstr]$. For the implication \cref{prop:equ_fibrancy_3} $\implies$ \cref{prop:equ_fibrancy_4}, note that using the identification $\sReal{\Delta^\J}$ with the (appropriately stratified) join $\sReal{\Delta^{\J_p}} \star \sReal{\Delta^{\J_q}}$, we can interpret $\sReal{\Delta^\J}$ as a stratified quotient space of 
    \[
    \sReal{\Delta^{[p<q]}} \otimes ( \Delta^{k} \times \Delta^{n-k-1} ),
    \]
    obtained by collapsing $\sReal{\Delta^{[p]}} \otimes (\Delta^{k} \times \Delta^{n-k-1})$ to $\sReal{\Delta^{[p]}} \otimes \Delta^{k}$ and $\sReal{\Delta^{[q]}} \otimes (\Delta^{k} \times \Delta^{n-k-1})$ to $\sReal{\Delta^{[q]}} \otimes \Delta^{n-k-1}$. Within $ \sReal{\Delta^{[p<q]}} \otimes ( \Delta^{k} \times \Delta^{n-k-1} )$, $\sReal{\Lambda^\J_k}$ corresponds to 
    \begin{diagram}\label{prop:equ_fib_join_inclusion}
         \sReal{\Delta^{[p<q]}} \otimes ( \Lambda_k^{k} \times \Delta^{n-k-1} ) \cup \sReal{\Delta^{[p]}} \otimes (\Delta^{k} \times \Delta^{n-k-1}) \cup  \sReal{\Delta^{q}} \otimes (\Delta^{k} \times \Delta^{n-k-1}) \arrow[d, hook] & \\
         \sReal{\Delta^{[p<q]}} \otimes ( \Delta^{k} \times \Delta^{n-k-1} )
    \end{diagram}
    % \begin{align}
    %     \hookrightarrow 
    % \end{align}
    Hence, it suffices to show that $\tstr$ has the extension property with respect to this stratum-preserving map over $\ptstr$. Denote $A = \Lambda_k^{k} \times \Delta^{n-k-1}$ and $B = \Delta^{k} \times \Delta^{n-k-1}$. 
    Under the identification $\sReal{\Delta^{[p<q]}} \cong \sReal{\Delta^{[p<q]}} \cup_{
    \sReal{\Delta^{[q]}}} \sReal{\Delta^{[q \leq q]}}
    $ we may decompose (\ref*{prop:equ_fib_join_inclusion}) into pushouts of inclusions
    \begin{align}\label{prop:equ_fib_dec_inclusion_1}
        \sReal{\Delta^{[p]}} \otimes B \cup_{\sReal{\Delta^{[p]}} \otimes A}\sReal{\Delta^{[p < q]}} \otimes A \hookrightarrow \sReal{\Delta^{[p < q]}} \otimes B
    \end{align}
    and
    \begin{align}\label{prop:equ_fib_dec_inclusion_2}
        \sReal{\partial \Delta^{[q \leq q]}} \otimes B \cup_{\sReal{\partial \Delta^{[q \leq q]}} \otimes A}\sReal{\Delta^{[q<q]}} \otimes A \hookrightarrow \sReal{\Delta^{[q\leq q]}} \otimes B.
    \end{align}
    The inclusion $A \hookrightarrow B$ is an anodyne map of simplicial sets. It follows that the inclusion (\ref*{prop:equ_fib_dec_inclusion_2}) which lies entirely over one stratum is given by an acyclic Quillen cofibration. Consequently, $\tstr$ has the right lifting property with respect to the inclusion (\ref*{prop:equ_fib_dec_inclusion_2}). Finally, lifting problems with respect to (\ref{prop:equ_fib_dec_inclusion_1}) are equivalent to lifting problems of the form 
    \begin{diagram}
        {}   \real{A} \arrow[d, hook] \arrow[r] & \HolIP[p<q](\tstr) \arrow[d] \\
        {} \real{B} \arrow[r] \arrow[ru, dashed] & \tstr_p \spacecomma
    \end{diagram}
    since $A \hookrightarrow B$ is an anodyne extension, the left-hand vertical is an acyclic cofibration in the Quillen model structure, showing the existence of a lift. \\ It remains to show that \cref{prop:equ_fibrancy_4} implies \cref{prop:equ_fibrancy_1}. It suffices to show that $\tstr_{\I}$ is categorically fibrant for all regular flags $\I =  [q_0 < \cdots < q_n]$ of $\ptstr$.
    We only need to cover the case where $n \geq 1$, that is, where $\J$ has at least two different elements (otherwise, we are in the trivially stratified case, that is, the classical case).
    First, note that \cref{prop:equ_fibrancy_1} already implies the case where $\I = [p <q]$. Indeed, let $\J = [p_0 \leq \cdots \leq p_{n_{\J}}]$ be a flag that degenerates from $\I$.
    When $p_k = p$ we may without loss of generality assume that $k$ is maximal with the property that $p_k = p$, by permuting the corners of $\sReal{\Delta^\J}$ by a stratum-preserving homeomorphism. Then we are in the situation of \cref{prop:equ_fibrancy_4}.
    In case when $p_k = q$, the inclusion $\sReal{\Lambda^\J_k \hookrightarrow \Delta^\J}$ admits a stratum-preserving retraction. (Consider a homeomorphism of $\sReal{\Delta^\J}$ with $D^n$, mapping the $k$-th face to the northern hemisphere. Then the $p$-stratum is entirely contained in the equator, and projecting vertically down to the southern hemisphere defines a retraction.) 
    We now proceed to show the horn filling property with respect to arbitrary $\J$, by induction over $n_{\I}$, the length of $\I$.
    We have already covered the case $n_{\I}=1$. Now for the inductive step let $\J = [p_0 \leq \cdots \leq p_m]$ be a flag degenerating from a subflag of $\I$, and $0 < k < m$.
    For such $\J$, denote by $t_{\J}$ the length of $\J_{p_0}$ and by $o_{\J}$ the length of $\J_{q_{n-1}}$.
    We proceed by double induction on $t_{\J}$ and $o_{\J}$, keeping $k$ flexible. If $o_{\J} = -1$, then $\J$ degenerates from a proper subflag of $\I$, and we are reduced to the inductive assumption in the induction on $n_{\I}$. If $p_0 \neq p_k$, which is in particular the case when $t_{\J} = 0$ (since $0 < k < n$), then $\{q_i \mid i \in [n_{\I }], q_i \geq p_k \}$ has cardinality smaller than $\I$, and by inductive assumption (in $n_{\I}$) we can apply \cref{lem:metric_reduction_lemma} proving the existence of a filler with respect to $\J,k$. So, suppose $t_{\J} \geq 1$, $p_0 = p_k$ and we have already proven the result for all cases $\tilde \J, \tilde k$ where either $t_{\tilde \J} < t_{\J}$ or $o_{\tilde \J} < o_{\J}$. Let $s \in  [n_{\J}]$ be such that $e_{s} \in \sReal{\Delta^\J}$ is the maximal vertex lying in the $q_{n_{\I}-1}$ stratum and $e_{n_{\J}} \in  \sReal{\Delta^\J}$ the $n_{\J}$-th vertex of $\sReal{\Delta^\J}$. We may assume that $s \geq 2$, as $t_{\J} \geq 1$ and $n_{\I} \geq 2$.
    Furthermore, we may assume that $e_{n_{\J}}$ lies in the $q_n$ stratum, as otherwise $\J$ degenerates from a proper subflag of $\I$, and we are reduced to the inductive assumption. Now, let 
    \[
    e' = \frac{1}{2}e_s + \frac{1}{2}e_{n_{\J}}
    \]
    be the halfway point between these two vertices, which by assumption also lies in the $q_{n_{\I}}$ stratum. 
    Consider the two affine stratum-preserving embeddings
    \begin{align*}
        i_1 \colon \sReal{\Delta^\J} \hookrightarrow \sReal{ \Delta^\J}  \\
        e_i \mapsto \begin{cases}
            e_i &  i  < n_{\J} \\
            e' & i = n_{\J}
        \end{cases}
    \end{align*}
    and 
    \begin{align*}
       i_2 \colon \sReal{\Delta^{\J'}} \hookrightarrow \sReal{ \Delta^\J}  \\
        e_i \mapsto \begin{cases}
            e_i &  i \neq s \\
            e' & i = s
        \end{cases}
    \end{align*}
    where $\J'$ is obtained from $\J$ by replacing the $s$-th entry by $q_{n_{\I}}$. The images of these two embeddings cover $\sReal{\Delta^\J}$ and intersect in the convex span 
    \[
    <\{ e', e_i \mid  i\in [n_{\J}], i \neq n_{\J},s \}>.
    \]
     Under $i_1$ this span corresponds to the $s$-th face of $\Delta^\J$, and under $i_2$ to the $n_{\J}$-th face of $\Delta^{\J'}$, which are both given by the flag $\J''$ obtained by removing $p_{s}$ from $\J$.
     We obtain an induced stratum-preserving homeomorphism (see \cref{fig:decomp_quinn_proof} for an illustration)
     \[
     \sReal{\Delta^\J \cup_{\Delta^{\J''}} \Delta^{\J'}} \cong \sReal{\Delta^\J} \cup_{\sReal{\Delta^{\J''}}} \sReal{\Delta^{\J'}} \xrightarrow{(i_1,i_2)} \sReal{\Delta^\J}.
     \]
     \begin{figure}[H]
         \centering
    \begin{tikzpicture}[scale=1]
    \coordinate (e3) at (2,0,0);
    \coordinate (e1) at (0,2,0);
    \coordinate (e0) at (0,0,2);
    \coordinate (e2) at (1.2,1.4,1.4);
    \coordinate (e') at ($0.5*(e2)+0.5*(e3)$);
        \fill[color = blue, color = blue, opacity = 0.5] (e1) -- (e2) -- (e3) -- cycle;

    \fill[fill = violet, opacity = 0.6] (e0) -- (e1) -- (e2) -- cycle;
        \fill[color = blue, color = blue, opacity = 0.6] (e0) -- (e2) -- (e3) -- cycle;
    \draw[thick, color = blue] (e1) -- (e3);
     \draw[thick, color = blue] (e0) -- (e3);
    \draw[thick, color = blue] (e2) -- (e3);
    \draw[thick, color = blue] (e0) -- (e2) -- (e3) -- cycle;
    \draw[thick, color = blue] (e0) -- (e');
    \draw[thick, color = blue] (e1) -- (e');
    \draw[thick, color = violet] (e1) -- (e2);
    \draw[thick, color = violet] (e0) -- (e2);
    \draw[thick, color = red] (e0) -- (e1);
    \node[circle, fill=red, inner sep=0.8pt] at (e0){};
        \node[circle, fill=red, inner sep=0.8pt] at (e1){};
       \node[circle, fill=blue, inner sep=0.8pt] at (e'){};
       \node[circle, fill=blue, inner sep=0.8pt] at (e3){};
           \node[circle, fill=violet, inner sep=0.8pt] at (e2){};

    \node[left] at (e2) {$\color{white} e_2$};
    \node[left] at (e1) {$e_1$};
    \node[below] at (e0) {$e_0$};
    \node[below] at (e3) {$e_3$};
     \node[below] at (e') {$\color{white} {e'}$};
\end{tikzpicture}         \caption{Depiction of $\Delta^\J \cup_{\Delta^{\J''}} \Delta^{\J'}$ with $\Delta^\J$ spanned by $e_0,e_1, e_2, e'$ and $\Delta^{\J'}$ spanned by $e_0,e_1,e',e_3$.}
         \label{fig:decomp_quinn_proof}
     \end{figure}
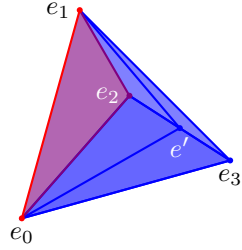
     Under this identification $\Lambda^\J_k$ corresponds to the subcomplex $A$ of $\Delta^\J \cup_{\Delta^{\J''}} \Delta^{\J'}=:B$ given by removing:
     \begin{itemize}
         \item The simplices corresponding to $\Delta^\J$ and $\Delta^{\J'}$;
         \item The $k$-th face of $\Delta^\J$ , $\Delta^{\J_k}$, and of $\Delta^{\J'}$, $\Delta^{\J_k'}$;
         \item The $s$-th face of $\Delta^\J$, $\Delta^{\J_s}$, which is equivalently the $n_{\J}$-th face of $\Delta^{\J'}$, or the top dimensional simplex of $\Delta^{\J''}$.
         \item The $k$-th face of $\Delta^{\J''}$, $\Delta^{\J''_k}$, which is equivalently the $(s-1)$-th face of $\Delta^{\J_k}$.
     \end{itemize}
     \begin{center}
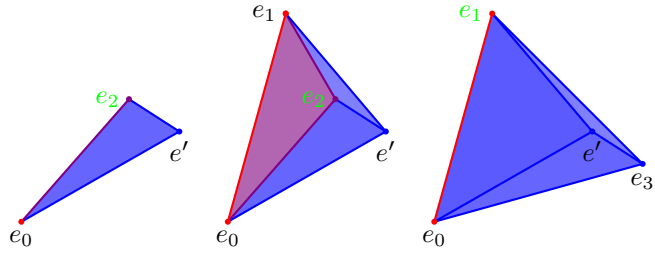
\begin{figure}[H]
    \centering
    \begin{tikzpicture}[scale=1]
    \coordinate (e3) at (2,0,0);
    \coordinate (e1) at (0,2,0);
    \coordinate (e0) at (0,0,2);
    \coordinate (e2) at (1.2,1.4,1.4);
    \coordinate (e') at ($0.5*(e2)+0.5*(e3)$);
        \fill[color = blue, color = blue, opacity = 0.6] (e0) -- (e2) -- (e') -- cycle;
    \draw[thick, color = blue] (e2) -- (e');
    \draw[thick, color = blue] (e0) -- (e');
    \draw[thick, color = violet] (e0) -- (e2);
    \node[circle, fill=red, inner sep=0.8pt] at (e0){};
       \node[circle, fill=blue, inner sep=0.8pt] at (e'){};
           \node[circle, fill=violet, inner sep=0.8pt] at (e2){};
    \node[left] at (e2) {$\color{green}{e_2}$};
    \node[below] at (e0) {$e_0$};
     \node[below] at (e') {$ {e'}$};
    %  \coordinate (sh) at (0,2,0);

    %     \coordinate (e0) at (0,0,2);
    % \coordinate (e2) at ($(e2) + (sh)$);
    % \coordinate (e0) at ($(e0) + (sh)$);
    % \coordinate (e') at ($(e') + (sh)$);
    % \draw[thick, color = blue] (e2) -- (e');
    % \draw[thick, color = violet] (e0) -- (e2);
    % \node[circle, fill=red, inner sep=0.8pt] at (e0){};
    %    \node[circle, fill=blue, inner sep=0.8pt] at (e'){};
    %        \node[circle, fill=violet, inner sep=0.8pt] at (e2){};
    % \node[left] at (e2) {$e_2$};
    % \node[below] at (e0) {$e_0$};
    %  \node[below] at (e') {${e'}$};

\end{tikzpicture}
 \begin{tikzpicture}[scale=1]
    \coordinate (e3) at (2,0,0);
    \coordinate (e1) at (0,2,0);
    \coordinate (e0) at (0,0,2);
    \coordinate (e2) at (1.2,1.4,1.4);
    \coordinate (e') at ($0.5*(e2)+0.5*(e3)$);
    % \coordinate (e'') at (e2);
    % \coordinate(e2) at (e');
    % \coordinate(e') at (e'');
    \fill[color = blue, color = blue, opacity = 0.6] (e0) -- (e2) -- (e') -- cycle;
        \fill[fill = violet, opacity = 0.6] (e0) -- (e1) -- (e2) -- cycle;
    \fill[color = blue, color = blue, opacity = 0.5] (e1) -- (e2) -- (e') -- cycle;

    \draw[thick, color = blue] (e2) -- (e');
    \draw[thick, color = blue] (e0) -- (e');
    \draw[thick, color = violet] (e0) -- (e2);
    \draw[thick, color = red] (e0) -- (e1);
        \draw[thick, color = violet] (e1) -- (e2);
    \draw[thick, color = blue] (e1) -- (e');

    \node[circle, fill=red, inner sep=0.8pt] at (e0){};
       \node[circle, fill=blue, inner sep=0.8pt] at (e'){};
           \node[circle, fill=violet, inner sep=0.8pt] at (e2){};
                   \node[circle, fill=red, inner sep=0.8pt] at (e1){};

    \node[left] at (e2) {$\color{green}e_2$};
    \node[below] at (e0) {$e_0$};
     \node[below] at (e') {${ e'}$};
        \node[left] at (e1) {$e_1$};
            \node[circle, fill=red, inner sep=0.8pt] at (e1){};
\end{tikzpicture}
\begin{tikzpicture}
    \coordinate (e3) at (2,0,0);
    \coordinate (e1) at (0,2,0);
    \coordinate (e0) at (0,0,2);
    \coordinate (e2) at (1.2,1.4,1.4);
    \coordinate (e') at ($0.5*(e2)+0.5*(e3)$);
        \fill[color = blue, color = blue, opacity = 0.5] (e1) -- (e') -- (e3) -- cycle;
    % \fill[fill = violet, opacity = 0.6] (e0) -- (e1) -- (e2) -- cycle;
        \fill[color = blue, color = blue, opacity = 0.6] (e0) -- (e') -- (e3) -- cycle;
                \fill[color = blue, color = blue, opacity = 0.65] (e0) -- (e') -- (e1) -- cycle;
    \draw[thick, color = blue] (e1) -- (e3);
     \draw[thick, color = blue] (e0) -- (e3);
    \draw[thick, color = blue] (e') -- (e3);
    % \draw[thick, color = blue] (e0) -- (e2) -- (e3) -- cycle;
    \draw[thick, color = blue] (e0) -- (e');
    \draw[thick, color = blue] (e1) -- (e');
    % \draw[thick, color = violet] (e1) -- (e2);
    % \draw[thick, color = violet] (e0) -- (e2);
    \draw[thick, color = red] (e0) -- (e1);
    \node[circle, fill=red, inner sep=0.8pt] at (e0){};
        \node[circle, fill=red, inner sep=0.8pt] at (e1){};
       \node[circle, fill=blue, inner sep=0.8pt] at (e'){};
       \node[circle, fill=blue, inner sep=0.8pt] at (e3){};
           % \node[circle, fill=violet, inner sep=0.8pt] at (e2){};

    % \node[left] at (e2) {$\color{white} e_2$};
    \node[left] at (e1) {$\color{green} e_1$};
    \node[below] at (e0) {$e_0$};
    \node[below] at (e3) {$e_3$};
     \node[below] at (e') {${e'}$};
\end{tikzpicture}
    \caption{Depiction of $\Delta^{\J_k}, \Delta^{\J}$ and $\Delta^{\J'}$ with the vertices respectively opposite to $\Delta^{\J''_k}, \Delta^{\J''}$ and $\Delta^{\J'_k}$ marked in green.
    }
    \label{fig:three_horns}
\end{figure}
     \end{center}
    
     Next, consider the inclusions
     \begin{enumerate}
         \item $j_1 \colon \Lambda^{\J_k}_{s-1} \hookrightarrow \Delta^{\J_k}$;
         \item $j_2 \colon \Lambda^{\J}_{s} \hookrightarrow \Delta^{\J}$;
         \item $j_3 \colon \Lambda_k^{\J'} \hookrightarrow \Delta^{\J'}.$
     \end{enumerate}
     Gluing $\Delta^{\J_k}$ to $A$ along $j_1$ adds in the missing simplices $\Delta^{\J_k}$ and $\Delta^{\J''_k}$. Then, gluing in $\Delta^{\J}$ along $j_2$ adds the missing simplices $\Delta^{\J}$ and $\Delta^{\J''}$ (consider \cref{fig:three_horns}, for an illustration). Finally, gluing in $\Delta^{\J'}$ along $j_3$ adds the missing simplices $\Delta^{\J'}$ and $\Delta^{\J'_k}$. We have thus exposed $A \hookrightarrow B$ as a composition of the pushouts of horn inclusions $j_1$, $j_2$ and $j_3$. Since $\sReal{-}$ preserves pushouts, it suffices to show that $\tstr$ has the horn filling property with respect to (the realizations of) $j_1,j_2$ and $j_3$. Since $s \geq 2$, $j_1$ is an inner horn inclusion, and furthermore since $p_0 = p_k$ we have $t_{\J_k} < t_{\J}$, which is covered by the inductive assumption.
     $j_2$ is an inner horn inclusion, in which the $s$-th entry is not minimal, which we have already covered above through \cref{lem:metric_reduction_lemma}. Finally, $j_3$ is an inner horn inclusion with $o_{\J'} < o_{\J}$, and hence is also covered by the inductive assumption. This finishes the induction.     
    \end{proof}
    \subsection{Cofibrant stratified spaces and cellularly stratified spaces}
    Next, let us investigate the class of cofibrant objects in the stratum-preserving and poset-preserving semi-model categories, both categorical and diagrammatic.
    \begin{proposition}\label{prop:char_cofibrants_wo_ref}
        Let $\tstr \in \StratN$. Then the following statements are equivalent:
            \begin{enumerate}
                \item $\tstr$ is a retract in $\TopPN[\ptstr]$ of a cellularly stratified space;
                \item $\tstr$ is cofibrant in all of the semi-model categories $\TopPDN[\ptstr], \TopPCN[\ptstr], \StratDN, \StratCN$;
                \item $\tstr$ is cofibrant in one of the semi-model categories $\TopPDN[\ptstr], \TopPCN[\ptstr], \StratDN, \StratCN$.
            \end{enumerate}
    \end{proposition}
    \begin{proof}
        Note that there is no need to differentiate between the diagrammatic and categorical scenarios, as the latter are obtained from the former via left Bousfield localization (\cref{thm:huge_overview}), and hence both settings have the same cofibrant objects.
        That being a retract of a cellularly stratified space is equivalent to being cofibrant over $\ptstr$ is immediate by the cofibrant generators of \cref{thm:ex_of_model_struct_strat_nonref}.
        It remains to see that $\tstr$ being cofibrant in $\StratDN$ is equivalent to $\tstr$ being cofibrant in $\TopPDN[\ptstr]$. 
         This follows immediately from \cref{prop:model_global_mod_struct_are_glued}.
        % To see that every such retract is cofibrant, it suffices to see that cellularly stratified spaces are cofibrant in $\StratDN$. 
        % Even more, the generating cofibrations of $\TopPDN$ are all pushouts of the generating cofibrations in $\StratDN$.
        % Indeed, given any flag in $\J \in \catFlag$, consider the map $f \colon [n_\J] \to \pos$ defining the latter, and the induced stratified simplicial map $\partial \Delta^{[n_{\J}]} \to \Delta^\J.$ 
        % These maps fit into a pushout diagram in $\sStratN$
        % \begin{diagram}
        %     {}\partial \Delta^{[n_{\J}]} \arrow[d] \arrow[r, hook] &  \Delta^{[n_\J]} \arrow[d]\\
        %  \partial \Delta^\J\arrow[r, hook] & \Delta^\J;
        % \end{diagram}
        % the realization of which shows the claim.
        % Conversely, 
    \end{proof}
    \begin{definition}
    We call a stratified topological space $\tstr$ that satisfies any of the equivalent conditions of \cref{lem:equ_fibrancy_cond} \define{triangularly cofibrant.}
\end{definition}
     \begin{remark}
         Clearly every stratified space that admits a triangulation compatible with the stratification (i.e., in particular is the realization of a stratified simplicial set) is triangularly cofibrant. In particular, this holds for piece-linear pseudomanifolds (see, for example, \cite{banagl2007topological}) or locally compact stratified spaces which are definable (with definable stratification) in some o-minimal structure on the reals (\cite[Thm. 1.7]{van1998tame}). Furthermore, by \cite{TriangulationsGoresky}, all Whitney stratified spaces are in this class.
     \end{remark}
     However, one should not restrict oneself to stratified spaces that admit a cell structure. Topological manifolds, for example, are (to the best of our knowledge) not known to admit CW structures in dimension $4$ (see \cite{SiebenMannHauptvermutung} for all other dimensions). However, every topological manifold is a Euclidean neighborhood retract (\cite{HannerENR}) and thus a retract of a CW complex. Supposing the existence of certain stratified cylinder neighborhoods, our theory also covers examples of stratified spaces with manifold strata (see \cref{prop:inductive_cofibrancy}).
     Let us begin by proving that stratified neighborhood retracts of triangularly cofibrant stratified spaces are triangularly cofibrant:
     \begin{proposition}\label{prop:stratified_neighborhood_retracts}
        Let $\tstr \in \StratN$ be a triangularly cofibrant stratified space. Let $A \subset \tstr$ be a closed subspace such that the following holds: There is a neighborhood $U \subset \tstr$ of $A$, such that $\tstr[A] = \big (A, \sstr(U), \sstr \mid_A \colon A \to \sstr(U) \big)$ is a stratum-preserving  retract of $\tstr[U] = \big (U, \sstr(U), \sstr \mid_U \colon U \to \sstr(U) \big)$.
        Then $A$ is also triangularly cofibrant.
\end{proposition}
\begin{proof}
    We will make frequent use of results on stratified cell complexes developed in \cite{HoLinksWa}.
    Let $\tstr[Y]$ be a cellularly stratified space and $i \colon \tstr \hookrightarrow \tstr[Y]$, $r \colon \tstr[Y] \to \tstr$ such that $r \circ i = 1_{\tstr}$. We claim that there is a subcomplex $\tstr[B]$ of a cell structure on $\tstr[Y]$ that contains $i(A)$ and is entirely contained in $r^{-1}(U)$. Then we can compose
    \[
    \tstr[B] \to \tstr[U] \to \tstr[A],
    \]
    where $\tstr[B]$ is considered a stratified space over $\ptstr[U]$, to obtain a retraction of $\tstr[A] \hookrightarrow \tstr[B]$. Fix a stratified cell structure $(\sigma_{i} \colon \sReal{\Delta^{\J_{i}}} \to \tstr[Y])_{i \in I}$ on $\tstr[Y]$ (see \cite[Subsec. \ref{hol:subsec:strat_cell_complex}]{HoLinksWa}, for a definition and basic properties). Every cell in $\tstr[Y]$ defined by a map $\sigma_{i}$ is contained in a finite subcomplex of $\tstr[Y]$ (\cite[Lem. \ref{hol:lem:finite_subcomplex}]{HoLinksWa}). We proceed to inductively construct a subdivision of $\tstr[Y]$ (in the sense of (in the sense of \cite[Def. \ref{hol:def:subdivision}]{HoLinksWa}), via induction over the minimal number of cells, $n_{i}$, required in a subcomplex that contains a cell $\sigma_{i}$. 
    Denote by $\tstr[Y]^{(n-1)}$ the subcomplex constructed from cells $\sigma_i$ with $n_{i} = n$.
    We construct a subdivision of the cell structure of $\tstr[Y]$, denoted $\tstr[Y]',$ through induction on $n$.
    For $n=1$, we simply take the cell structure on $\tstr[Y]^{(0)}$ inherited from $\str[Y]$. 
    Now, suppose that we have already defined a cell structure $\tstr[Y]^{'(n)}$ on $\tstr[Y]^{(n)}$ for $n>1$, and a subcomplex $\tstr[B]^n \subset \tstr[Y]^{(n)}$ of this cell structure, fulfilling the following properties:
            \begin{enumerate}
                \item $r^{-1}(A) \cap \tstr[Y]^{(n)} \subset \tstr[B]^n$ and $\tstr[B]^n$ is a $\Delta_P$-neighborhood of $r^{-1}(A) \cap \tstr[Y]^{(n)}$ (in the sense of \cite[Def. \ref{hol:def:DeltaP_neighborhood})]{HoLinksWa};
                \item $\tstr[B]^n$ is contained in the interior of $r^{-1}(U)$.
            \end{enumerate}
    Next, let us construct $\tstr[Y]^{'(n+1)}$ and $\tstr[B]^{n+1}$, i.e., we now add cells $\sigma_i$ with $n_i = n+2$ to $\tstr[Y]^{(n)}$.
    For any cell $\sigma_i$ with $n_i = n+2$, we can barycentrically subdivide $\sReal{\Delta^{\I_i}}$ (sufficiently many times) so that for every $\tau \subset \sReal{\Delta^{\J_i}}$ that is contained in a simplex $\tau'$ intersecting $(r \circ \sigma)^{-1}(A)$ we have the following. \begin{enumerate}
        \item $\tau$ is entirely contained in the interior of $(r \circ \sigma)^{-1}(U)$;
        \item $\tau \cap \sReal{ \partial \Delta^{\I_i}}$ is entirely contained in $\sigma^{-1}(\tstr[B]^{n})$.
    \end{enumerate}
    Indeed, the first property is possible, since $(r \circ \sigma)^{-1}(A)$ is compact and contained in the interior of $(r \circ \sigma)^{-1}(U)$ by (inductive) assumption. It follows that $(r \circ \sigma)^{-1}(A)$ has a positive distance (for any choice of compatible metric on $\sReal{\Delta^{\J_i}}$), to the complement of the interior of $(r \circ \sigma)^{-1}(U)$, which implies that simplices of sufficiently small diameter cannot intersect $(r \circ \sigma)^{-1}(A) \cup \sigma^{-1}( \tstr[B]^n)$ and the complement of the interior of $(r \circ \sigma)^{-1}(U)$.
    That the second property can be assumed is argued similarly. 
    Having chosen such a subdivision for any cell $\sigma_i$ with $n_i = n+2$, we then obtain an induced cell structure $\tstr[Y]^{'(n+1)}$ on $\tstr[Y]^{(n+1)}$ given by the cells of $\tstr[Y]^{'n}$ and the cells corresponding to (open) simplices in the interior of the subdivisions of $\sReal{\Delta^{\I_i}}$, for $n_i = n+2$.
    Finally, we let $\tstr[B]^{n+1}$ be the subcomplex given by adding to $\tstr[B]^{n}$ all such cells of $\tstr[Y]^{'(n+1)} \setminus \tstr[Y]^{'(n)}$, $\sigma_i$, which correspond to an (open) simplex in the interior of some $\sReal{\Delta^{\I_i}}$, contained in a closed simplex intersecting $(r \circ \sigma)^{-1}(A)$. 
    Note that this does indeed define a cell complex, since the intersection of the boundary of such a simplex with the boundary of $\sReal{\Delta^{\I_i}}$ is assumed to be mapped to $\tstr[B]^{n}$. 
    $\tstr[B]^{n+1}$ defined in this fashion fulfills the inductive assumptions by construction and \cite[Lem. \ref{hol:lem:char_of_nbhd_in_cell_complex}]{HoLinksWa}. Then, finally, set $\tstr[B] = \bigcup_{n \in \mathbb N} \tstr[B]^{n}$, with the induced cell structure.
\end{proof}

     \begin{construction}[Stratified mapping cylinders]
    Given $f \colon \tstr \to \tstr[Y]$ in $\StratN$, the \define{stratified mapping cylinder} of $f$, $M_{s}(f)$, is defined as the pushout
        \begin{diagram}
            \tstr \times \sReal{\Delta^{0}} \arrow[r, hook, "1 \times i_0"]  \arrow[d, "f"]& \tstr \times \sReal{\Delta^{[1]}} \arrow[d] \\ 
            \tstr[Y] \arrow[r, hook] & M_{s}(f).
        \end{diagram}
    Under the (nonstratified) identification $\sReal{\Delta^{[1]}} = [0,1]$, this equips the classical mapping cylinder of $f$ with an alternative stratification over the lower right corner of the following pushout diagram of posets.
    \begin{diagram}
            \ptstr   \arrow[r, hook, "1 \times i_0"]  \arrow[d, "f"] & \ptstr \times [1] \arrow[d] \\ 
            \ptstr[Y] \arrow[r, hook] & \ptstr[Y] \sqcup_{\leq_f} \ptstr \spaceperiod
        \end{diagram}
    Here, $\ptstr[Y] \sqcup_{\leq_f} \ptstr$ is explicitly given by equipping the disjoint union $\ptstr[Y] \sqcup  \ptstr  $, with an additional relation $p \leq q$, for $p \in \ptstr[Y]$ and $q \in \ptstr[X]$, whenever $p \leq f(q)$. In the special case where $\ptstr[Y]$ is a poset with one element, this amounts to adjoining a minimal element to $\ptstr$.
\end{construction}
\begin{remark}
    Note that $M_{s}(f)$ is generally not a mapping cylinder of $f$ with respect to any of the model structures on $\StratN$. Indeed, the inclusion $\tstr[Y] \hookrightarrow M_{s}(f)$ can only be a weak equivalence (in any of the model structures on $\StratN$) if $\tstr$ is empty. From the perspective of $\infty$-categories of exit paths, no point $x \in \tstr \times \sReal{\Delta^{\{1\}}}\subset M_{s}(f)$ can lie in the essential image corresponding to the map $\tstr[Y] \hookrightarrow{} M_{s}(f)$, as this would imply a non-stratified path from a point $y \in \tstr[Y]$ to $x$.
\end{remark}
% \begin{proposition}\label{prop:prod_is_bifibr}
%     Let $\TopN$ (as in \Lukas{add ref}) be cartesian closed. Furthermore, consider $\StratN$ as equipped with any of the model structures in \Lukas{add ref} and let $\sStratN$ be equipped with its simplicial counterpart.
%     The functor 
%     \[
%     - \times \sReal{-} \colon \StratN \times \sStratN \to \StratN 
%     \]
%     is a Quillen bifunctor of (semi)model categories, for any of the semi-model structures in \Lukas{add ref}. That is:
%     \begin{enumerate}
%         \item For $i \colon \tstr \hookrightarrow \tstr[Y] \in \StratN$, $j \colon A \hookrightarrow B \in \sStratN$ cofibrations the induced pushout product
%         \[
%         i \boxtimes j \colon  \tstr[Y] \times \sReal{A} \cup_{\tstr[X] \times \sReal{A}} \tstr[X] \times \sReal{B} \hookrightarrow \tstr[Y] \times \sReal{B} 
%         \]
%         is a cofibration;
%         \item If in addition, $i$ has cofibrant source, and either $i$ or $j$ is a weak equivalence, then $i \boxtimes j$ is a weak equivalence.
%     \end{enumerate}
% \end{proposition}
\begin{corollary}\label{cor:map_cyl_res}
    Suppose that $\tstr[L], \tstr[X], \tstr[Y] \in \StratN$ are triangularly cofibrant. Let $f \colon \tstr[L] \to \tstr [Y], g \colon \tstr [L] \to \tstr$ be stratified maps. Then the stratified double mapping cylinder $M_{s}(f) \cup_{g} \tstr $ is triangularly cofibrant.  
\end{corollary}
\begin{proof}
    $M_{s}(f) \cup_{g} \tstr $ fits into a pushout diagram of stratified spaces
        \begin{diagram}
            \tstr [L] \times \sReal{\partial \Delta^{[1]}} = \tstr [L] \sqcup \tstr [L] \arrow[d, "f\sqcup g"]\arrow[r, hook] & \tstr [L] \times \sReal{ \Delta^{[1]}} \arrow[d] \\
            \tstr [Y] \sqcup \tstr [X] \arrow[r, hook] & M_{s}(f) \cup_{g} \tstr \spaceperiod
        \end{diagram}
    Note that we may, without loss of generality, assume that we are in the $\Delta$-generated scenario, as every triangularly cofibrant stratified space is $\Delta$-generated. By \cref{cor:cart_closed} the upper horizontal is a cofibration, from which it follows that the lower horizontal is a cofibration. $\tstr [Y] \sqcup \tstr [X]$ is cofibrant by assumption, which, together with the cofibrancy $\tstr [Y] \sqcup \tstr [X] \hookrightarrow $ of $M_{s}(f) \cup_{g} \tstr$, implies the claim.
\end{proof}
\begin{definition}
    Let $U \subset \tstr$ and $A \subset U$ be a closed subset contained in the interior of $U$. Denote $\tstr[A] = ( \sstr \mid_A \colon A \to \sstr(A) )$ and $\tstr[U] = ( \sstr \mid_{U} \colon U \to \sstr(U) )$.
    We say that $U$ is a \define{(closed) stratified mapping cylinder neighborhood of $A$}, if the following holds: \\
    There exists a stratified space $\tstr[L]$ together with a stratified map $f \colon \tstr[L] \to \tstr[A]$, as well as a stratum-preserving homeomorphism $\phi \colon M_s(f) \hookrightarrow  \tstr[U]$ under $\tstr[A]$, such that $\phi^{-1}( \partial U) = \tstr[L] \times \{1\}$.
\end{definition}
\begin{definition}
    Let $\pos$ be a poset and $p \in \pos$. We call
    \[
    \sup \{ n \in \mathbb N \mid \exists p_0, \cdots, p_n \colon p = p_0 < p_1 < \cdots < p_n \}
    \]
    the \textit{depth} of $p$ in $\pos$. We call the supremum over the depths of all $p \in \pos$ the \define{depth} of $\pos$. Given $k \in \mathbb N$ and $\tstr \in \StratN$, we denote by $\tstr_{k}$ the stratified subspace given by restricting $\tstr$ along \[
    (\ptstr)_k:= \{ p \in \pos \mid p \textnormal{ has depth $k$} \} \hookrightarrow \pos. 
    \]
    Similarly, we denote by $\tstr_{\underline{k}}$ the stratified subspace of $\tstr$ obtained by restricting $\tstr$ along \[
    (\ptstr)_{\underline k}:=\{ p \in \pos \mid p \textnormal{ has depth $\leq k$} \} \hookrightarrow \pos. 
    \]
\end{definition}
\begin{proposition}\label{prop:inductive_cofibrancy}
    Let $\tstr$ be a stratified space with $\pstr$ finite depth. Suppose that the following holds: 
    \begin{enumerate}
        \item Each stratum of $\tstr$ is cofibrant in the Quillen model structure, i.e. a retract of a non-stratified absolute cell complex. 
        \item For each $k \in \mathbb N$, there is a family of pairwise disjoint subsets of $\tstr$, $(U^p)_{p\in (\ptstr)_k}$, indexed over strata of depth $k$ such that, for each $p$, $U^p$ is a stratified mapping cylinder neighborhood of $\tstr_p$.
    \end{enumerate}
    Then $\tstr$ is triangularly cofibrant.
\end{proposition}
\begin{proof}
    For ease of notation, we denote $\pos := \ptstr$.
    We proceed via induction over the depth of $\pos$, denoted $d$. In the case $d=0$, $\tstr$ is simply a disjoint union of trivially stratified spaces, each of which is cellularly stratified by assumption. Now, for the inductive step $d$ to $d+1$: By assumption, there are stratified maps $f^p \colon \tstr[L]^{p} \to \tstr_p $, for each $p \in \pos_{d+1}$, together with an injective (on the space level) stratified map 
    \[
    \bigsqcup_{p \in (\pos_{d+1})} M(f^p)  \hookrightarrow \tstr, \\
    \]
    which restricts to an open inclusion on the open cylinders (obtained by removing $\str L^{p} \times \{1\}$ and denoted by $\mathring{M}(f^p)$). Setting, $\tstr[L] = \bigsqcup_{p \in \pos_{d+1}} \tstr[L]_p$ and $f = \bigsqcup_{p \in \pos_{d+1}}f^p$, we obtain a stratified open inclusion
    \[
    i \colon \bigsqcup_{p \in (\pos_{d+1})} \mathring M(f^p) \cong \mathring M_s(f) \hookrightarrow \tstr,
    \]
    which defines a neighborhood of $\tstr_{d+1}$. Note that $i$ is not necessarily an inclusion on the poset level. 
    Now, consider $\mathring M_s(f)$ as reparametrized (over $[0,2)$), denoted $\mathring M'_s(f)$, such that we may consider $M_s(f)$ (with the usual parametrization) as a closed stratified subspace of $\mathring M_s(f)$. 
    With this new notation, we have inclusion $\mathring M_s(f) \hookrightarrow M_s(f) \hookrightarrow \mathring M'_s(f) \hookrightarrow \tstr[X]$.
    Then $\mathring M_s(f) \subset \mathring M'_s(f)$ is an open neighborhood of $\tstr_{d+1}$ with boundary $\tstr[L] \times \{1\}$. 
    We obtain a commutative diagram of stratified maps
   \begin{diagram}\label{diag:double_mapping_inductive}
       \tstr[L] \times \{1\}  \arrow[d, hook] \arrow[r, "g"]& \tstr \setminus \mathring M_s(f) \arrow[d]\\
       M_s(f) \arrow[r, hook] & \tstr \spacecomma
   \end{diagram}
   where $ \tstr \setminus \mathring M_s(f) $ is stratified over $\pos_{\underline{d}}$.
   Denote by $\tilde \tstr$ the pushout in $\StratN$. Now, on the level of topological spaces, this diagram is clearly pushout. Therefore, the induced map $g \colon \tilde \tstr \to \tstr$ is a homeomorphism on the underlying spaces. Consider the commutative diagram 
   \begin{diagram}
       \tilde \tstr \arrow[r] \arrow[rd] & P(g)_! \tilde{\tstr} \arrow[d] \\
                              &  \tstr \spaceperiod
   \end{diagram}    
    Since $g$ is a homeomorphism on the underlying spaces, it follows that the right vertical is an isomorphism. The upper vertical is always a cofibration (\cref{prop:model_global_mod_struct_are_glued}), hence it follows that $\tilde \tstr \to \tstr$ is a cofibration. Therefore, it suffices to show that $\tilde{ \tstr}$ is cofibrant. For this, in turn, it suffices to show that the left-hand vertical is a cofibration and $\tstr \setminus \mathring M_s(f)$ is cofibrant. The latter is a retract of $\tstr_{\underline d}$. Hence, cofibrancy follows by inductive assumption. For the left vertical, by \cref{cor:map_cyl_res}, it suffices to see that $\tstr[L]_p$ and $\tstr_p$ are cofibrant. The latter is cofibrant by assumption.
   % On the level of posets, the diagram looks as follows 
   % \begin{diagram}
   %    \ptstr[L]  \arrow[d, hook] \arrow[r, hook]& (\ptstr)_{\underline d} \arrow[d]\\
   %    (\ptstr)_{d+1} \sqcup_{\leq_f} \ptstr[L] \arrow[r, hook ] & \ptstr[X].
   % \end{diagram}
   % On the level of sets this diagram is clearly pushout. On the level of posets, it suffices to see that every relation in $\ptstr[X]$ $p \leq q$ with $p \in (\ptstr)_{d+1}$, $q \in (\ptstr)_{\underline d}$ can be written in the form $(\ptstr)_{d+1} = f(q')$, $g(q') \leq q$. Indeed, this is the case since we assumed the minimal elements $q$ over $p$, to admit $q' \in \ptstr[L]$ with $g(q') = q$ and $f(q') = p$. Consequently, we have exposed $\tstr$ as a double mapping cylinder of $f$ and $g$. Using \cref{cor:map_cyl_res}, it suffices to show that $\tstr[L], \tstr_{d+1}$ and $ \tstr \setminus \mathring M_s(f)$ are retracts of cellularly stratified spaces. For $\tstr_{d+1}$ this holds by the case $d=0$. For $\tstr[L]$, note that $\tstr[L] = \bigsqcup_{p \in (\ptstr)_{d+1}} \tstr[L]^p$ and it suffices to show the claim separately for each component. 
   % We may without loss of generality assume that $\tstr[L]$ has no empty strata (since postcomposing with a poset map preserves being a retract of a cellularly stratified space by \Lukas{add ref}),
   By construction, $\tstr[L]^p$ embeds into $\tstr$ with an open (stratified) neighborhood of the form $\tstr[L]^p \times (0,1)$. Clearly, $\tstr[L]^p$ is a retract of the latter, making $\tstr[L]^p$ a stratified neighborhood retract of $\tstr_{\underline{d}}$, which is cofibrant by inductive assumption. Hence, the case of $\tstr[L]^p$ follows by \cref{prop:stratified_neighborhood_retracts}.
 \end{proof}
\begin{remark}
Every topological manifold is a Euclidean neighborhood retract (\cite{HannerENR}) and since every open subset of Euclidean space can be triangulated, it follows that every topological manifold is cofibrant in the Quillen model structure. Consequently, it follows by \cref{prop:inductive_cofibrancy} that every stratified space with manifold strata, each of which admits appropriate stratified mapping cylinder neighborhoods, is cofibrant. \cite[Prop. 8.2.3]{LocalStructOnStrat} asserts the existence of stratified mapping cylinder neighborhoods for conically smooth stratified spaces, which would make the latter triangularly cofibrant. 
More generally, for homotopically stratified spaces with manifold strata there are obstructions to the existence of (pairwise) stratified mapping cylinder neighborhoods (\cite[Thm. 1.7]{quinn1988homotopically}).
\end{remark}
\subsection{Refined stratified spaces}\label{subsec:refined_strat_spaces}
Let us take a more detailed look at the cofibrant objects in the semi-model categories $\StratDR$ and $\StratCR$. In particular, our aim is to relate them to more classical properties of stratified spaces.
 \begin{definition}\label{def:ref_strat_sp}
        A stratified space $\str \in \StratN$ is called \define{refined} if the following holds:
        \begin{itemize}
            % \item $\str$ carries the $\Delta_s$ topology.
            \item $\str$ is \define{surjectively stratified}, that is, $\sstr \colon \utstr \to \ptstr$ is surjective.
            \item For each pair $x,y\in \str$, the relation $\sstr(x) \leq \sstr(y)$ holds if and only if there is a finite sequence of stratified maps $\gamma_i \colon \sReal{\stratSim[1]} \to \str$, $i=1, \dots, n$, with $\gamma_1(0) = x$, $\gamma_n(1) =y$ and $\gamma_{i}(1) = \gamma_{i+1}(0)$, for $i < n$. 
        \end{itemize}
\end{definition}
We can think of being refined as the poset $\ptstr$ being completely reflected in the stratified paths of $\ptstr$. 
Next, note the following elementary properties about refinedness.
\begin{proposition}\label{prop:basic_props_ref}
Let $\tstr,\tstr[Y] \in \StratN$ and $\str[A] \in \sStratN$. Then the following holds: 
    \begin{enumerate}
        \item $\tstr$ is refined, if and only if $\SingS (\tstr)$ is refined;
        \item If $f \colon \tstr \to \tstr[Y]$ is a poset-preserving categorical equivalence, then $\tstr$ is refined, if and only if $\tstr[Y]$ is refined;
        \item $\str[A]$ is refined, if and only if $\sReal{\str[A]}$ is refined;
        \item If $\tstr[A] \in \sStratCN$ is fibrant, i.e. $\ustr[A]$ a quasi category and $\ustr[A] \to \pos$ a conservative functor, then $\str[A]$ is refined if and only if $\str[A]$ is $0$-connected, in the sense of \cite{Exodromy}, as an abstract stratified homotopy type.
    \end{enumerate}
\end{proposition}
It follows that being refined is really purely a property of the homotopy type in $\AltStratD$ defined by $\tstr$.
% We first need the following lemma:
% As a consequence, corollary of \cref{thm:overview_over_all_hypothesis} holds:
% \begin{lemma}
%     \label{lem:ref_inv_under_real}
%     Let $\str[A] \in \sStrat$ be a stratified simplicial set. Then $\str[A]$ is refined, if and only if $\sReal{\str[A]}$ is refined.
% \end{lemma}
% \begin{proof}
%     Note that being refined is a homotopy theoretic property in $\sStratC$, as it is defined through a derived functor, $\rpstr: \sStrat \to \Pos$. It thus follows by \cref{thm:overview_over_all_hypothesis} that $\str[A]$ is refined, if and only if $\SingS \sReal{\str[A]}$ is refined. Since by definition $\sReal{A}$ carries the $\Delta_s$-generated topology, it follows by \cref{prop:equ_char_of_ref_top}, that  $\SingS \sReal{\str[A]}$ is refined if and only if $\sReal{\str[A]}$ is refined.
% \end{proof}
The model categories $\StratCRN$ and $\StratDN$ may now be interpreted as the respective right Bousfield localizations presenting the full sub-$\infty$-categories of refined stratified spaces.
Even more, we may essentially construct the left adjoints of this adjunction on the $1$-categorical level.
 \begin{definition}\label{def:ref_topology}
        A stratified space $\tstr[X] \in \StratN$ is said to carry the $\Delta_s$ topology if one of the following equivalent conditions holds:
        \begin{enumerate}
            \item $\utstr[X]$ has the final topology with respect to the set of stratified maps $\sReal{\stratSim[n]} \to \tstr[X]$, for $n \in \mathbb N$.
            \item $\utstr[X]$ has the final topology with respect to the set of stratified maps $\sReal{\stratSim[1]} \to \tstr[X]$.
        \end{enumerate}
    \end{definition}
    \begin{proof} Let us show that these conditions are equivalent.
  We need to see that every stratified simplex $\sReal{\stratSim[n]}$, for $n \in \mathbb N$, carries the final topology with respect to stratified maps with source $\sReal{\stratSim[1]}$.
  Note that it suffices to see that for any convergent sequence in $\sReal{\stratSim[n]}$, there exists a subsequence $(x_i)_{i \in \mathbb N}$, as well as a stratified map $\gamma \colon [0,1] \cong \sReal{\stratSim[1]} \to \sReal{\stratSim[n]}$, with $\gamma(\frac{1}{2^i}) = x_i$, for $i \in \mathbb N$. Indeed, the topology on $\sReal{\stratSim[n]}$ is entirely determined by convergent sequences, and the latter conditions mean that these are detected by $\sReal{\stratSim[1]}$. So suppose that we are given a convergent sequence $(\hat x_i)_{i \in \mathbb N}$ in $\sReal{\stratSim[n]}$. By passing to a subsequence, we may assume that the sequence is contained in a single stratum $j \in [n]$. Then, the limit point lies in some stratum $k \leq j$. 
    Now define $\gamma \colon [0,1] \cong \sReal{\stratSim[1]} \to \sReal{\stratSim[n]}$ by setting $\gamma(\frac{1}{2^i}) := x_i$, convexly interpolating between $x_i$ and $x_{i+1}$ and sending $0$ to the limit point of $x_i$. This map is continuous. Furthermore, as the strata of $\sReal{\stratSim[n]}$ are convex, it is also stratified.
    \end{proof}
    \begin{definition}\label{def:strong_ref_strat_sp}
        A stratified space $\str \in \StratN$ is called \define{strongly refined}, if $\tstr$ is refined and carries the $\Delta_s$-topology. 
    \end{definition}
    % \begin{remark}\label{rem:equ_cond_for_ref_pos}
    %     Note that the third condition in \cref{def:ref_strat_sp} may equivalently be replaced by all strata being path connected, and there being a relation $\sstr(x) < \sstr(y)$, if and only if there is a finite sequence of stratified maps $\gamma_i \colon \sReal{\stratSim[1]} \to \str$, $i=1, \dots, n$, with $\gamma_1(0) = x$, $\gamma_n(1) =y$ and $\gamma_{i}(1) = \gamma_{i+1}(0)$, for $i < n$, at least one of which is not constant on strata. 
    % \end{remark}
    \begin{construction}\label{con:top_refinement}
        For a stratified space $\tstr[X] \in \StratN$, its \define{refined poset}, denoted $\rpstr$, is the poset generated by the following elements and relations:
        \begin{itemize}
            \item An element for each $x \in \tstr[X]$.
            \item A relation $x \leq y$, whenever there is a sequence of stratified paths from $x$ to $y$, as in \cref{def:ref_strat_sp}.
        \end{itemize}
        Equivalently, $\rpstr$ is given by the set of path components of strata of $\tstr[X]$, together with a generating relation whenever there is an exit path from one component to another.
        Consider the (generally non-continuous) map $\utstr[X] \to \rpstr$, given by $x \mapsto [x]$. It provides a factorization of $\ststr[X]$ through the map of posets \begin{align*}
            \rpstr &\to \ptstr[X] \\
            [x] &\mapsto \ststr[X](x)
        \end{align*}
        as follows:
        \begin{diagram}
            {\utstr[X]} \arrow[rd, "{\ststr[X]}"'] \arrow[r, dashed]  & \rpstr \arrow[d] \\
            & {\ptstr[X]}
        \end{diagram}
        Although $\utstr[X] \to \rpstr[X]$ is not necessarily continuous, its precomposition with any stratified map $\sReal{\stratSim[1]} \to \tstr[X]$ is continuous, by construction.
        We denote by $\tstr[X]^{\ared}$ the stratified space given by $\utstr[X^{\ared}] \to \rpstr[X]$, where $\utstr[X^{\ared}]$ has the same underlying set as $\utstr[X]$
        and is equipped with the final topology with respect to stratified maps $\sReal{\stratSim[1]} \to \str[X]$. By construction, $\utstr[X^{\ared}] \to \rpstr[X]$ is indeed continuous and $\tstr[X]^{\ared}$ is a strongly refined stratified space. This construction induces a simplicial functor from $\StratN$ into the full subcategory of strongly refined stratified spaces, called \define{the refinement functor}, which exposes the latter as a full coreflective subcategory of $\StratN$. The counit of adjunction, given by the commutative squares
        \begin{diagram} 
          {  \utstr[X^{\ared}]}  \arrow[d, "{\ststr[X^{\ared}]}"] \arrow[r]& {\utstr[X]} \arrow[d, "{\ststr[X]}"]             \\
            {\rpstr[X]} \arrow[r] & {\ptstr[X]} \spacecomma
        \end{diagram}
         is called the \define{refinement map}. 
    \end{construction}
    As a consequence of the existence of the refinement functor, which makes the inclusion of strongly refined stratified spaces a left adjoint, we obtain:
    \begin{corollary}\label{cor:colims_of_refs_are_ref}
        Every colimit of strongly refined stratified spaces in $\StratN$ is again strongly refined.
    \end{corollary}
    The refinement construction of \cref{con:top_refinement} is really the topological analogue of \cite[Def. \ref{comb:def:refined}]{ComModelWa}.
    \begin{proposition}[See {\cite[Subsec. \ref{comb:subsec:refining_strat_sset}]{ComModelWa}}, for notation]\label{prop:ref_com_sings}
        For any $\tstr \in \StratN$, applying $\SingS$ to the refinement map induces a natural transformation (dashed) that makes the diagram.
        \begin{diagram}
            (\SingS \tstr[X])^\ared \arrow[r, dashed] \arrow[rd] & \SingS ( \tstr[X]^{\ared}) \arrow[d]  \\
            {}& \SingS (\tstr[X])
        \end{diagram}
        commute. The dashed transformation is an isomorphism.
    \end{proposition}
    \begin{proof}
        $\SingS \tstr[X]^{\ared}$ is refined by \cref{prop:basic_props_ref}. Hence, the dashed map is induced by the fact that (simplicial) refinement is right adjoint to the inclusion of refined stratified simplicial sets. That this map is an isomorphism on posets follows immediately by the construction of $\pstr[X^{\ared}]$ in \cref{con:top_refinement} and \cite[Prop. \ref{comb:prop:explicit_rp}]{ComModelWa}. Furthermore, since every stratified simplex $\sReal{ \stratSim[n]}$ is refined, note that the map is an isomorphism on simplicial sets, given by
        \[
        (\SingS \tstr[X])^\ared [n] = (\SingS \tstr) = \StratN(\sReal{ \stratSim[n]}, \tstr[X]) =  \StratN(\sReal{ \stratSim[n]}, \tstr[X]^\ared ) = \SingS (\tstr[X]^{\ared})[n].
        \]
    \end{proof}
        We may think of \cref{prop:ref_com_sings} as stating that $\SingS$ sends the topological refinement map to the simplicial refinement map. As an immediate corollary, using \cref{rec:ref_mod_structs}, we obtain the missing part of \cref{thm:overview_over_all_hypothesis}:
    \begin{corollary}
        The simplicial semi-model categories $\StratDR$ and $\StratCR$ are obtained, respectively, from $\StratD$ and $\StratC$ by right Bousfield localizing at the refinement maps $\tstr{}^{\ared} \to \tstr$.
    \end{corollary}
In particular, we have the following description of cofibrant objects in $\StratDRN$ and $\StratCRN$:
\begin{proposition}\label{prop:rel_cof_and_ref}
    Let $\tstr \in \StratN$. Then the following statements are equivalent:
    \begin{enumerate}
        \item $\tstr$ is cofibrant in $\StratDRN$;
        \item $\tstr$ is cofibrant in $\StratCRN$;
        \item $\tstr$ is triangularly cofibrant and refined;
        \item $\tstr$ is a retract of a refined, cellularly stratified space.
    \end{enumerate}
\end{proposition}
\begin{proof}
    % Note first that every cofibrant object of $\StratC$ carries the $\Delta_s$-topology, as by \cref{prop:cof_fib_nonref_top,rem:cell_are_delta_gen} it is a retract of a stratified space with the $\Delta_s$-topology. Hence, in all of the characterizations above, the spaces carry the $\Delta_s$-topology, and we really only need to verify that the remaining conditions of \cref{def:ref_topology} are fulfilled.
    The equivalence between the first two statements is immediate from both semi-model categories having the same generating cofibrations.
    Since every cellularly stratified space carries the $\Delta_s$-topology (by definition), it follows that a refined cellularly stratified space is strongly refined. Since the latter form a full coreflective subcategory, by \cref{con:top_refinement}, any retract of strongly refined spaces is strongly refined. In particular, the fourth statement implies the third.
    % Since being refined and being cofibrant is stable under retracts (as a consequence of the characterization in terms of the refinement map being an isomorphism in \cref{prop:equ_char_of_ref_top}), it follows that the fourth condition implies the third. 
    % Since 
    % By \cref{rem:ref_can_det_on_sing_lvl}, the third condition implies the second. 
    To see that the third implies the second, note that cofibrant objects $\tstr$ in $\StratCRN$ are characterized by $\emptyset\hookrightarrow \tstr$ having the left lifting property, with respect to all stratified maps $f \colon \tstr[Y] \to {\tstr[Z]}$, which are acyclic fibrations in $\StratCRN$. By \cref{prop:ref_com_sings} and \cite[Thm. \ref{comb:prop:ex_red_struct}]{ComModelWa}, this, in turn, is equivalent to $f^{\ared}$ being an acyclic fibration in $\StratCN$. Since $\tstr$ is assumed to be refined, a lifting diagram
    \begin{diagram}
      {} & {\tstr[Y]} \arrow[d, "f"]\\
        \tstr \arrow[ru, dashed]   \arrow[r]& {\tstr[Z]}
    \end{diagram}
    admits a solution if and only if the induced diagram
     \begin{diagram}
       {}& {\tstr[Y]^{\ared}} \arrow[d, "f^{\ared}"]\\
        \tstr \arrow[ru, dashed]   \arrow[r]& {\tstr[Z]}^{\ared}.
    \end{diagram}
    admits a solution. The latter holds, as $\tstr$ was assumed cofibrant in $\StratCN$. Finally, to see that the first characterization implies the fourth, note that every cofibrant object in $\StratCRN$ is a retract of an absolute cell complex with respect to the stratified boundary inclusions $\{\sReal{\stratBound \hookrightarrow \stratSim}  \mid n \in \mathbb N\}$. It follows by \cref{prop:basic_props_ref}, stratified simplices carrying the $\Delta_s$-topology and \cref{cor:colims_of_refs_are_ref}, that every such absolute cell complex is a refined cellularly stratified space.
\end{proof}
\subsection{Frontier conditions and refinement}\label{subsec:frontier_cond_and_refinement}
It turns out that being refined is strongly related to the way the poset structure on the strata of a stratified space is classically constructed:
\begin{recollection}\label{recol:classical_strat}
    Classically, stratifications often arise from the so-called \define{frontier condition} (see, for example, \cite{mather1970notes}). Namely, one starts with a topological space $X$ and a locally finite decomposition into nonempty locally closed pieces $(X_i)_{i \in I}$. One assumes that for any $i \in I$, the closure $\overline{X_i}$ is given by the disjoint union $(X_j)_{j \in J}$, for some subset $J \subset I$. Then, $I$ naturally carries the structure of a poset, setting $i 
 \leq j$, whenever $X_i \subset \overline{X_j}$, and $X \to I$  defines a stratification of $X$. 
\end{recollection}
In our framework, such poset structures induced by closure relations can be constructed as follows.
\begin{construction}\label{con:top_strat}
     Let $\tstr \in \StratN$, and denote \[
    I = \{ S \subset \utstr \mid S \neq \emptyset \exists p \in \ptstr \colon S = \tstr_p.\}
    \]
    We consider $I_{\tstr}$ as equipped with the structure of a poset, by equipping it with the relation generated by 
    \[
    S \leq S' \iff S \cap \overline{S} \neq \emptyset.
    \]
    It follows immediately by construction that mapping $S$ to the unique $p \in \ptstr$, for which $S = \tstr_p$, induces a map of posets \[
    I_{\tstr} \to \ptstr.
    \]
\end{construction}
\begin{definition}
We say that $\tstr$ is \define{weakly frontier stratified}, if the induced map $I_{\tstr} \to \ptstr$ is an isomorphism of posets. We say that $\tstr$ is \define{frontier stratified}, if in addition to this it fulfills \[
\tstr_p \cap \overline{\tstr_q} \neq \emptyset \implies \tstr_p \subset \overline{\tstr_q},
\]
for all $p,q \in \ptstr$.
\end{definition}
\begin{remark}\label{rem:different_notions_of_stratified_homs}
    If one describes stratified spaces $\mathcal{X}$ and $\str[Y]$ as the data of spaces $X$ and $Y$ equipped, respectively, with decompositions $(X_i)_{i \in I}$ and $(Y_j)_{j \in J}$ into non-empty pieces, then classically the morphisms which are considered between such spaces are given by continuous maps $f \colon X \to Y$, such that for each $i\in I$ there exists a $j\in J$ with $f(X_i) \subset Y_j$ \cite{HughesPathSpaces}.
    Let us call such objects decomposition spaces, and such maps decomposed maps, and denote the corresponding category by $\cat[D]$.
    Homotopies in this setting are defined through the cylinder given by equipping $X \times [0,1]$ with the decomposition $(X_i \times [0,1])_{i \in I}$.
    There is an obvious forgetful functor from the category $\mathcal{D} \colon \StratN \to \cat[D]$, given by equipping a stratified space with its decomposition into nonempty strata. Now, if $\tstr$ is weakly frontier stratified, then it is not hard to see that for any $\tstr[Y] \in \StratN$ the induced map
    \[
    \StratN( \tstr, \tstr[Y]) \to \cat[D]( \tstr, \tstr[Y])
    \]
    is a bijection.
    Furthermore, the forgetful functor commutes with cylinders. It follows that as long as one restricts to bifibrant objects (in $\StratDN, \StratCN$) that are weakly frontier stratified then the resulting homotopy theory agrees with the classical homotopy theory of stratified spaces as studied in \cite{quinn1988homotopically,HughesPathSpaces,miller2013}.
\end{remark}
    It turns out that the triangularly cofibrant objects that have path-connected strata and are weakly frontier stratified are precisely the cofibrant objects in the refined setting (see \cref{prop:ref_vs_frontier} below). \\
Note that for strongly frontier stratified spaces the generating relations of \cref{con:top_strat} already define a partial order.
This leads to the following alternative characterization of refinedness, for specific stratified spaces.
\begin{proposition}\label{prop:ref_vs_frontier}
Let $\str \in \StratN$, and consider the following conditions:
     \begin{enumerate}
        \item \label{enum:char_via_frontier_cond_2} $\str$ is refined.
        \item \label{enum:char_via_frontier_cond_3} $\tstr$ has path-connected strata and is frontier stratified.
        \item \label{enum:char_via_frontier_cond_1} $\str$ has path-connected strata and is weakly frontier stratified.
    \end{enumerate}
     The first and the second property imply the third. Furthermore, if $\tstr$ is triangularly cofibrant, then the third property implies the first. Finally, if $\tstr$ is categorically fibrant, then the first property implies the third. In particular, for bifibrant stratified spaces in $\StratCN$ all three properties are equivalent.
    % Let $\str \in \StratCN$ be triangularly cofibrant stratified space. Then the following two statements are equivalent:
    % % \labelcref{enum:char_via_frontier_cond_1} below implies Condition \labelcref{enum:char_via_frontier_cond_2}.
    % \begin{enumerate}
    %     \item \label{enum:char_via_frontier_cond_1} $\str$ has path connected strata and is weakly frontier stratified.
    %     \item \label{enum:char_via_frontier_cond_2} $\str$ is refined.
    % \end{enumerate}
    % Conversely, if $\str$ is categorically fibrant, then these two statements are equivalent to 
    % \begin{enumerate}
    %     \setcounter{enumi}{2}
    % \end{enumerate}
    % In particular, a stratified space $\tstr \in \StratCRN$ is bifibrant, if and only if it is categorically fibrant, triangularly cofibrant and frontier stratified.
\end{proposition}
\begin{proof}
    The first two implications are obvious, using the fact that for any continuous map $f \colon X \to Y$, $x \in \overline{S}$ implies $f(x) \in \overline{f(S)}$.
   Suppose the third condition holds and $\str$ is triangularly cofibrant. To see that $\str$ is refined, we need to expose for any $x,y \in \str$ with $p := \sstr(x) < \sstr(y):= q$ a sequence of stratified paths $\gamma_i \colon \sReal{\stratSim[1]} \to \str$ from $x$ to $y$, as in \cref{def:ref_strat_sp}.
   Since $\tstr$ is weakly frontier stratified, it suffices to show that for any pair $p,q \in \ptstr$ and $x \in \tstr_p, y \in \tstr_q$ with
   \[
   \tstr_p \cap \overline{\tstr_q} \neq \emptyset,
   \]
   there is a stratified path from $x$ to $y$. Furthermore, since strata are assumed to be path-connected, it suffices to construct a path from any element in $\tstr_p$ to any element in $\tstr_q$. 
   Now, let $\tstr \xrightarrow{i}  \tstr[Y] \xrightarrow{r} \tstr$ expose $\tstr$ as a retract of a cellularly stratified space $\tstr[Y]$ (\cref{prop:char_cofibrants_wo_ref}).  
   Any nonempty closure intersection
   \[
   \tstr_{p} \cap \overline{\tstr_q} \neq \emptyset
   \]
   in $\tstr$ implies a non-empty closure intersection 
   \[
    \tstr[Y]_{i(p)} \cap \overline{\tstr[Y]_{i(q)}} \neq \emptyset.
   \]
    Conversely, every sequence of stratified paths in $\tstr[Y]$, starting and ending in $\tstr$, descends to a sequence of stratified paths in $\tstr$ with the same starting and end points. Hence, it suffices to show that a closure relation
   \[
   \tstr[Y]_{i(p)} \cap \overline{\tstr[Y]_{i(q)}} \neq \emptyset,
   \]
   implies the existence of a concatenable sequence of stratified paths starting in $\tstr[Y]_{i(p)}$ and ending in $\tstr[Y]_{i(q)}$.
   Choose a cell structure for $\str[Y]$ (see \cite[Def.\ref{hol:def:cell_struct}]{HoLinksWa}), with open cells $e_{i}$, $i \in I$. 
   Note that each open cell is entirely contained in only one stratum. We write $\sstr[Y](e_i) \in \pstr[Y]$ to denote the latter.
   Furthermore, note that if $e_i$ intersects the closure of $e_j$, then $\sstr[Y](e_i) \leq \sstr[Y](e_j)$.
   Suppose that $x\in \tstr[Y]_{i(p)} \cap \overline{\tstr[Y]_{i(q)}}$. In particular, $x$ is contained in the minimal subcomplex of $\tstr[Y]$ that contains $\str[Y]_{i(q)}$. This implies that there is a sequence of cells $e_0, e_1,e_2,\dots,e_n$, such that $x \in e_0$ and $e_n \subset \str[Y]_q$, and for each $i \in [n-1]$, $e_i$ intersects the closure of $e_{i+1}$. For $i \in [n-1]$, let $x_i \in e_{i} \cap \overline{e_{i+1}}$. Since $\overline{e_{i+1}}$ is the quotient of a stratified simplex $\sReal{\Delta^{\J}}$ over $\pstr[Y]$, there is a stratified path $\gamma'_i \colon \sReal{\stratSim[1]} \to \str[Y]$, starting in $x_i$, immediately entering and staying in $e_{i+1}$, and ending in $x_{i+1}$.\\ The converse implication is immediate.
    Finally, assume that $\tstr$ is categorically fibrant and refined (and is thus weakly frontier stratified). Then, whenever $\sstr(x) \leq \sstr(y)$,  any concatenable sequence of stratified paths from $x$ to $y$ induces a stratified path, $\gamma \colon \sReal{\stratSim[1]} \to \tstr$ from $x$ to $y$, which implies \[
    x = \gamma(0) \subset \overline{\gamma(0,1]} \subset \overline{\tstr_{\sstr(y)}}.\] In particular, it follows that 
    \[
       \tstr_{p} \cap \overline{\tstr_q} \neq \emptyset \implies \tstr_{p} \subset \overline{\tstr_q}.
    \]
\end{proof}  
The characterization of refinedness in \cref{prop:ref_vs_frontier} allows us to represent all homotopy types in $\StratCRN$ through the following particularly convenient stratified spaces.
\begin{definition}\label{def:perfect_strat_sp}
    A stratified space $\str \in \StratN$ is called \define{CFF stratified} (C for cellular, first F for frontier, second F for fibrant), if it fulfills the following conditions:
    \begin{enumerate}
        \item $\str$ is cellularly stratified.
        \item $\str$ has nonempty, connected strata and is frontier stratified.
        \item $\str$ has the horn filling property with respect realizations of inner stratified horn inclusions $\sReal{ \stratHorn \hookrightarrow \stratSim}$, $0<k<n$.
    \end{enumerate}
    We denote the full simplicial subcategory of $\Strat$ given by CFF stratified spaces by $\CFF$.
\end{definition}
CFF stratified spaces can be seen as an analogue to CW complexes in the classical scenario. However, note that no assumptions are made that the attaching maps of cells only map to cells of certain dimensions. This is necessary for the small object argument to be able to produce CFF stratified spaces. 
\begin{example}\label{ex:of_classical_sp}
    Every stratified space that admits a PL structure that is compatible with its stratification is cellularly stratified. Furthermore, in increasing order of generality, every Whitney stratified space, Thom-Mather stratified space, pseudomanifold, or conically stratified space (in the sense of \cite[A.5]{HigherAlgebra}) has the inner horn filling property, by \cite[Thm. A.6.5]{HigherAlgebra}. Hence, it follows that if we equip such spaces with the refined stratification (induced by the frontier condition, and by taking path components of strata) and additionally assume a piecewise linear structure, then they provide examples of CFF stratified spaces.
    Note also that for Thom-Mather (and hence for Whitney stratified spaces) a piecewise linear structure (compatible with the stratification) always exists (\cite{TriangulationsGoresky}). 
\end{example}
Let us now characterize the bifibrant objects of $\sStratCR$ in terms of CFF stratified spaces. Using the small object argument on the generating classes in \cref{thm:ex_mod_struct_red}, together with \cref{prop:rel_cof_and_ref,prop:ref_vs_frontier}, we obtain the following corollary.
\begin{corollary}\label{cor:perf_strat_sp_are_essentially_surj}
    A stratified space $\str \in \StratCR$ is bifibrant if and only if it is a retract of a CFF stratified space.
    Furthermore, every bifibrant stratified space is stratified homotopy equivalent to a CFF stratified space.
    Every stratified space is categorically equivalent to a CFF stratified space.
\end{corollary}
In particular, it follows that the homotopy theory defined by $\StratCRN$ may equivalently be interpreted in terms of CFF stratified spaces:
\begin{corollary}\label{cor:CFF_equ_to_hotheory}
        Denote by $H_s$ the class of stratified homotopy equivalences between CFF stratified spaces.
        The inclusion $\CFFN \to \StratN$ induces an equivalence of $\infty$-categories
            \[
            \CFFN[H_s^{-1}] \xrightarrow{\simeq} \AltStratCR.
            \]
\end{corollary}
Combining this result with \cref{thm:equ_bifib_and_layered_qc}, we obtain the following version of the stratified homotopy hypothesis:
\begin{corollary}
    Luries exit-path construction induces an equivalence of $\infty$-categories
            \[ \CFFN[H_s^{-1}] \xrightarrow{\simeq} \iCatO\]
    between CFF stratified spaces localized at stratified homotopy equivalences and layered $\infty$-categories.
\end{corollary}
\subsection{Stratified homotopy link fibrations}\label{subsec:stratified_homotopy_link_fibrations}
Much of the literature on stratified spaces takes an inductive approach to the study of stratified spaces.
This follows the observation that in many geometric scenarios a stratified space $\str$ over a finite linear poset $\pos$ with minimal element $p$ can be decomposed into a diagram 
\[
\str_{p} \leftarrow \str[E] \rightarrow \str[Y]
\]
where $\str[E]$ and $\str[Y]$ are stratified over $\pos_{>p}$ and $\str[E] \to \str_p$ is a kind of \textit{fibration with stratified fiber} or a retraction associated to some stratified notion of a block bundle. See \cite{weinberger1994topological} for a good overview of such phenomena and  \cite{thom1969ensembles,Stone1972,weinberger1994topological,HughesPathSpaces}, for some examples of this approach.
If one keeps inductively repeating this kind of procedure with $\str[E]$ and $\str[Y]$, one ultimately ends up with a diagram of stratified spaces indexed over the poset of linear flags in $\pos$, $\sd(\pos)$. Whatever geometric construction one uses to decompose one's stratified space, at least from a homotopy-theoretic perspective, one ultimately ends up with the associated diagram of generalized homotopy links $\HolIP[](\str) \in \FunC ((\sd(\pos))^\op, \iSpaces)$
\footnote{This is more of a meta theorem, which at least holds in all of the cases known to the author. The reader should think of it as a heuristic to motivate the rigorous mathematics performed below.}.
In order to capture the homotopy-theoretic essence of these types of constructions, Hughes (see \cite{HughesPathSpaces}) used a stratified version of the pairwise homotopy links of \cite{quinn1988homotopically}, thus obtaining a functorial, homotopy-theoretic version of such decompositions. In this subsection, we replicate this construction using the cartesian structure on $\StratN$, and derive a series of homotopy-theoretic consequences. Most of these are probably known to the expert, at least in the alternative framework of homotopically stratified spaces. However, we think they may also help to connect the classical inductive approach to stratified algebraic topology with our results on stratified homotopy theory.
\begin{notation}
    In the following, we will often cover both the categorical as well as the diagrammatic cases in one statement. When we add the prefix diagrammatic or categorical to names of stratum-preserving maps, such as fibrations, we mean that they are fibrations in the respective model structure on $\TopPN$. We will then often add the alternative prefix in parenthesis to indicate that both cases hold.
\end{notation}
\begin{construction}\label{con:decomp_into_strat_hol}
    Let $\pos$ be some poset, and let $\str \in \TopPN$. Consider the evaluation map 
    \[
    \str^{\sReal{\Delta^{[1]}}} \to \str^{ \sReal{\partial \Delta^{[1]}}} = \str \times \str,\] where the left component of this map is given by evaluation of a stratified path at $0$ and the right component by evaluation at $1$. We denote by $\HolS (\str)$ the $\pos_{>p}$ stratified space, obtained via the following diagram of pullback squares in $\StratN$.
    % https://q.uiver.app/#q=WzAsNixbMSwwLCJcXHN0cl57XFxzUmVhbHtcXERlbHRhXntbMV19fX0iXSxbMSwxLCJcXHN0ciBcXHRpbWVzIFxcc3RyIl0sWzAsMSwiXFxzdHJfcCBcXHRpbWVzIFxcc3RyX3s+cH0iXSxbMCwwLCJcXEhvbFMgKFxcc3RyKSJdLFsxLDIsIlxccG9zIFxcdGltZXMgXFxwb3MiXSxbMCwyLCJcXHBvc197PnB9IFxcY29uZyBcXHtwXFx9IFxcdGltZXMgXFxwb3Nfez5wfSAiXSxbMiwxLCIiLDAseyJzdHlsZSI6eyJ0YWlsIjp7Im5hbWUiOiJob29rIiwic2lkZSI6InRvcCJ9fX1dLFszLDJdLFswLDFdLFszLDBdLFszLDEsIiIsMSx7InN0eWxlIjp7Im5hbWUiOiJjb3JuZXIifX1dLFsxLDRdLFsyLDVdLFs1LDQsIiIsMSx7InN0eWxlIjp7InRhaWwiOnsibmFtZSI6Imhvb2siLCJzaWRlIjoidG9wIn19fV0sWzIsNCwiIiwxLHsic3R5bGUiOnsibmFtZSI6ImNvcm5lciJ9fV1d
\begin{diagram}
	{\HolS (\str)} & {\str^{\sReal{\Delta^{[1]}}}} \\
	{\str_p \times \str_{>p}} & {\str \times \str} \\
	{\pos_{>p} \cong \{p\} \times \pos_{>p} } & {\pos \times \pos}
	\arrow[from=1-1, to=1-2]
	\arrow[from=1-1, to=2-1]
	\arrow["\lrcorner"{anchor=center, pos=0.125}, draw=none, from=1-1, to=2-2]
	\arrow[from=1-2, to=2-2]
	\arrow[hook, from=2-1, to=2-2]
	\arrow[from=2-1, to=3-1]
	\arrow["\lrcorner"{anchor=center, pos=0.125}, draw=none, from=2-1, to=3-2]
	\arrow[from=2-2, to=3-2]
	\arrow[hook, from=3-1, to=3-2]
\end{diagram}
    In other words, we may either think of $\HolS (\str)$ as a restriction of $\sReal{\str}^{\sReal{\Delta^1}}$ to $\{p\} \times \str_{>p}$, or as the stratified space of such stratified paths in $\str$ that start in the $p$-stratum, and immediately exit.
    We call this stratified space the \define{$p$-th stratified homotopy link of $\str$.}
    We may then compose ${\HolS (\str)} \to \str_p \times \str_{>p}$ with the projections to $\str_{p}$ and $\str_{>p}$, to obtain a diagram in $\StratN$
   \[
    \str_p \twoheadleftarrow {\HolS (\str)} \rightarrow \str_{>p}
    \]
    where the two maps are respectively given by evaluating a stratified path in $\str$ at the start and end-point.
    This construction extends to a functor from $\StratN_{\pos}$ into the category of spans in $\StratN$
    \[
   B \leftarrow \str[E] \rightarrow \str[Y]
     \]
      where $T$ is trivially stratified, and $\str[E] \to \str[Y]$ is a stratum-preserving map over $\pos_{>p}$.
     % which can be identified with the fiber product of categories 
%     % https://q.uiver.app/#q=WzAsNSxbMSwwLCJcXEZ1bkMoXFxEZWx0YV4xLCBcXFRvcFBOKSJdLFswLDIsIlxcRnVuQyhcXERlbHRhXjEsIFxcVG9wTikiXSxbMCwwLCJcXEZ1bkMoXFxEZWx0YV4xLCBcXFRvcE4pXFx0aW1lc197XFxUb3BQTn1cXEZ1bkMoXFxEZWx0YV4xLCBcXFRvcFBOKSJdLFsxLDEsIlxcVG9wUE4iXSxbMSwyLCJcXFRvcE4iXSxbMiwxXSxbMiwwXSxbMCwzLCJcXGV2XzAiLDFdLFszLDRdLFsxLDQsIlxcZXZfMCIsMV0sWzIsNCwiIiwxLHsic3R5bGUiOnsibmFtZSI6ImNvcm5lciJ9fV1d
% \begin{diagram}
% 	{\FunC(\Delta^1, \TopN)\times_{\TopN}\FunC(\Delta^1, \TopPN)} & {\FunC(\Delta^1, \TopPN)} \\
% 	& \TopPN \\
% 	{\FunC(\Delta^1, \TopN)} & \TopN
% 	\arrow[from=1-1, to=1-2]
% 	\arrow[from=1-1, to=3-1]
% 	\arrow["\lrcorner"{anchor=center, pos=0.125}, draw=none, from=1-1, to=3-2]
% 	\arrow["{\ev_0}"{description}, from=1-2, to=2-2]
% 	\arrow[from=2-2, to=3-2]
% 	\arrow["{\ev_0}"{description}, from=3-1, to=3-2]
% \end{diagram}
\end{construction}
As a corollary of the cartesianity of the semi-model structures on $\StratN$, one obtains:
\begin{lemma}
    Let $\str \in \TopPN$ be diagrammatically fibrant (categorically fibrant). Then the starting point evaluation map
    \[
    \HolS (\str) \to \str_p
    \]
    is a diagrammatic (categorical) fibration.
\end{lemma}
\begin{proof}
    It follows from the cartesianity of the model structures on $\StratN$, that the right vertical in the pullback square 
    % https://q.uiver.app/#q=WzAsNCxbMSwwLCJcXHN0cl57XFxzUmVhbHtcXERlbHRhXntbMV19fX0iXSxbMSwxLCJcXHN0ciBcXHRpbWVzIFxcc3RyIl0sWzAsMSwiXFxzdHJfcCBcXHRpbWVzIFxcc3RyX3s+cH0iXSxbMCwwLCJcXEhvbFMgKFxcc3RyKSJdLFsyLDEsIiIsMCx7InN0eWxlIjp7InRhaWwiOnsibmFtZSI6Imhvb2siLCJzaWRlIjoidG9wIn19fV0sWzMsMl0sWzAsMV0sWzMsMF0sWzMsMSwiIiwxLHsic3R5bGUiOnsibmFtZSI6ImNvcm5lciJ9fV1d
\begin{diagram}
	{\HolS (\str)} & {\str^{\sReal{\Delta^{[1]}}}} \\
	{\str_p \times \str_{>p}} & {\str \times \str}
	\arrow[from=1-1, to=1-2]
	\arrow[from=1-1, to=2-1]
	\arrow["\lrcorner"{anchor=center, pos=0.125}, draw=none, from=1-1, to=2-2]
	\arrow[from=1-2, to=2-2]
	\arrow[hook, from=2-1, to=2-2]
\end{diagram}
    is a diagrammatic (categorical) fibration. Consequently, so is the left vertical. As $\str$ is diagrammatically (categorically) fibrant, the projection map $\str_{p} \times \str_{>p} \to \str_{p}$ is also a diagrammatic (categorical) fibration. The evaluation map in question is now the composition of these two diagrammatic (categorical) fibrations. 
\end{proof}
Let us quickly compare our framework to the stratified path spaces studied in \cite{HughesPathSpaces}. This first requires another remark on spaces with decompositions:
\begin{remark}    \cref{rem:different_notions_of_stratified_homs} also holds in a self-enriched sense, providing an answer to \cref{iQ:ExponentialObjects}: In \cite{HughesPathSpaces} the author studied stratified notions of mapping space, given by equipping the set of decomposition maps $\cat{D}(\tstr, \tstr[Y])$ with the subspace topology of the compact-open topology, and the decomposition over $\pstr[Y]^{\pstr}$ induced by mapping a stratified space to its underlying map of posets.
    If we work in the setting where $\TopN$ is the category of $\Delta$-generated spaces or compactly generated spaces, and thus replace the subspace topology with its respective Kelleyfication, then this construction defines the internal mapping space in the category of decomposition spaces.
    There is an obvious map of decomposition spaces
    \[
    \mathcal{D}( \tstr[Y]^{\tstr} ) \to \mathcal{D}(\tstr[Y])^{\mathcal{D}(\tstr)},
    \]
    given on the set level by the map $\StratN( \tstr, \tstr[Y]) \to \cat[D]( \tstr, \tstr[Y])$. Therefore, whenever $\tstr[X]$ is refined, the comparison map of stratified mapping spaces above is bijective. Furthermore, by construction of $\tstr[Y]^{\tstr}$ and \cite[Cor. 2.2.11]{may2016finite}, it is a decomposition preserving homeomorphism, whenever $\pstr$ and $ \pstr[Y]$ are finite.
    In other cases, it is the map that refines the topology on $\mathcal{D}(\tstr[Y])^{\mathcal{D}(\tstr)}$ such that the induced map \begin{equation*}\label{map:comp_map_mapping_spaces}
        \mathcal{D}(\tstr[Y])^{\mathcal{D}(\tstr)} \to \pstr[Y]^{\pstr}
    \end{equation*}
    is continuous.
\end{remark}
    \begin{remark}
    In \cite{HughesPathSpaces}, Hughe's main object of study was the space of stratified paths in a stratified space $\tstr$ with finitely many strata, which start in a closed union of strata $\tstr[A] \subset \str$, which we denote by $\textnormal{Path}_{nsp}(\str, \str[A])$. In the special case when $\str[A] = \str_p$, for some $p \in \pos$, then we can think of $\textnormal{Path}_{nsp}(\str, \str[A])$ as the union of $\HolS(\str)$ with the space of paths entirely contained in $\str[X]_p$.
    One of the main results of \cite{HughesPathSpaces} is that when $\tstr$ is a homotopically stratified space, then the starting point evaluation map $\textnormal{Path}_{nsp}(\str, \str[A]) \to \str[A]$ lifts stratified homotopies. 
    Using the cartesian closedness of the semi-model structures on $\StratN$ we may recover a version of this result in $\StratN$:
    Given a categorically (respectively, diagrammatically) fibrant stratified space $\tstr \in \StratN$, the starting point evaluation map $\tstr^{\sReal{\Delta^{[1]}}} \to \tstr$ is a stratified fibration. Now, for any subspace $\tstr[A] \subset \tstr$, equipped with the induced stratification,  we may consider the pullback diagram
    \begin{diagram}
    \tstr^{\sReal{\Delta^{[1]}}} \times_{\ev_0}\tstr[A] \arrow[r] \arrow[d] &\tstr[X]^{\sReal{\Delta^{[1]}}} \arrow[d, "\ev_0"]\\
    \tstr[A] \arrow[r] & \tstr[X] \spaceperiod
    \end{diagram}
    It follows that the right vertical is a categorical (or, respectively, diagrammatic) fibration.
    In particular, it has the lifting property with respect to stratified homotopies with cofibrant source. Since $\sReal{\Delta^{[1]}}$ is clearly strongly frontier stratified, the natural comparison map in (\ref{map:comp_map_mapping_spaces}) induces a natural continuous bijection \[
    \mathcal{D}(\tstr^{\sReal{\Delta^{[1]}}} \times_{\ev_0}\tstr[A] ) \to \textnormal{Path}_{nsp}(\str, \str[A]),\] which refines the topology on $\textnormal{Path}_{nsp}(\str, \str[A])$ in order to turn it into a poset-stratified space (and make it compactly or $\Delta$-generated). Hence, one obtains a Serre (as opposed to Hurewicz-style homotopy theory) version of Hughes's result, replacing homotopically stratified spaces with the weaker condition of diagrammatic fibrancy and not requiring that $\tstr[A]$ is a closed union of strata, but only obtaining the homotopy lifting property with respect to triangularly cofibrant stratified spaces.
\end{remark}

Much like generalized homotopy links, stratified homotopy links can be computed in terms of appropriately stratified regular neighborhoods. Let us first observe that stratified homotopy links can be computed locally:
\begin{lemma}
    Let $\str \in \TopPN$, $p \in \pos$ and $\str[N] \subset \str$ be a neighborhood of the $p$-stratum $\str_p \subset \str$. Then the induced map 
    \[
    \HolS \str[N] \to \HolS \str 
    \]
    is a diagrammatic equivalence.
\end{lemma}
\begin{proof}
    This is a consequence of \cite[Prop. \ref{hol:prop:computing_links_using_nbhds}]{HoLinksWa} and \cref{lem:identifying_iterated_holink_with_gen_holink} below.
\end{proof}
Furthermore, one obtains the following stratified analogue of the two strata case of \cite[Prop. \ref{hol:prop:computing_holinks_via_aspire}]{HoLinksWa}.
\begin{lemma}\label{lem:strat_aspir}
    Let $\str[N]$ be a stratified space, and suppose we are given a stratified map
    \[
    R \colon \str[N]_{\geq p} \times \sReal{\Delta^{[1]}} \to \str[N]_{\geq p}
    \]
    such that 
    \begin{enumerate}
        \item $R(x,t) = x$, for $x \in \str[N]_p$ or $t = 1$;
        \item $R(x,0) \in \str[N]_p$, for all $x \in \str$.
    \end{enumerate}
    Then the evaluation map 
    \[
    \ev_1 \colon \HolS (\str[N]) \to \str[N]_{>p}
    \]
    is a stratum-preserving homotopy equivalence.
\end{lemma}
\begin{proof}
    A section $s \colon \str[N]_{>p} \to  \HolS (\str[N]) $ of this map is provided by
    \[
    x \mapsto \{ t \mapsto R(x,t) \}.
    \]
    A homotopy between $s \circ \ev_1$ and $1_{\HolS (\str[N])}$ is constructed exactly as in the two strata case of \cite[Prop. \ref{hol:prop:properties_of_star_homotopy}]{HoLinksWa}. One easily verifies that the homotopy given there is stratum-preserving.
\end{proof}
Combining these two results, one obtains:
\begin{corollary}
    Let $\str \in \TopPN$, $p \in \pos$ and $\str[N] \subset \str$ be a neighborhood of the $p$-stratum $\str_p \subset \str$. Suppose that it admits a stratified map $R \colon \str[N]_{\geq p} \times \sReal{\Delta^{[1]}} \to \str[N]_{\geq p}$ as in \cref{lem:strat_aspir}. Then there is a zig-zag of diagrammatic equivalences \[
    \str[N]_{>p} \xleftarrow{\ev_1} \HolS \str[N] \xrightarrow{\simeq } \HolS \str.
    \]
\end{corollary}
In the introduction, we already alluded to the fact that one can think of generalized homotopy links as arising as the iterated stratified homotopy links of a stratified space. This follows from the following observation:
\begin{lemma}\label{lem:identifying_iterated_holink_with_gen_holink}
     Let $\I = [p_1 < \dots < p_n] \in \sd(\pos_{>p})$ be a regular flag in $\pos$ of length $n-1$, containing only elements larger than $p \in \pos$ and let $\str \in \TopPN$.
    There is a natural weak homotopy equivalence of simplicial sets
    \[
    \HolIP(\HolS (\str)) \simeq \HolIP[{\{p\} \cup \I }] (\str).
    \]
\end{lemma}
\begin{proof}
      By definition of $\HolS (\str)$ as a restriction of $\str^{\sReal{\Delta^{[1]}}}$ to certain strata, we can identify the $\I$-th homotopy link $\HolIP[\I] (\HolS (\str))$ with the $\I'$-th homotopy link of $\str^{\sReal{\Delta^{[1]}}}$, where $\I'$, is the flag 
    \[
    [ (p,p_1) < \dots <( p,p_n)]
    \]
    in $\pos \times \pos$.
    Observe that under the natural isomorphism 
    \[
    \Strat( \sReal{\Delta^{[n-1]}}, \str^{\sReal{\Delta^{[1]}}}) \cong \Strat (\sReal{ \Delta^{{[n-1]}}} \times \sReal{\Delta^{[1]}} , \str)  \cong \Strat (\sReal{ \Delta^{[n-1]} \times \Delta^{[1]}} , \str) 
    \]
    the component of $ \Strat( \sReal{\Delta^{{[n-1]}}}, \str^{\sReal{\Delta^{[1]}}})$ which is given by $\HolIP[{\I'}] (\str^{\sReal{\Delta^{[1]}}})$ is identified with the  component of $\Strat (\sReal{ \Delta^{[n-1]} \times \Delta^{[1]}} , \str)$ given by such stratified maps whose underlying map of posets is given by \begin{align*}
        (k,i) \mapsto \begin{cases}
            p_k & \textnormal{, if $i = 1$} \\
            p & \textnormal{otherwise.}
        \end{cases}
    \end{align*}
    We may equivalently identify this component with the simplicial mapping space
    \[
    \TopP ( \str[Y], \str)
    \]
    where $\str[Y]$ is the stratified space obtained by equipping 
    $\real{\Delta^{n-1} \times \Delta^1}$ with the stratification 
    \[
    (x,t) \mapsto \begin{cases}
        s_{\sReal{\Delta^\I}} (x) & \textnormal{, if $t >0$} \\
        p \textnormal{, otherwise.} 
    \end{cases}
    \]
    In this way, we have obtained an identification
    \[
    \HolIP(   \HolS (\str))  \cong \TopP ( \str[Y], \str)
    \]
    Now, consider $\sReal{\Delta^\I}$ as embedded into $\sReal{\Delta^{\{p\} \cup \I }}$ as its $\I$-face.
    The collapsing map
    \begin{align*}
        \tstr[Y] &\to \sReal{\Delta^{\{p\} \cup \I }} \\
        (x,t) &\mapsto ((1-t)e_0 + tx)
    \end{align*}
    fits into a pushout square
    % https://q.uiver.app/#q=WzAsNCxbMCwxLCJcXHN0cltEXSJdLFswLDAsIlxcc1JlYWx7XFxEZWx0YV5ufSBcXHRpbWVzIFxce3AgXFx9Il0sWzEsMCwiXFxzdGFyIl0sWzEsMSwiXFxzUmVhbHtcXERlbHRhXntwIFxcY3VwIFxcSX19Il0sWzEsMCwiIiwwLHsic3R5bGUiOnsidGFpbCI6eyJuYW1lIjoiaG9vayIsInNpZGUiOiJ0b3AifX19XSxbMSwyXSxbMCwzXSxbMiwzXSxbMywxLCIiLDEseyJzdHlsZSI6eyJuYW1lIjoiY29ybmVyIn19XV0=
\begin{diagram}
	{\real{\Delta^{n-1}} \times \{p \}} & \{p\} \\
	{\str[Y]} & {\sReal{\Delta^{ \{p\} \cup \I}}}
	\arrow[from=1-1, to=1-2]
	\arrow[hook, from=1-1, to=2-1]
	\arrow[from=1-2, to=2-2]
	\arrow[from=2-1, to=2-2]
	\arrow["\lrcorner"{anchor=center, pos=0.125, rotate=180}, draw=none, from=2-2, to=1-1]
\end{diagram}
    of cofibrant objects in $\TopPDN$, with verticals given by cofibrations. In particular, this square is a homotopy pushout. As the upper horizontal is a weak diagrammatic equivalence, so is the lower horizontal. As both stratified spaces in the lower horizontal are diagrammatically fibrant and triangularly cofibrant, it follows that the lower horizontal is a stratified homotopy equivalence. Hence, precomposing with this map defines a homotopy equivalence of simplicial sets
    \[
    \TopP ( \str[Y], \str) \simeq \TopP ( \sReal{\Delta^{ \{p \} \cup \I }}, \str) = \HolIP[{\{p\} \cup \I}](\str).
    \]
    This completes the proof.
\end{proof}
\begin{remark}
    Applying \cref{lem:identifying_iterated_holink_with_gen_holink} inductively, one obtains that the generalized homotopy link $\HolIP (\str)$ (modeled by a topological mapping space instead of a simplicial set, as in \cite{douteauwaas2021}) associated to a stratified space $\str \in \TopPN$ and a flag $\I = [p_0 < \dots <p_n]$ can be homotopically identified with the iterated stratified homotopy link
    \[
    (\HolS[p_n] \circ \dots \circ \HolS[p_0]) (\str). 
    \]
\end{remark}
As an immediate corollary of \cref{lem:identifying_iterated_holink_with_gen_holink}, one obtains:
\begin{corollary}\label{cor:detecting_strat_equ_on_strat_links}
    Let $w \colon \str \to \str[Y] \in \TopPN$ be a stratum-preserving map. $w$ is a diagrammatic equivalence if and only if, for each $p \in \pos$, the following holds:
    \begin{enumerate}
        \item $w_p \colon \str \to \str[Y]$ is a weak homotopy equivalence;
        \item The induced map of stratified homotopy links $\HolS (\str) \to \HolS(\str[Y])$ is a diagrammatic equivalence.
    \end{enumerate}
\end{corollary}
The true power of this criterion lies in the fact that it allows one to think of stratum-preserving maps between fibrant objects as having a \textit{tangential component} along the strata, as well as a \textit{normal component}, given by the fibers of stratified homotopy links. To explain this more rigorously, we need the following lemma:
\begin{lemma}\label{lem:fiberwise_equ_of_strat_fib}
    Suppose we are given a commutative square in $\StratN$
    % https://q.uiver.app/#q=WzAsMyxbMCwwLCJcXHN0cltFXSciXSxbMiwwLCJcXHN0cltFXSJdLFsxLDEsIkIiXSxbMCwxLCJcXHRpbGRlIHciXSxbMCwyXSxbMSwyXV0=
\begin{diagram}
	{\str[X]} && {\str[Y]} \\
	& B
	\arrow["{\tilde w}", from=1-1, to=1-3]
	\arrow[from=1-1, to=2-2]
	\arrow[from=1-3, to=2-2]
\end{diagram}
    with $B$ trivially stratified and $\str \to \str[Y]$ a stratum-preserving map over $\pstr[X] = \pstr[{Y}]$. Suppose that both diagonals are diagrammatic (categorical) fibrations. 
    \begin{enumerate}
        \item For every $x \in B$, the induced map on fibers $\str[X]_x \to \str[Y]_{x}$ is a diagrammatic (categorical) equivalence;
        \item For every path component of $B$, there exists a representative $x \in B$ such that the induced map on fibers $\str[X]_x \to \str[Y]_{x}$ is a diagrammatic (categorical) equivalence;
        \item $\tilde{w}$ is a diagrammatic (categorical) equivalence.
    \end{enumerate}
\end{lemma}
\begin{proof}
    Observe that we only need to prove the diagrammatic case. The categorical case then follows by the Whitehead theorem for left Bousfield localizations.
    That the first condition implies the second is trivial. 
    That the last condition implies the first follows by observing that the
    two squares
    % https://q.uiver.app/#q=WzAsOCxbMCwwLCJcXHN0cltYXV94Il0sWzEsMCwiXFxzdHJbWF0iXSxbMywwLCJcXHN0cltZXV94Il0sWzQsMCwiXFxzdHJbWV0iXSxbMywxLCJcXHN0YXIiXSxbNCwxLCJCIl0sWzEsMSwiQiJdLFswLDEsIlxcc3RhciJdLFs0LDUsIngiLDAseyJzdHlsZSI6eyJ0YWlsIjp7Im5hbWUiOiJob29rIiwic2lkZSI6InRvcCJ9fX1dLFs3LDYsIngiLDAseyJzdHlsZSI6eyJ0YWlsIjp7Im5hbWUiOiJob29rIiwic2lkZSI6InRvcCJ9fX1dLFsxLDZdLFswLDddLFswLDFdLFsyLDRdLFsyLDNdLFszLDVdXQ==
\begin{diagram}
	{\str[X]_x} & {\str[X]} && {\str[Y]_x} & {\str[Y]} \\
	\star & B && \star & B
	\arrow[from=1-1, to=1-2]
	\arrow[from=1-1, to=2-1]
	\arrow[from=1-2, to=2-2]
	\arrow[from=1-4, to=1-5]
	\arrow[from=1-4, to=2-4]
	\arrow[from=1-5, to=2-5]
	\arrow["x", hook, from=2-1, to=2-2]
	\arrow["x", hook, from=2-4, to=2-5]
\end{diagram}
    are not just pullback, but (by the assumption on fibrancy of the right verticals and since every trivially stratified space is fibrant) even homotopy pullback, in the sense that the underlying square in the associated $\infty$-category $\AltStratDR$ is pullback. 
    (There are many ways to see this. One of them is applying the equivalence of $\infty$-categories induced by $\SingS$ together with the dual of \cite[Thm. 4.2.4.1.]{HigherTopos}.)
    Thus it follows that the induced 
    map on fibers $\str[X]_x \to \str[Y]_{x}$ is an isomorphism in $\AltStratDR$ and hence a diagrammatic equivalence. To see the final remaining implication, observe that since diagrammatic equivalences can be verified on connected components, we may without loss of generality assume that $B$ is connected. Now, let $n \in \mathbb{N}$. Observe that $\AltHolIP( \star) = \Delta^0$. Applying $\AltHolIP$ to the pullback squares above, and using the commutativity of $\AltHolIP[n]$ with pullbacks, we thus obtain a morphism of fiber sequences
    % https://q.uiver.app/#q=WzAsNixbMCwyLCJcXFN0cmF0KCBcXHN0cltaXSxCKSJdLFsxLDIsIlxcU3RyYXQoIFxcc3RyW1pdLEIpIl0sWzAsMSwiXFxTdHJhdCggXFxzdHJbWl0sXFxzdHJbWF0pIl0sWzEsMSwiXFxTdHJhdCggXFxzdHJbWl0sXFxzdHJbWV0pIl0sWzAsMCwiXFxTdHJhdCggXFxzdHJbWl0sXFxzdHJbWF1feCkiXSxbMSwwLCJcXFN0cmF0KCBcXHN0cltaXSxcXHN0cltZXV94KSJdLFswLDEsIiIsMCx7ImxldmVsIjoyLCJzdHlsZSI6eyJoZWFkIjp7Im5hbWUiOiJub25lIn19fV0sWzIsMF0sWzMsMV0sWzIsM10sWzQsMl0sWzUsM10sWzQsNV1d
\begin{diagram}
	{\AltHolIP[n](\str[X]_x)} & {\AltHolIP[n](\str[Y]_x)} \\
	{\AltHolIP[n](\str[X])} & {\AltHolIP[n](\str[Y])} \\
	{\AltHolIP[n](B)} & {\AltHolIP[n](B)}.
	\arrow[from=1-1, to=1-2, "\simeq"]
	\arrow[from=1-1, to=2-1]
	\arrow[from=1-2, to=2-2]
	\arrow[from=2-1, to=2-2]
	\arrow[from=2-1, to=3-1]
	\arrow[from=2-2, to=3-2]
	\arrow[Rightarrow, no head, from=3-1, to=3-2]
\end{diagram}
Observe that by simpliciality of the model structure on $\StratN$, both lower verticals are Kan-fibrations. Hence, the left and the right vertical sequence are homotopy fiber sequences. Furthermore, as $B$ is trivially stratified, we have that $\AltHolIP[n](B) \simeq \Top (\real{\Delta^n}, B) \simeq \SingS (B)$ is path-connected. Hence, it follows by the classical fiberwise characterization of weak equivalences between Kan-fibrations that $\AltHolIP[n](\str) \to \AltHolIP(\str[Y])$ is a weak homotopy equivalence. 
\end{proof}
\begin{notation}
    Given a diagrammatically fibrant stratified space $\str$, we call the (homotopy) fiber of $\HolS \str \to \str_{p}$ at $x \in \str_p$ the \define{local homotopy link} of $\str$ at $x$.
\end{notation}
As a corollary, combining \cref{lem:fiberwise_equ_of_strat_fib,cor:detecting_strat_equ_on_strat_links}, as well as the stability of weak equivalences between fibrant objects under pullbacks along fibrations, one obtains the following detection criterion for diagrammatic equivalences.
\begin{corollary}\label{cor:char_local_strat_ho_equ}
    Let $w \colon \str \to\str[Y] \in \TopPN$ be a stratum-preserving map of diagrammatically fibrant stratified spaces. Then the following are equivalent:
    \begin{enumerate}
        \item $w$ is a diagrammatic equivalence;
        \item For each $p \in \pos$, the induced map $w_p \colon \str_p \to \str[Y]_p$ is a weak homotopy equivalence, and for each $x \in \str_p$ (or just for a representative system of path components) the induced map on local homotopy links
        \[
        \HolS (\str)_{x} \to \HolS (\str[Y])_{w(x)}
        \]
        is a diagrammatic equivalence.
    \end{enumerate}
    The analogous equivalence holds for categorically fibrant spaces and categorical equivalences. Even more, in this case it follows from the d\'ecollage condition that these two statements are furthermore equivalent to 
    \begin{enumerate}
     \setcounter{enumi}{2}
        \item For each $p \in \pos$, the induced map $w_p \colon \str_p \to \str[Y]_p$ is a weak homotopy equivalence, and for each $x \in \str_p$ (or just for a representative system of path components) the induced map on local homotopy links
        \[
        \HolS (\str)_{x} \to \HolS (\str[Y])_{w(x)}
        \]
        induces weak equivalences on strata.
    \end{enumerate}
\end{corollary}
Using an analogous (but significantly easier) argument to the proof of \cite[Prop. \ref{hol:prop:computing_links_using_nbhds}]{HoLinksWa}, one can show that local homotopy links are indeed \textit{local}, in the following sense.
\begin{lemma}\label{lem:local_holinks_are_local}
    Let $\str$ be a diagrammatically fibrant stratified space. Let $p \in \pos$, $x \in \str_{p}$ and let $\str[N]$ be a neighborhood of $x \in \str$. Then the inclusion of local homotopy links
    \[
    \HolS (\str[N])_{x} \to \HolS (\str)_{x}
    \]
    is a diagrammatic equivalence.
\end{lemma}
\begin{remark}
   Together with \cref{cor:char_local_strat_ho_equ}, \cref{lem:local_holinks_are_local} allows us to verify stratified homotopy equivalence between bifibrant stratified spaces in terms of a two-step procedure, involving a global computation, which is only concerned with strata, and a purely local computation, concerning the local homotopy links. In the case of categorically fibrant stratified spaces, one even only needs to verify strata-wise weak equivalence on the local level. 
\end{remark}
In many cases of geometric stratified spaces, such as topological pseudo manifolds, local stratified homotopy links admit explicit geometric models.
\begin{remark}
    Let $\str \in \TopPN$ be a conically stratified space (see \cite{HigherAlgebra}). Recall that this means that, for every $p \in \pos$ and every $x \in \str_{p}$, there exists a neighborhood $\str[N]$ of $x$, that is locally stratum-preserving homeomorphic to a product $U_x \times C\str[L]_x$, where $U_x$ is a trivially stratified space and $C \str[L]_x$ is the stratified cone on a stratified space $\str[L]_x \in \TopPN[{\pos_{>p}}]$ obtained by equipping the (teardrop) cone on $L_x$ with the stratification 
    \[
    [y,t] \mapsto \begin{cases}
        p & \textnormal{, if $t=0$} \\
        s_{\str[L]_x}(y) & \textnormal{, otherwise.}
    \end{cases}
    \]
    Suppose, furthermore, that the strata of $\str$
    are locally (weakly) contractible, such that we may choose $U_x$ to be (weakly) contractible.
    Then, combining \cref{lem:local_holinks_are_local,lem:strat_aspir} we obtain a sequence of diagrammatic equivalences over $\pos$:
    \begin{align*}
         \HolS(\str)_x &\simeq \HolS(\str[N])_x \\  &\cong \HolS(U_x \times C\str[L]_x)_x \\
         &\simeq \HolS(U_x \times C\str[L]_x) \\
         &\simeq U_x \times (C\str[L]_x)_{>p} \\
         &\simeq (C\str[L]_x)_{>p} \\
         & \cong (0,1] \times \str[L]_x 
         \\ &\simeq \str[L]_x.
    \end{align*}
    It follows that the local homotopy link can be homotopically identified with what is often referred to as the (local) link of a topological pseudomanifold (see, for example \cite{banagl2007topological}).
\end{remark}

%% file: 6Appendix.tex
\section{Refined stratified spaces and Nand-Lal's homotopy theory}\label{subsec:ref_and_nand_theory}
    In this section, we relate the homotopy theory $\AltStratCR$ to the work of \cite{nand2019simplicial}:
    \begin{recollection}[\cite{nand2019simplicial}]\label{recol:nand_setup}
    To avoid dealing with empty strata, \cite{nand2019simplicial} introduced the notion of a \define{surjectively stratified space}, defined as follows: Denote by $\StratN_{s}$ the full subcategory of $\StratN$ given by such stratified spaces $\tstr$, for which the stratification $\ststr \colon \utstr \to \ptstr$ is surjective. In other words, such stratified spaces whose stratification poset contains no redundant elements (but possibly redundant relations). Such a stratified space is called \define{surjectively stratified}.
    In \cite{nand2019simplicial}, Nand-Lal constructs a homotopy theory for $\StratN_{s}$ via transfer along the composition 
    \[
    \StratN_{s} \hookrightarrow \StratN \xrightarrow{\SingS} \sStratN \xrightarrow{\forget} \sSetN,
    \]
    using the Joyal model structure on $\sSetN$. In other words, a weak equivalence in the resulting homotopy theory is a map of stratified spaces $f \colon \tstr \to \tstr[T]$, for which the underlying simplicial map of the stratified simplicial map $\SingS(f)$ is a categorical equivalence. It follows that this class of weak equivalences is precisely the class of categorical equivalences with source and target in $\StratN_s$.
    We denote by $\RelStratsCR$ the relative category given by $\StratN_s$, together with such weak equivalences.
    \end{recollection}
    Let us now relate $\RelStratsCR$ to $\RelStratCR$. 
    \begin{proposition}\label{prop:eq_of_sur_strat}
        The inclusion $\StratN_s \hookrightarrow \StratN$ induces a homotopy equivalence of relative categories \[ \RelStratsCR \xrightarrow{\simeq} \RelStratCR ,\]
        and hence an equivalence of the corresponding $\infty$-categories. 
    \end{proposition}
\begin{proof}
         The functor $(-)^{\ared} \colon \StratN \to \StratN$ has image in $\StratN_s$.
         By \cref{prop:ref_com_sings}, there is a natural isomorphism $\SingS \circ (-)^{\ared} \cong (-)^{\ared} \circ \SingS$, which shows that $(-)^{\ared}$ preserves refined categorical equivalences.
         Furthermore, since (again by \cref{prop:ref_com_sings}) the refinement map is a categorical equivalence, it follows that $(-)^{\ared}$ defines a homotopy inverse to the inclusion of relative categories $ \RelStratsCR 
   \hookrightarrow \RelStratCR$. One homotopy is given by the refinement map and one by the identity.
\end{proof}
\section{An elementary extension lemma}
In this section, we give a proof of \cref{lem:extension_lemma}.
Before we give a proof, let us introduce some notation and quickly illustrate where metrizability comes into play. 
It is not hard to see that there is a continuous map $i\colon D^{n+1} \times [0,1] \to D^{n+1} \times [0,1) \cup S^{n} \times [0,1]$, which is the identity on $D^{n+1} \times \{ 0 \} \cup S^{n} \times [0,1]$ and maps $\mathring D^{n+1} \times [0,1]$ into $\mathring D^{n+1} \times [0,1)$. Hence, the obvious thing to do would be to simply glue $f$ and $f'$ along $D^{n+1} \times \{ 0 \} \cup S^{n} \times [0,1)$ to a map $\hat f \colon D^{n+1} \times [0,1) \cup S^{n} \times [0,1] \to X$ and then set $\tilde f = \hat f \circ i$. The issue with this approach is that the map $\hat f$ obtained in this way will generally not be continuous. Indeed, it is obtained by gluing two maps defined respectively on an open and a closed subset of $D^{n+1} \times [0,1) \cup S^{n} \times [0,1]$. To illustrate this point a little better,
consider the homeomorphism 
\begin{align*}
    D^{n+1} & \to \faktor{S^{n} \times [0,1]}{S^n \times \{1\}} \\
    y &\mapsto [\frac{y}{\norm{y}}, 1- \norm{y}]
\end{align*}
mapping $0$ to the point given by $S^{n} \times \{1\}$. We obtain a change of coordinates $y \hateq {} [x,s]$. 
By setting 
\[
f'_{x,t} (s) = f'([x,s],t)
\] we can interpret the data of a map $f' \colon D^{n+1} \times [0,1] \to X$ as a continuous family of paths $f'_{x,t}\colon [0,1] \to X$, indexed over $(x,t) \in S^n \times [0,1)$, which fulfill \[
\gamma_{x,t}(1) = \gamma_{x',t}(1),\] for all $x,x' \in S^n, t\in [0,1]$.
If we want $\hat f = f \cup f'$ to be continuous, then we precisely need convergence 
\[
f'_{x_n,t_n} (s_n) \to f([x,0],1)
\]
for sequences $([x_n,s_n],t_n) \to ([x,0],1)$. Let $d \colon X \times X \to [0,\infty)$ denote a metric that induces the topology on $X$ and suppose that there is a uniform bound
\[
d(f'_{x,t} (s), f'_{x,t} (0) ) = d(f'_{x,t} (s), f([x,0],t) )  \leq \varphi(s)
\]
by some continuous function $\varphi \colon [0,1] \to [0,1]$ with $\varphi(0) = 0$, at least for $(s,t) \in [0,\frac{1}{2}] \times [\frac{1}{2},1]$. It follows from the triangle inequality that
\[
d(f'_{x_n,t_n} (s_n), f([x,0],1) ) \leq \varphi(s_n) + d( f([x_n,0],t_n), f([x,0],1)),
\]
for $(s_n, t_n)$ close to $(0,1)$.
In particular, it then follows from the continuity of $f$ that this expression converges to $0$.
If the speed at which the paths $f'_{x_n,t_n}$ leave $f'_{x_n,t_n}(0)$ becomes arbitrarily large as $(x_n,t_n) \to (x,1)$, then such a bound $\varphi$ may generally not exist. We may, however, use the metrizability of $X$ to reparametrize the paths $f'_{x,t}$ in a way to enforce such global bounds.
We proceed as follows.
\begin{proof}[Proof of \cref{lem:extension_lemma}]\label{appendix:elementary_lem}
    Using the notation above, $i$ is given by
    \begin{align*}
    i: D^{n+1} \times [0,1] &\to   D^{n+1} \times [0,1) \cup S^{n} \times [0,1]\\
    ([x,s],t) &\mapsto ([x,s], s\frac{t}{2} + (1-s)t).
    \end{align*}
    It remains to show that we may without loss of generality assume that $f'$ admits a bounding function $\varphi \colon [0,1] \to [0,1]$ as above. 
    For two topological spaces $T,T'$, denote by $T^{T'}$ the mapping space equipped with the compact open topology. If $T'$ is locally compact Hausdorff, then this construction defines the right adjoint to the functor $-\times T'$ (\cite{KelleyMappingSpaces}). Furthermore, if $T$ is metrizable, then the topology on $T^{T'}$ is easily seen to be the topology of uniform convergence on every compactum. Finally, denote by $\mathrm{Aut}(\mathbb R_{\geq 0})$ the subspace of self-homeomorphisms of $\mathbb R_{\geq 0}$.
    Consider the following three maps.
    \begin{align*}
     M^{[0,1]} &\to\mathrm{Aut}(\mathbb R_{\geq 0})\\
     \gamma &\mapsto  \sigma_{\gamma} :=\{ s \mapsto \sup_{t \leq s,1}\dmet{\gamma(t)}{\gamma(0)} + s \};
    \end{align*}
    \begin{align*}
    \mathrm{Aut}(\mathbb R_{\geq 0}) &\to\mathrm{Aut}(\mathbb R_{\geq 0}) \\
    \sigma &\mapsto \sigma^{-1};
    \end{align*}
    \begin{align*}
    \mathrm{Aut}(\mathbb R_{\geq 0}) &\to [0,1]^{[0,1]} \\
    \sigma &\mapsto \{ s \mapsto \min\{ \sigma(s),1 \} \}.
    \end{align*}
    It is not hard to see, using the topology of uniform convergence, that the first and last of these maps are continuous. 
    That inversion is continuous follows from \cite[Thm. 4]{homeogroup}. Now, let 
\begin{align*}
     \rho: M^{[0,1]} \to [0,1]^{[0,1]}.
\end{align*}
    be the composition of these three maps. Then $\rho$ has the following properties:
  
    \begin{RhoProps}
        \item\label{property:rho_extension_lemma_1} $\rho(\gamma)(0) = 0$ if and only if $s = 0$;
        \item\label{property:rho_extension_lemma_2} $d\big (\gamma(\rho(\gamma)(s)), \gamma(0) \big ) \leq s ,$
    \end{RhoProps}
    for all $s \in [0,1], \gamma \in X^{[0,1]}$. The first property is immediate by construction of $\rho$. The second inequality is obtained from 
    \begin{align*}
        d\big (\gamma(\rho(\gamma)(s)), \gamma(0) \big ) & \leq  \sup_{t \leq \rho(\gamma)(s), 1} d\big (\gamma(t), \gamma(0) \big ) + \rho(\gamma)(s) = \sigma_\gamma ( \rho(\gamma)(s))  \\
                                                       & 
                                                         \leq \sigma_\gamma (  \sigma_\gamma^{-1}(s)) \\
                                                         &= s.
    \end{align*}
    Next, denote by $\phi: [0,1]^2 \to [0,1]$ a function fulfilling
    %%%%
    
    %%
    \begin{phiProps}
    \item\label{property:phi_extension_lemma_1} $\phi(s,t) =  0$, $(s,t) \in [0,\frac{1}{2}] \times [\frac{1}{2},{1}]$ ;
    \item\label{property:phi_extension_lemma_2} $\phi(s,t) = 1$, $(s,t) \in [0,1] \times \{0\} \cup \{1\} \times [0,1]$, 
    \end{phiProps}
    and, define \begin{align*}
    \Phi: D^{n+1} \times [0,1) &\to D^{n+1} \times [0,1)  \\
    ([x,s],t) &\mapsto ([x, (1-\phi(s,t))\rho(f'_{x,t},s) + \phi(s,t)s], t).
\end{align*}
Note first that $\Phi$ is well defined, that is, its value is independent of $x \in S^n$ when $s = 1$, since then also $ (1-\phi(s,t))\rho(f'_{x,t},s) + \phi(s,t)s = 1$. 
Furthermore, $\Phi$, has the following properties:
%%%

    %%
\begin{PhiProps}
    \item \label{property:Phi_extension_lemma_1} $\Phi([x,s],t) = ([x,s],t)$, for $([x,s],t) \in D^{n+1} \times \{ 0 \} \cup S^{n} \times [0,1)$;
    \item \label{property:Phi_extension_lemma_2} $\Phi ( \mathring D^{n+1} \times [0,1) )  \subset  \mathring D^{n+1} \times [0,1) $;
    \item \label{property:Phi_extension_lemma_3} $\Phi([x,s],t) = \rho(f'_{x,t},s)$, for $(s,t) \in [0,\frac{1}{2}] \times [\frac{1}{2},{1}]$.
\end{PhiProps}
\cref{property:Phi_extension_lemma_1} follows from \cref{property:phi_extension_lemma_1,property:rho_extension_lemma_1}. \cref{property:Phi_extension_lemma_2} follows from \cref{property:rho_extension_lemma_1}, and finally \cref{property:Phi_extension_lemma_3} follows from \cref{property:phi_extension_lemma_2}. We may then replace $f'$ by $f'' =  f' \circ \Phi$, obtaining
\begin{fProps}
    \item \label{property:f_extension_lemma_1} $f''([x,s],t)= f'([x,s],t) = f([x,s],t)$, for $([x,s],t) \in D^{n+1} \times \{ 0 \} \cup S^{n} \times [0,1)$;
    \item \label{property:f_extension_lemma_2} $f''( \mathring D^{n+1} \times [0,1) ) \subset f'(\mathring D^{n+1} \times [0,1))$;
    \item \label{property:f_extension_lemma_3} $d(f'_{x,t} (s), f'_{x,t} (0) ) \leq s$, for $(s,t) \in [0,\frac{1}{2}] \times [\frac{1}{2},{1}]$.
\end{fProps}
\cref{property:f_extension_lemma_1,property:f_extension_lemma_2} are immediate by \cref{property:Phi_extension_lemma_1,property:Phi_extension_lemma_2} respectively, while \cref{property:f_extension_lemma_3} follows from \cref{property:Phi_extension_lemma_3} and \cref{property:rho_extension_lemma_2}.
\end{proof}
\section{Remaining part of the proof of the nonexistence proposition}
In \cref{prop:counter_example_mod} (using the notation there) we claimed that a path $\gamma$ from $a$ to $b$ that (described as starting from $b$) ascends monotonously in height and passes to the left of $b_1$, to the right of $c_2$, to the left of $b_3$ etc., as illustrated, cannot lie in the image of $r_*$ (even up to stratified homotopy). Here, we provide a rigorous proof of this statement.
\begin{proof}\label{appendix:proof_of_nonexistence_of_oath}
 This insight can be formalized as follows: For $n \in \mathbb N$, denote by $\tstr_{b_n}$ the stratified spaces obtained from $\tstr$ by taking only the point $b_n$ as the $p$-stratum.
    The identity at the space level $1_{\utstr}$ does not induce stratified maps $\tstr, \tstr[Y] \to  \tstr_{b_n}$. However, it nevertheless induces (non-)stratified maps of stratified singular simplicial sets. Furthermore, by the same argument as above, $\SingS(\tstr_{b_n})$ are quasi-categories.
    In addition to this, on $\tstr_{b_n}$, $r$ is stratified homotopic to the identity relative to $a$ and $b$.
    It follows that we obtain a commutative diagram  
        \begin{diagram}
            \ho  \SingS(\tstr) (a,b) \arrow[rr, "r_*"]  \arrow[rd, "\tau_{c_n}"] & & \ho  \SingS(\tstr[Y]) (a,b) \arrow[ld, "{\sigma_{b_n}}"] \\ 
                & \ho \SingS (\tstr_{b_n})(a,b) &
        \end{diagram}
        Now, $\SingS (\tstr_{b_n})(a,b)$ is simply the set of homotopy classes of paths in $X \setminus \{b_n\}$, which is homotopy equivalent to $S^1$. Under post-composition with the straight line path from $b$ to $a$, we may thus identify \[
        \ho \SingS (\tstr_{b_n})(a,b) \cong \ho \SingS (\tstr_{b_n})(a,a) \cong \pi_1(S^1) \cong \mathbb Z .\]
        We may proceed mutatis mutandis replacing $b_n$ with $c_n$. Then, under this identification with the integers, the sequences $\sigma_{b_n}([\gamma]), \sigma_{c_n}([\gamma]) \in \mathbb Z$ are (modulo signs, originating from choices of orientation) given by:
        \begin{align*}
            &(1, 0, 1, 0, 1, \cdots) \\
            &(0, 1, 0 , 1 ,0, \cdots). 
        \end{align*}
        However, for every morphism $g \in \ho  \SingS(\tstr) (a,b)$ at least one of the sequences $\tau_{b_n}(g)$ or $ \tau_{c_n}(g)$ has only finitely many non-zero values.
        Indeed, assume that $g$ is such that both $\tau_{b_n}(g)$ and $\tau_{c_n}(g)$ have infinitely many non-zero values. Every such $g$ can be represented by a path $\gamma'$ which stays in $\tstr_p$ for $[0,\frac{1}{2}]$ and then exits into $\tstr_q$.
        For $\tau_{b_n}(g)$ ($\tau_{c_n}(g)$) to be non-zero, $\gamma'$ must intersect the line segment connecting $a$ and $b_n$ ($c_n$), as otherwise $\gamma'$ maps into a contractible subspace of $\tstr_{b_n}$ ($\tstr_{c_n})$. Since these intersection points necessarily converge to $a$ as $n \to \infty$, it follows that $a \in \gamma[\frac{1}{2},1]$ and therefore $\gamma(\frac{1}{2}) = a$. 
        It follows that $\gamma'(t)$ lies strictly above $a'$ (in direction of the $y$-axis), for all $ t> \frac{1}{2}$ greater than some $t_h \in \frac{1}{2}$. Furthermore, it follows from the stratification of $\gamma$ and the assumption on infinitely many intersection points with the line segment $b_n$ and $a$ above, that we find $t_b, t_c \in  (\frac{1}{2},t_h)$, such that $\gamma'(t_b)$ lies on a line segment between $a$ and $b_n$ for some $n >0$, and $\gamma'(t_c)$ lies on a line segment between $a$ and $c_m$, for some $m >0$. To see this (in the case of $(b_n)$), consider a sequence $t_i$ such that each $\gamma'(t_i)$ lies on the line segment between $a$ and $b_i$. By compactness of $[\frac{1}{2},1]$, we may assume that $t_i$ converges to some $t' \in [\frac{1}{2},1]$. Since the intersection points of the line segments from $a$ to $b_i$ converge to $a$, it follows that $\gamma'(t') = a$. But $\frac{1}{2}$ is the only value in $[\frac{1}{2},1]$ with $\gamma'(t') \in \tstr_p$, hence $t' = \frac{1}{2}$, and it follows that $t_i$ converges to $\frac{1}{2}$.
        Finally, the set of points strictly above $a'$ in $\tstr[X]_q$ is disconnected, with the line segments between $a$ and $b_m$ and $a$ and $c_n$ lying in different components, for all $m,n$, in contradiction to the connectedness of $\gamma'(\frac{1}{2},t_h)$.   
\end{proof}
\section{On the proof of \texorpdfstring{\cite[Thm. 0.1.1/3.2.4]{haine2018homotopy}}{Haine's homotopy hypothesis}}\label{sec:haines_proof_of_strat_ho_hy}
In this section, we take a detailed look at the proof of \cite[Thm. 0.1.1/3.2.4]{haine2018homotopy}, pointing out a gap in the latter which stems from different uses of the terminology \textit{Segal space}. 
\begin{remark}\label{rem:haines_proof}On the proof of \cite[Thm. 0.1.1/3.2.4]{haine2018homotopy}:
Let $\TopPN^{\textnormal{ex}}$ denote the full subcategory of $\TopPN$ given by the fibrant stratified spaces in the categorical model structure and denote by $W_{\textnormal{ex}}$ the class of categorical equivalences between the latter. 
It was already asserted in \cite[Thm. 0.1.1/3.2.4]{haine2018homotopy} that $\SingS$ induces an equivalence of quasi-categories \newcommand{\ex}{\textnormal{ex}}
    \[
   \TopPN^{\textnormal{\ex}}[W^{-1}_{\ex}] \xrightarrow{\simeq} \AbStrP.
    \]
While this result is a consequence of the existence of the categorical semi-model structure on $\TopPN$ together with \cref{thm:overview_over_all_hypothesis}, \cite{haine2018homotopy} suggested a proof that did not assume the existence of such a structure. We want to point out two gaps within this proof, which seem to require a notion of fibrant replacement in $\TopPCN$ to be closed. Let us first sketch the proof in \cite{haine2018homotopy}. Some of the equivalences of $\infty$-categories in \cite{haine2018homotopy} are not made explicit. We will choose explicit models which, to the best of our knowledge, present the intended functors of $\infty$-categories.
\begin{enumerate}
    \item Denote by $\sSetPN^{\ex}$ the full subcategory of $\sSetPN$ given by the (bi)fibrant objects in $\sSetPCN$ (i.e. quasi-categories with a conservative functor into $\pos$) and by $H$ the class of equivalences of quasi-categories over $\pos$, between the latter. The quasi-category $\sSetPN^{\ex}[H^{-1}]$ is one possible model for the $\infty$-category of abstract stratified homotopy types over $\pos$, $ \AbStrP$. One may then consider a diagram of $1$-categories (commutative up to natural isomorphism)
    \begin{diagram}{}
       \TopPN^{\ex} \arrow[d, hook] \arrow[r] & \sSetPN^{\ex} \arrow[d, "{\HolIPS[]}"]\\
       \TopPN \arrow[r, "{\HolIP[]}"] & \FunC ( \sd(\pos)^{\op}, \sSetN) \spaceperiod
    \end{diagram}
    Now, let $W_{\dec}$ be the class of morphisms in $\FunC ( \sd(\pos)^{\op}, \sSetN)$ that are weak equivalences in the model structure presenting d\'ecollages over $\pos$. 
    The class of categorical equivalences in $\TopPN$, denoted $W_{\mathfrak{c}}$ is precisely the inverse image of $W_{\dec}$ under $\HolIP[]$. Furthermore, \cite{haine2018homotopy} proves that $H$ is precisely the inverse image of $W_{\dec}$ under the right vertical, and that $W_{\ex}$ is precisely the inverse image of $H$ under $\SingS$. Hence, there is an induced diagram of quasi-categories
    \begin{diagram}\label{diag:Haines_proof}{}
       \TopPN^{\ex}[W_{\ex}^{-1}] \arrow[d, hook] \arrow[r, "\SingS"] & \sSetPN^{\ex}[H^{-1}] \arrow[d, "{\HolIPS[]}"] \simeq \AbStrP \\
      \AltTopPC = \TopPN{}[W_{\mathfrak c}^{-1}] \arrow[r, "{\HolIP[]}"] & \FunC ( \sd(\pos)^{\op}, \sSetN)[W_{\dec}^{-1}] \spacecomma
    \end{diagram}
    commutative up to natural isomorphism. 
    It follows from \cite[Thm. 3]{douteauEnTop}, that the lower horizontal is an equivalence. Furthermore, it follows either from \cite[\cref{comb:prop:equ_decol_haine}]{ComModelWa} or from \cite[Thm 2.7.4]{Exodromy} that the right vertical is an equivalence. This already proves one of the central results of \cite{haine2018homotopy}, namely the existence of an equivalence $\AltTopPC \simeq \AbStrP$. In addition to this result, \cite{haine2018homotopy} aims to show that the upper horizontal in \cref{diag:Haines_proof} is also an equivalence. This follows from commutativity of \cref{diag:Haines_proof}, if one can show the following two additional claims:
    \item The left vertical in \cref{diag:Haines_proof} 
    \[
     \TopPN^{\ex}[W_{\ex}^{-1}] \to \TopPN{}[W_{\mathfrak c}^{-1}] =  \AltTopPC
    \]
    is fully faithful: 
    This is not commented on further in \cite{haine2018homotopy}. 
    Note that there is generally no reason to assume that in a category with weak equivalences $(\cat[C],W)$ with full subcategory $\cat[A]$ the natural functor $\cat[A][\cat[A] \cap W^{-1}] \to \cat[C][W^{-1}]$ is an equivalence of quasi-categories.
    Consider, for example, the flattening of a quasi-category $\mathcal{C} \in \sSet$, given by the $1$-category $(\Delta_{/\mathcal C})^{\op}$.
    The flattening of $\mathcal{C}$ has a full discrete subcategory $\mathcal C_0$ given by the objects of $\mathcal{C}$. 
    If we denote by $W$ the class of morphisms in $(\Delta_{/\mathcal C})^{\op}$ that map $0$ to $0$, there is a natural equivalence of quasi-categories, $(\Delta_{/\mathcal C})^{\op} [W^{-1}] \simeq \mathcal{C}$, which fixes the objects in $\mathcal{C}_0$ (see, for example, \cite[Thm. 3.3.8]{Land}). However, $\mathcal{C}_0[W \cap \mathcal{C}_0^{-1}]$ remains a discrete category.
    A classical condition to even ensure equivalence between $\cat[C][W^{-1}]$ and  $\cat[A][\cat[A] \cap W^{-1}]$ is the existence of an endofunctor $F \colon \cat[C] \to \cat[C]$, with image in $\cat[A]$, together with a zig-zag of natural transformations $1 \leftrightarrow F$, given by morphisms in $W$.
    In the situation of \cref{diag:Haines_proof} this amounts to constructing a fibrant replacement in $\TopPCN$, which we have done in this work, but whose existence was not yet known at the time of writing of \cite{haine2018homotopy}.
    \item The diagonal in \cref{diag:Haines_proof} 
    \[
    \TopPN^{\ex}[W_{\ex}^{-1}] \to \FunC ( \sd(\pos)^{\op}, \sSetN)[W_{\dec}^{-1}] 
    \]
    (and hence also the left vertical) is essentially surjective: Again, this claim is close to requiring fibrant replacements in $\TopPCN$. 
    Recall that $\AltDiag$ denotes the localization of $ \FunC( \sd(\pos)^{\op}, \sSetN) $ at pointwise-weak homotopy equivalences of diagrams.
    It follows from the equivalence $\AltTopPD \simeq \AltDiag$, that every d\'ecollage in $\AltDiag$ lies in the image of 
    \[ \TopPN \xrightarrow{\HolIP[]} \FunC( \sd \pos, \sSetN) \to  \AltDiag. \]
    In other words, if we denote by $\TopPN^{\dec}$ the full subcategory of stratified spaces $\str$ for which $\HolIP[](\str)$ is a d\'ecollage, then the restricted functor 
    \[
    \TopPN^{\dec} \to \FunC ( \sd(\pos)^{\op}, \sSetN)[W_{\dec}^{-1}]\]
    is essentially surjective. \cite{haine2018homotopy} then claims that $\TopPN^{\dec} = \TopPN^{\ex}$, from which essential surjectivity of the diagonal would follow. 
    With the definition of d\'ecollages we used here and which is also used in \cite{haine2018homotopy} this is incorrect. 
    The confusion arises from two possible definitions of d\'ecollages and complete Segal spaces: One is expressed purely intrinsically to the $\infty$-categories $\FunC ( \sd(\pos)^{\op}, \iGrpd)$ and $\FunC ( {\Delta}^{\op}, \iGrpd)$ and one expressed in terms of injective model structures on the $1$-categories $\FunC ( \sd(\pos)^{\op}, \sSetN)$ and $\FunC ( {\Delta}^{\op}, \sSetN)$. 
    Here, we have defined d\'ecollages intrinsically to $\FunC ( {\sd(\pos)}^{\op}, \iGrpd)$ and complete Segal spaces in \cite{haine2018homotopy} are also defined intrinsically to $\FunC ( {\Delta}^{\op}, \iGrpd)$. 
    Alternatively, one can define d\'ecollages as fibrant objects in the model category $\FunC ( {\sd(\pos)}^{\op}, \sSetN)^{\dec}$ (defined in \cite[Subsec. \ref{comb:subsec:decollages}]{ComModelWa}, as a left Bousfield localization of the injective model structure) and, as it is classically done, complete Segal spaces as (certain) fibrant objects in the injective model structure on $\FunC ( \Delta^\op, \sSetN)$.
    Let us call the former definition $\infty$-categorical and the latter $1$-categorical.
    \cite{haine2018homotopy} makes use of a result of Joyal and Tierney concerning the relationship of quasi-categories and $1$-categorical Segal spaces \cite[Cor. 3.6]{JoyalQCvsSS}, from which it follows that $\TopPN^{\ex}$ consists exactly of such stratified spaces for which $\HolIP[](\tstr)$ is a $1$-categorical d\'ecollage. Due to the overlap in terminology, \cite{haine2018homotopy} then reasons that $\TopPN^{\dec} = \TopPN^{\ex}$. This is false. Indeed, every stratified space with two strata, $\tstr$ over $\pos = \{ p <q \}$, has the property that $\HolIP[](\tstr)$ is an $\infty$-categorical d\'ecollage (there are no conditions to be verified). But for $\HolIP[](\tstr)$ to be injectively fibrant the map $\HolIP(\tstr) \to \tstr_p$, for $\I = [p < q]$,  needs to be a Serre fibration. For example, $\sReal{\Lambda^{\J}_0}$ with $\J = [p \leq p <q]$ fails to have this property.  
\end{enumerate}
\end{remark}